\documentclass[11pt]{article}
\usepackage{amsthm,amsmath,latexsym,amssymb}

\newcommand{\bP}{{\rm |\kern-.15em P}}
\newcommand{\Q}{\kern.3em\rule{.07em}{.65em}\kern-.3em{\rm Q}}
\newcommand{\R}{{\rm I\kern-.15em R}}
\newcommand{\D}{{\rm |\kern-.15em D}}
\newcommand{\h}{{\rm |\kern-.15em H}}
\newcommand{\C}{\kern.3em\rule{.07em}{.65em}\kern-.3em{\rm C}}
\newcommand{\T}{{\rm T\kern-.35em T}}

\theoremstyle{plain}
\newtheorem{theorem}{Theorem}[section]
\newtheorem{lemma}[theorem]{Lemma}

\newtheorem{proposition}[theorem]{Proposition}
\newtheorem{corollary}[theorem]{Corollary}

\theoremstyle{definition}
\newtheorem{definition}[theorem]{Definition}

\theoremstyle{remark}
\newtheorem{remark}[theorem]{Remark}

\newcommand\blfootnote[1]{%
  \begingroup
  \renewcommand\thefootnote{}\footnote{#1}%
  \addtocounter{footnote}{-1}%
  \endgroup
}

\begin{document}
\title{On the automorphic group of an entire function}
\author{Ronen Peretz}
 
\maketitle

\begin{abstract}
This paper develops further the theory of the automorphic group of non-constant entire functions. This theory
has already a long history that essentially started with two remarkable papers of Tatsujir\^o Shimizu
that were published in 1931. The elements $\phi(z)$ of the group are defined by the automorphic equation
$f(\phi(z))=f(z)$, were $f(z)$ is entire. Tatsujir\^o Shimizu also refers to the functions of this group
as those functions that are determined by $f^{-1}\circ f$. He proved many remarkable properties of those
automorphic functions. He indicated how they induce a beautiful geometric structure on the complex plane.
Those structures were termed by Tatsujir\^o Shimizu, the system of normal polygonal domains, and the more
refined system of the fundamental domains of $f(z)$. The last system if exists tiles up the complex
plane with remarkable geometric tiles that are conformally mapped to one another by the automorphic
functions. In the Ph.D thesis of the author, those tiles were also called the system of the maximal domains
of $f(z)$. One can not avoid noticing the many similarities between this automorphic group and its accompanying
geometric structures and analytic properties, and the more tame discrete groups that appear in the theory
of hyperbolic geometry and also the arithmetic groups in number theory. This paper pursues further the
theory initiated by Tatsujir\^o Shimizu, towards understanding global properties of the automorphic group, rather
than just understanding the properties of the individual automorphic functions. We hope to be able in sequel
papers to generalize arithmetic and analytic tools such as the Selberg trace formula, to this new setting.
\end{abstract}

\tableofcontents
\newpage

\section{Some background and the contribution of Tatsujir\^o Shimizu}\label{sec1}

\blfootnote{\textup{2010} \textit{Mathematics Subject Classification}: \textup{30B40,30C15,30D05,30D20,30D30,30D99,
30F10,30F35,32D05,32D15}}
\blfootnote{\textit{Key Words and Phrases:} \textup{entire functions, integral functions, meromorphic functions, 
fundamental domains, automorphic functions of a meromorphic function, the automorphic group of a meromorphic function, 
Weierstrass Factorization, Gronwall-Hahn Factorization}}

In 1931 Tatsujir\^o Shimizu published two remarkable papers having the titles: On the Fundamental Domains and the
Groups for Meromorphic Functions. I and II. \cite{s1,s2}. Here are quotations from Shimizu's \cite{s1} that
describe few of the main notions of the theory: \\
"{\bf 1. (p. 179)} We call an open domain of the Riemann surface of the inverse function $z=f^{-1}(w)$ of an
integral or a meromorphic function $w=f(z)$ "a leaf" if it satisfies the following three conditions: \\
1) It covers almost all the whole $w$-plane without leaving any complementary domain. \\
2) It does not cover the $w$-plane more than once in any part. \\
3) Each part of the boundaries is common to certain domains of the surface which is exterior to the considered
domain. \\
{\bf 2. (p. 179)} To each leaf so defined on the Riemann surface of the meromorphic function $w=f(z)$ there 
corresponds an open domain on the $z$-plane which we call a "polygonal domain". In the polygonal domain of
a meromorphic function the function is mono-valued and it takes all values except a set of values not forming
a domain. We call such a leaf whose boundary consists of only accessible points from the inside of it "a leaf
with accessible boundary" and the corresponding polygonal domain "a polygonal domain with accessible
boundary". By mapping the sequence of leaves with accessible boundaries on the $z$-plane, the $z$-plane
of the meromorphic function $w=f(z)$ is divided into a system of polygonal domains for the function, whose
boundaries consist only of accessible points from the inside of each domain respectively. \\
{\bf 3. (p. 185)} I will here call that the $z$-plane is divided into "a system of normal polygonal 
domains", if the $z$-plane is divided into a system of polygonal domains whose boundaries consist of
only accessible points from the inside of them so that an infinite number of the boundaries of different
polygonal domains may not accumulate in the finite part of the $z$-plane.

Further we call that the $z$-plane is divided into "a system of fundamental domains" if the plane is
divided into a system of normal polygonal domains so that the boundary of each normal polygonal domain 
may be all transformed into all the boundary of another normal polygonal domain by the transformation 
of the group of the function which I shall consider in the section VII, that is, the transformation defined
by $f(z')=f(z)$, where $z$ runs over the boundary of some polygonal domain. \\
{\bf 4. (p. 185)} We show that: \\
{\bf Theorem 1.} {\it For any meromorphic function $f(z)$ we can divide the $z$-plane into a system
of normal polygonal domains, that is, we can divide the Riemann surface for the inverse function of
$f(z)$ into a system of leaves without leaving any elements (belonging to the surface) except
point sets containing no domain such that each leaf covers almost all the $w$-plane except point
sets not forming domains and the boundary of each leaf consists of only accessible points from the
inside and further, when all the leaves thus obtained are mapped conformally on the $z$-plane, there
exists no point set in the finite part of the $z$-plane which is a limiting set of an infinite number
of boundaries of the images of the leaves.}" \\
In this basic theorem of the theory, Shimizu demonstrates that any meromorphic function carries with it the
geometric structure of a system of normal polygonal domains. However as he later on proves there are entire
functions that have no system of fundamental polygonal domains. Gross constructed an entire function whose
set of all asymptotic values is the whole of $\mathbb{C}$. Shimizu proves that Gross' function has no system
of fundamental polygonal domains. Thus any meromorphic function induces those remarkable tilings of the
complex plane by systems of normal polygonal domains. But there are entire functions for which the boundaries
of the different tiles are mapped to one another (infinity included) by the automorphic functions, in a rather
complicated manner. \\
A large portion of Tatsujir\^o Shimizu's papers was dedicated to understand the analytic and the geometric 
properties of the individual elements of the group that are defined by 
\begin{equation}\label{eq1} f(z')=f(z). \end{equation}
In our paper we will call this defining equation, the automorphic equation of $f(z)$. 
\begin{remark}\label{rem1}
In our manuscript we will call this group of Shimizu, "the automorphic group of $f(z)$" and we
will use the notation ${\rm Aut}(f)$ to designate it. The binary operation is composition of mappings.
\end{remark}
\begin{remark}\label{rem2}
As mentioned in remark \ref{rem1} we will use the term "automorphic function of $f(z)$" instead of Shimizu's "fundamental function
with respect to $f(z)$".
\end{remark}
\noindent
For example, very simple such groups are ${\rm Aut}(z^n)=\{e^{2\pi ik/n}z|\,k=0,\ldots,n-1\}$
and ${\rm Aut}(e^z)=\{z+2\pi ik|\,k\in\mathbb{Z}\}$. Possible tilings of the complex plane that correspond
to these groups (and functions) are $\Omega_j(z^n)=\{z\in\mathbb{C}|\,2\pi ij<\arg z<2\pi i(j+1)\}$, $j=0,\ldots,n-1$ and
$\Omega_j(e^z)=\{z\in\mathbb{C}|\,2\pi j<\Im z<2\pi (j+1)\}$, $j\in\mathbb{Z}$ respectively. However, these
two examples are exceptional, having all of the automorphic functions entire. A remarkable property proved by
Shimizu asserts that the only possible automorphic functions which are entire, have the form $e^{i\theta}z+b$,
where $\theta\in 2\pi\cdot\mathbb{Q}$. More complicated entire functions do not qualify being automorphic.
In general those automorphic functions are multi-valued or "leaves" thereof with a complicated structure.
Further research on this topic was carried on, for example in \cite{r} and \cite{r1}. There are papers
that computed systems of fundamental domains (and their automorphic groups) for specific important functions,
in particular in number theory such as the Riemann Zeta function and the Gamma function. It is clear that a lot of further research
is needed in order to better understand the automorphic functions. In particular we clearly have to understand
more global properties of the groups ${\rm Aut}(f)$ and of their induced normal (or fundamental if exist) polygonal domains
$\{\Omega_j(f)\}_j$. For example it is clearly important to understand if tools parallel to Selberg Trace Formula could be extended
to the automorphic groups of entire functions.

Here is a brief summary of the results and the ideas in the paper. In section \ref{sec3} we use the Weierstrass
representation as (generically) an infinite product for $f(w)-f(z)$. Here $w\in\mathbb{C}$ is the variable while
the parameter $z$ lies in $\mathbb{C}-f^{-1}(f(0))$. We have:
$$
f(w)-f(z)=\exp{\left(g(w,z)\right)}\prod_{n=1}^{\infty}E\left(\frac{w}{\phi_{0n}(z)},\lambda_n\right)=
$$
$$
=\exp{\left(g(w,z)\right)}\prod_{n=1}^{\infty}\left(1-\frac{w}{\phi_{0n}(z)}\right)e^{Q_{\lambda_n}(w/\phi_{0n}(z))},
$$
where if $\lambda_n>0$, then:
$$
Q_{\lambda_n}\left(\frac{w}{\phi_{0n}(z)}\right)=\left(\frac{w}{\phi_{0n}(z)}\right)+\frac{1}{2}
\left(\frac{w}{\phi_{0n}(z)}\right)^2+\ldots +\frac{1}{\lambda_n}\left(\frac{w}{\phi_{0n}(z)}\right)^{\lambda_n},
$$
and $Q_0(w/\phi_{0n}(z))\equiv 0$. Weierstrass representation parameters are the function $g(w,z)$ which is entire
in $w$ and $z$-holomorphic off $f^{-1}(f(0))$, and the non-negative integers $\lambda_n$ that depend on $z$. Clearly
$f(w)-f(z)$ is $z$-${\rm Aut}(f)$ invariant. But this happens due to a complicated interaction of the infinite
product and the exponential $\exp{\left(g(w,z)\right)}$. In the case that $f$ has a finite order it follows that
the infinite product as well as the exponential part are separately $z$-${\rm Aut}(f)$ invariant. Thus in this
case the behavior of the Weierstrass representation is tame, for the group invariance is not requiring any interaction
between the two parts of the Weierstrass representation. This is proved in Proposition \ref{prop3}. We prove that
the description of the Weierstrass representation of $f(w)-f(z)$ can be refined in that the exponential part
has the form $exp(F(w,f(z)))$ were $F(w,t)$ is holomorphic in each variable separately. This depends among other things on
Lemma \ref{lem1}. This lemma also implies the cycle relation, Corollary \ref{cor3} and the chain relation, Corollary
\ref{cor4}. However the proof of Lemma \ref{lem1} follows by a result of Eremenko and Rubel, \cite{er}, which makes
a use of the monodromy principle. 

In section \ref{sec4} we indicate what conclusions can be reached when we have no monodromy. In particular the proof
of Corollary \ref{cor05} defines the mapping $T:\,{\rm Aut}(f)\rightarrow\mathbb{Z}$, which will later on be used in
section \ref{sec7} (for example in Corollary \ref{cor11}). Most of the results in this section deal with the arithmetic
of the compositions of automorphic functions. The mapping $T$ provides means to induce from any such a composition
an appropriate factorization of a natural number over $\mathbb{Z}^+$. Thus we can use the multiplicative theory of numbers
in order to deduce results on factorizations of an automorphic function into a composition of other automorphic functions.
Theorem \ref{thm0} gives the general picture by describing the finiteness of the decomposition of automorphic functions.

In section \ref{sec5} the cycle relation and the chain relation are discussed in the general case where no assumption on
the finiteness of the order is assumed. Corollary \ref{cor53} deals with the cycle relation while Corollary \ref{cor54}
and Corollary \ref{cor55} deal with the chain relation.

In section \ref{sec6} we compute the example of the exponential function. We use our method of computation to arrive at a
general result, Theorem \ref{thm2} which indicates how to construct the entire function $f(w)$ from its automorphic
group ${\rm Aut}(f)$. This construction depends on the assumption that has no justification at the moment, that we have
some summation method for the infinite series:
$$
\sum_{n=1}^{\infty}Q_{\lambda_n}\left(\frac{w}{\phi_{0n}(z)}\right).
$$
The right summation method for this infinite sum of polynomials in $w$ which are multi-valued functions of $z$ is
an open problem.

Section \ref{sec7} gives, among other things other types of reconstruction formulas both to $f(z)$ and to $f'(z)$ in terms
of approximating automorphic functions. These are the automorphic functions of the partial sums of the power series expansion
of $f(w)$. Proposition \ref{prop4} gives the formulas: $f(z)=f(w)-\lim_{n\rightarrow\infty}a_n\prod_{j=1}^n(w-\phi_{j}{(n)}(z))$, and
$f'(z)=\lim_{n\rightarrow\infty}a_n\prod_{j=1}^{n-1}(z-\phi_{j}^{(n)}(z))$. The remarkable thing here is that $a_n\rightarrow 0$
as a sequence of numbers while both products blow up but as sequences of functions, but just in the right pace so that the
limits converge and reconstruct the function and its derivative. As mentioned above this section gives also properties of $\ker(T)$,
where the mapping $T$ was defined in section \ref{sec4} within the proof of Corollary \ref{cor05}. For example Corollary \ref{cor11} indicates  relations between $\ker(T)$ and the automorphic group. In fact when $f$ has a finite order as an entire function then
$\ker(T)={\rm Aut}_z(g(w,z))$, and all the automorphic functions in ${\rm Aut}(f)-{\rm Aut}_z(g(w,z))$ have infinite order in the
sense of group members. We recall that ${\rm Aut}(f(z))\subseteq{\rm Aut}_z(\exp(g(w,z))$. So ${\rm Aut}(f(z))$ is "trapped"
between ${\rm Aut}_z(g(w,z))$ and ${\rm Aut}_z(\exp(g(w,z))$ and the automorphic functions of $f$ outside the smaller group
${\rm Aut}_z(g(w,z))$ all have infinite order as group elements of the automorphic group of $f$, ${\rm Aut}(f)$.

In section \ref{sec8} we show how the function $g(w,z)-g(0,z)$, where $g(w,z)$ is the function that participates in the Weierstrass
representation of $f(w)-f(z)$ is determined by negative moments of the automorphic functions. Theorem \ref{thm3} determines
$\frac{\partial^kg}{\partial w^k}(0,z)$ in terms of $\sum(\phi_{0n}(z))^{-k}$. In fact for $k=1,2,3,\ldots $ we have the identities:
$$
\frac{1}{k!}\frac{\partial^k g}{\partial w^k}(0,z)=-\frac{1}{k}\sum_{\left\{\begin{array}{l} n \\ \lambda_n\ge k\end{array}\right.}
\left(\frac{1}{\phi_{0n}(z)}\right)^k.
$$
The left hand side is the $k+1$'st Maclaurin coefficient in the expansion of $g(w,z)-g(0,z)$. The right hand
side is $(-1/k)$ multiplying the $(-k)$-moment of all the relevant automorphic functions of $f$. That explains the
title of this section. The whole argument is based on the assumption that we have some summation method for the infinite series:
$$
\sum_{n=1}^{\infty}Q_{\lambda_n}\left(\frac{w}{\phi_{0n}(z)}\right).
$$
This assumption was also needed for deducing an essential part of Theorem \ref{thm2} in section \ref{sec6}. As mentioned
in section \ref{sec6} this summation problem is an open problem.

In section \ref{sec9} an infinite product representation of $f'(w)$ in terms of automorphic functions is given.
Proposition \ref{prop5} shows that ${\rm Aut}(f)\subseteq{\rm Aut}(g)$ implies that $\exists\,G(w,z)$, entire in $w$ and
holomorphic in $z\in\mathbb{C}-g^{-1}(g(0))-f^{-1}(f(0))$ such that $g(w)-g(z)=(f(w)-f(z))\cdot G(w,z)$. In particular
$\exists\,H(z)$, an entire function such that $g'(z)=H(z)\cdot f'(z)$. Lemma \ref{lem3} points out to a relation between the
fixed points of the non-identity automorphic functions ($\phi_{0n}(w)=w$ for $\phi_{0n}\not\equiv {\rm id.}$), and the
zeros of $f'(z)$, i.e. $Z(f')$. Accordingly $Z(f')$ might contain on the top of these fixed-points also elements from the fiber
$f^{-1}(f(0))$. Theorem \ref{thm4} gives a formula for $Z(f')$ in terms of the fixed-point sets of the non-identity automorphic
functions of $f(w)$. Theorem \ref{thm5} uses the Laguerre Theorem on separation of zeros and the formula of Theorem \ref{thm4}
to show the reality and the separation property of ${\rm Fix}({\rm Aut}(f))$ by $Z(f)$. 

Section \ref{sec10} deals with entire functions of the form $f(z)=P(z)e^{g(z)}$ where $P(z)\in\mathbb{C}[z]$ and where $g\in E$.
Let $d:=\deg p>0$ and $Z(p)=\{\alpha_1,\ldots,\alpha_d\}\subseteq\mathbb{C}$. Then $\forall\,j=1,\ldots,d$, $\alpha_j$ is a common zero of almost all the reciprocals of the automorphic functions of $f(z)$.

The main issue in section \ref{sec11} are formulas for the derivatives 
of the automorphic functions. Theorem \ref{thm7} give a kind of integral formula for $\sum_{|\phi_{0n}(z)|<R}\phi_{0n}^{(k)}(z)$.
Proposition \ref{prop6} gives a kind of a partial fractions expansion in terms of $1/(w-\phi_{0n}(z))$ to the $w$-logarithmic
derivative of $f(w)-f(z)$. In Proposition \ref{prop7} a parallel expansion is given for the $z$-logarithmic derivative
of $f(w)-f(z)$ (recall that $w$ is the function's variable while $z\not\in f^{-1}(f(0))$ is a parameter).

In section \ref{sec12} we apply the Jensen Theorem to compute the absolute value of products of automorphic functions in terms
of an integral of:
$$
\log|f(|\phi_{0n}(z)|e^{i\theta})-f(z)|d\theta
$$
These products are further discussed in section \ref{sec14}.

Section \ref{sec15} deals with order and type estimates of $f(w)$ in terms of the convergence exponent of the automorphic
group ${\rm Aut}(f)$. See Theorem \ref{thm8}. Some of the results are related to low order (less than $1$) (Theorem \ref{thm9}).
There are in this section also density estimates for ${\rm Aut}(f)$ for an entire $f(w)$ of a finite order. Theorem \ref{thm12}
and Theorem \ref{thm13} deal with entire functions of a finite and non-integral order and tie the convergence exponent of the
${\rm Aut}(f)$-orbits to this order.

Section \ref{sec20} is preparing for a future research on extending scattering theory, Selberg Trace formula etc... to
the setting of the discrete groups ${\rm Aut}(f)$. Theorem \ref{thm14} suggests what should be some of the counterparts of the
classical theory in the setting of ${\rm Aut}(f)$. This is far from being final and conclusive!

In the short section \ref{sec27} we bring the basics of the notion of local groups. This notion is clearly relevant to the
theory of the automorphic group of an entire function. The material is mostly taken from Terrence Tao's book \cite{tao}.

In section \ref{sec28} we prove the remarkable identities
$$
\lim_{j\rightarrow\infty}\sum_{\phi_{0n}(w)|<R_j}\phi_{0n}^{(k)}(w)\equiv 0,\,\,\,\,\,\forall\,w\in\mathbb{C},
$$
for certain sequences $0<R_1<R_2<\ldots<R_n<\ldots (R_n\rightarrow\infty)$. Each of these sequences fit simultaneously all
the values of $k\in\mathbb{Z}^+$. This is done for low order functions ($0<\rho<\frac{1}{2}$). The main tools used in the
proof are Wiman's-$\cos\pi\rho$ Theorem and our integral formulas for $\sum_{|\phi_{0n}(z)|<R}\phi_{0n}^{(k)}(z)$ in Theorem \ref{thm7}.
See Theorem \ref{thm15}. Few examples are elaborated to demonstrate the sharpness of our results here.

In section \ref{sec29} we give in Theorem \ref{thm16} some density estimates on  $\{|\phi_{0n}(z)|\}_n$ for functions of
a low order ($0<\rho<\frac{1}{2}$). Again Wiman's-$\cos\pi\rho$ Theorem is a main tool in our proof.

Vieta type formulas for ${\rm Aut}(f)$, $0\le\rho<1$, are given in section \ref{sec30}.
Theorem \ref{thm17} gives a formula for the Maclaurin coefficients of $f(w)$ in terms of $f(0)-f(z)$ and ${\rm Aut}(f)$.

A reasonable approach to try and extend the classical scattering theory results and the Selberg Trace formula to ${\rm Aut}(f)$,
is to naturally embed the automorphic group in a larger group, the way ${\rm SL}_n(\mathbb{Z})$ is embedded in 
${\rm GL}_n(\mathbb{R})$. In section \ref{sec31} and in section \ref{sec32} we try such an approach. Our initial embedding uses an ascending sequence of automorphic groups and a major problem is to try and understand what is the direct limit that is constructed. A typical example originates in the Tuen Wai NG construction of entire functions which have factorizations of unlimited number of prime factors, \cite{tuen}.
Theorem \ref{thm18} gives the structure of the direct limit group of the ascending automorphic groups. That is done in a certain important
case where a Tuen Wai NG function underlies the direct limit. Non-trivial consequences follow in Theorem \ref{thm19} and in
Theorem \ref{thm20}.

Section \ref{sec33} gives continuity relations between the group ${\rm Aut}(f)$ and
the sequence of groups ${\rm Aut}(f_n)$, where $f_n\rightarrow f$ uniformly on compact subsets of $\mathbb{C}$. The proof of that
theorem (Theorem \ref{thm24}) is tricky. It uses the elementary Newton's identities for moments and symmetric functions
of finite sets, and it uses one of the partial fractions expansions we found before for the $w$-logarithmic
derivative $f'(w)/(f(w)-f(z))$, in section \ref{sec11}, Proposition \ref{prop6}.

Those results are used in section \ref{sec34} to prove some results on the amenability
of ${\rm Aut}(f)$. In Theorem \ref{thm26} the assumptions are analytical while in Theorem \ref{thm28} the assumptions are
geometrical, i.e. they use the generations counting functions for systems of fundamental domains of $f(w)$. Some applications to
cases where we have control on the growth of the generations counting functions are given in Corollary \ref{cor13}.

\section{The Weierstrass representation of the automorphic group of an entire function, and the extra properties
in the case of a finite order}\label{sec3}

We will use (and for no particular reasons) the following two books: \cite{h}, Chapter IV, page 56,
and \cite{rc}, Chapter 15, page 87. This material is classical.

Let $f(z)$ be a non-constant entire function, and let $\{\Omega_i\}$ be a normal system of maximal domains 
of $f(z)$ ("fundamental domains" in Shimizu's terminology), and $\{\phi_{ij}\}$ is the corresponding automorphic group.
We view the difference $f(w)-f(z)$ as an entire function in $w$, and we view $z$ as a complex parameter. We have
the power series expansion $f(w)=a_0+a_1w+a_2w^2+\ldots $. Hence $f(w)-f(z)= a_1w+a_2w^2+\ldots-(a_1z+a_2z^2+\ldots)$.
As a function of $w$, it has a zero at the origin, $w=0$, if and only if $a_1z+a_2z^2+\ldots=f(z)-f(0)=0$ for the
particular value $z$ of the complex parameter. This is the case if $z=0$ (the trivial case). In all
other cases (where $f(z)-f(0)\ne 0$), the function $f(w)-f(z)$ (of $w$) does not vanish at the origin,
$w=0$. Hence the Weierstrass factorization theorem implies the following: \\
1) If $f(z)-f(0)\ne 0$, then there is a function $g(w,z)$, entire in $w$ and there are non-negative
integers $\lambda_n(w,z)$ which we will sometimes denote by $\lambda_n$, such that:
$$
f(w)-f(z)=\exp{\left(g(w,z)\right)}\prod_{n=1}^{\infty}E\left(\frac{w}{\phi_{0n}(z)},\lambda_n\right)=
$$
$$
=\exp{\left(g(w,z)\right)}\prod_{n=1}^{\infty}\left(1-\frac{w}{\phi_{0n}(z)}\right)e^{Q_{\lambda_n}(w/\phi_{0n}(z))},
$$
where if $\lambda_n>0$, then:
$$
Q_{\lambda_n}\left(\frac{w}{\phi_{0n}(z)}\right)=\left(\frac{w}{\phi_{0n}(z)}\right)+\frac{1}{2}
\left(\frac{w}{\phi_{0n}(z)}\right)^2+\ldots +\frac{1}{\lambda_n}\left(\frac{w}{\phi_{0n}(z)}\right)^{\lambda_n},
$$
and $Q_0(w/\phi_{0n}(z))\equiv 0$. \\
2) If $f(z)-f(0)=0$, then there is a natural number $m$, and there is an entire function in $w$, $h(w,z)$
(depending on each zero of $f(z)-f(0)$) and there are non-negative numbers $\lambda_n'(z)$,such that:
$$
f(w)-f(z)=w^m\exp{\left(h(w,z)\right)}\prod_{n=1,\phi_{0n}(z)\ne 0}^{\infty}E\left(\frac{w}{\phi_{0n}(z)},\lambda_n'\right),
$$
where $z$ satisfies: (a) $f(z)-f(0)=0$, (b) $\phi_{0n}(z)\ne 0$. \\
\\
Next, we flip the roles of the variable $w$ and the complex parameter $z$. We obtain: \\
1) If $f(w)-f(0)\ne 0$, then there is a function $g_1(z,w)$, entire in $z$ and there are non-negative
integers $\mu_n(z,w)$ which we will sometimes denote by $\mu_n$, such that:
$$
f(z)-f(w)=\exp{\left(g_1(z,w)\right)}\prod_{n=1}^{\infty}E\left(\frac{z}{\phi_{0n}(w)},\mu_n\right)=
$$
$$
=\exp{\left(g_1(z,w)\right)}\prod_{n=1}^{\infty}\left(1-\frac{z}{\phi_{0n}(w)}\right)e^{Q_{\mu_n}(z/\phi_{0n}(w))},
$$
2) If $f(w)-f(0)=0$, then with exactly the same values as in case 2 for $f(w)-f(z)$ above we have:
$$
f(z)-f(w)=z^m\exp{\left(h(z,w)\right)}\prod_{n=1,\phi_{0n}(w)\ne 0}^{\infty}E\left(\frac{z}{\phi_{0n}(w)},\lambda_n'\right),
$$
where $w$ satisfies: (a) $f(w)-f(0)=0$, (b) $\phi_{0n}(w)\ne 0$. \\
\\
Cases 1 are the generic cases (because cases 2 apply either to a discrete set of $z$, or to a discrete set of $w$).
By $f(w)-f(z)=-(f(z)-f(w))$ we obtain:
\begin{proposition}\label{prop1}
If $(f(w)-f(0))(f(z)-f(0))\ne 0$, then:
$$
\exp{\left(g(w,z)\right)}\prod_{n=1}^{\infty}\left(1-\frac{w}{\phi_{0n}(z)}\right)e^{Q_{\lambda_n}(w/\phi_{0n}(z))}=
$$
$$
=-\exp{\left(g_1(z,w)\right)}\prod_{n=1}^{\infty}\left(1-\frac{z}{\phi_{0n}(w)}\right)e^{Q_{\lambda_n'}(z/\phi_{0n}(w))},
$$
where $g(w,z)$ is entire in $w$, and $g_1(z,w)$ is entire in $z$. Moreover by the discussion in \cite{s1} that starts
on page 229 we may assume that all the automorphic functions $\phi_{0n}$ are holomorphic in the (interior) maximal
domains $\Omega_i$ of the system we fixed. Hence $g(w,z)$ is $z$-holomorphic in the $\Omega_i$'s and $g_1(z,w)$ is
$w$-holomorphic there.
\end{proposition}

\begin{proposition}\label{prop2}
If $(f(w)-f(0))(f(z)-f(0))\ne 0$, then:
$$
\frac{\partial g(w,z)}{\partial w}+\sum_{n=1}^{\infty}\frac{(w/\phi_{0n}(z))^{\lambda_n}}{w-\phi_{0n}(z)}=
\frac{\partial g_1(z,w)}{\partial w}+\sum_{n=1}^{\infty}\left(\frac{\phi_{0n}'(w)}{\phi_{0n}(w)}\right)z
\frac{(z/\phi_{0n}(w))^{\lambda_n'}}{\phi_{0n}(w)-z}.
$$
\end{proposition}
\noindent
{\bf Proof.} \\
Take the logarithm of the two sides in the identity of proposition \ref{prop1}, and then $\partial/\partial w$
both sides and simplify. We note that both $g(w,z)$ and $g_1(z,w)$ are holomorphic (usually not entire) in both
variables in the appropriate domains of the $\mathbb{C}\times\{\mathbb{C}-{\rm discrete}\,\,{\rm set}\}$. $\qed $ \\

\begin{remark}\label{rem11}
If we $\partial/\partial w$ both sides of $f(\phi_{0n}(w))=f(w)$, then we obtain $f'(\phi_{0n}(w))\phi_{0n}'(w)=f'(w)$
in the appropriate domain $\Omega $. By Shimizu this domain, is such that $\mathbb{C}-\Omega$ contains no continuum.
\end{remark}

\begin{remark}\label{rem12}
If we take $\partial/\partial z$ instead of $\partial/\partial w$, we get the symmetric identity:
$$
\frac{\partial g(w,z)}{\partial z}+\sum_{n=1}^{\infty}\left(\frac{\phi_{0n}'(z)}{\phi_{0n}(z)}\right)w
\frac{(w/\phi_{0n}(z))^{\lambda_n}}{\phi_{0n}(z)-w}=
$$
$$
=\frac{\partial g_1(z,w)}{\partial z}+\sum_{n=1}^{\infty}\frac{(z/\phi_{0n}(w))^{\lambda_n'}}{z-\phi_{0n}(w)}.
$$
\end{remark}
\noindent
The entire $w$-functions $g(w,z)$ in the generic Weierstrass factorization of $f(w)-f(z)$ are special
regarding their relation to the action of the automorphic group of $f(z)$ (note that here we take $z$ as
the variable).

\begin{proposition}\label{prop3}
If $f(z)-f(0)\ne 0$, then there is a function $g(w,z)$, entire in $w$ and there are non-negative integers
$\lambda_n=\lambda_n(w,\phi_{0n}(z))$ such that:
$$
f(w)-f(z)=\exp{\left(g(w,z)\right)}\prod_{n=1}^{\infty}\left(1-\frac{w}{\phi_{0n}(z)}\right)e^{Q_{\lambda_n}(w/\phi_{0n}(z))}.
$$
If $\forall\,n$ we have $\lambda_n=\lambda$, a constant value independent of $n$, then
the function $\exp{\left(g(w,z)\right)}$ is $z$-invariant with respect to the action of the automorphic
group of $f(z)$. This means that $\exp{\left(g(w,\phi_{ij}(z))\right)}=\exp{\left(g(w,z)\right)}$ for every
element $\phi_{ij}$ in the automorphic group. Our assumption on the $\lambda_n\equiv\lambda$ is valid whenever 
the entire function $f(w)$ has a finite order.
\end{proposition}
\noindent
{\bf Proof.} \\
The first part is just Weierstrass factorization theorem applied to $f(w)-f(z)$ (as an entire
function of $w$). Let $\phi_{ij}$ be any automorphic function of $f(z)$. This means that
$f(\phi_{ij}(z))=f(z)$. We note that:
$$
\prod_{n=1}^{\infty}\left(1-\frac{w}{\phi_{0n}(\phi_{ij}(z))}\right)
e^{Q_{\lambda_n(w,\phi_{0n}(\phi_{ij}(z)))}(w/\phi_{0n}(\phi_{ij}(z)))}=
$$
$$
=\prod_{n=1}^{\infty}\left(1-\frac{w}{\phi_{0n}(z)}\right)e^{Q_{\lambda_n}(w/\phi_{0n}(z))},
$$
because the left hand side product is a product of a permutation of the factors of the right hand side. This
follows by the assumption on the $\lambda_n\equiv\lambda$, independent of $n$.
By the convergence uniformly on compacta the two products are equal to one another. Hence the quotient function:
$$
\left(f(w)-f(z)\right)\left/\prod_{n=1}^{\infty}\left(1-\frac{w}{\phi_{0n}(z)}\right)e^{Q_{\lambda_n}(w/\phi_{0n}(z))}\right.
$$
is invariant with respect to the action of the elements $\{\phi_{ij}\}$ of the automorphic group of $f(z)$. But
by the Weierstrass identity above, this quotient function equals $\exp{(g(w,z))}$. $\qed $ \\
\\
We now indicate conclusions of this proposition which are of a different character than the ones
above. We will use the following:

\begin{lemma}\label{lem1}
If $u(z)$ and $v(z)$ are non-constant entire functions and if ${\rm Aut}(v)$ is a subgroup of
${\rm Aut}(u)$, then there exists a function $h(w)$, holomorphic on the image of $v$ (i.e. on
$v(\mathbb{C})$) such that $u(z)=h(v(z))$.
\end{lemma}
\noindent
{\bf Proof.} \\
That follows using the methods in \cite{r}, however, we will prove it using a result in \cite{er}. Namely,
we will make use of the result that appears on the last paragraph on page 334 and continues on the first
paragraph on the next page, 335 in \cite{er}. Thus we claim that $v\le u$, where the partial order is
defined in \cite{er} (where we indicated). To prove that we need to show that $v(z)=v(w)$ implies that
$u(z)=u(w)$, $\forall\,z,w\in\mathbb{C}$. Thus assuming that $v(z)=v(w)$, it follows that $\exists\,\Phi\in{\rm Aut}(v)$,
such that $w=\Phi(z)$. By an assumption we have, it follows that $\Phi\in{\rm Aut}(u)$. Hence it indeed follows
that $u(z)=u(w)$, and we proved that $v\le u$. The claim in our lemma now follows by \cite{er}. $\qed $ \\
\\
Using the Shimizu's \cite{s1} (or \cite{r}) we can re-write Lemma \ref{lem1} in a geometric manner:

\begin{lemma}\label{lem2}
If $u(z)$ and $v(z)$ are non-constant entire functions and if ${\rm Aut}(v)$ is a subgroup of
${\rm Aut}(u)$, then a normal system of maximal domains $\{\Omega_j\}$ of $v(z)$, is composed 
of maximal domains $\Omega_j$ each of which is tiled by some elements of the same normal system of
maximal domains of $u(z)$.
\end{lemma}
\noindent
Specializing Lemma \ref{lem2} to the setting of Proposition \ref{prop3} we obtain:

\begin{corollary}\label{cor1}
Under the assumptions of Proposition \ref{prop3}: The normal system of maximal domains $\{\Omega_j\}$
of $f(z)$ which induces the Weierstrass factorization of $f(w)-f(z)$, is composed of maximal domains
$\Omega_j$, each of which is tiled by some elements of the same normal system of maximal domains of
$\exp{\left(g(w,z)\right)}$ as a function of $z$ and for a fixed $w$.
\end{corollary}
\noindent
Specializing Lemma \ref{lem1} to the setting of Proposition \ref{prop3} we obtain:

\begin{corollary}\label{cor2}
Under the assumptions of Proposition \ref{prop3}:
There is a function $F(w,w_1)$, holomorphic in $\mathbb{C}\times f(\mathbb{C})$ such that
$\exp{\left(g(w,z)\right)}=\exp{\left(F(w,f(z))\right)}$.
\end{corollary}
\noindent
So the Weierstrass factorization of $f(w)-f(z)$ that is described in Proposition \ref{prop3} is
special, and we add the extra information in the following:

\begin{theorem}\label{thm1}
Let $f$ be a non-constant entire function of a finite order, and assuming that $f(z)-f(0)\ne 0$ we have the following
expansion:
$$
f(w)-f(z)=\exp{\left(F(w,f(z))\right)}\prod_{n=1}^{\infty}\left(1-\frac{w}{\phi_{0n}(z)}\right)e^{Q_{\lambda}(w/\phi_{0n}(z))}.
$$
In particular, the canonical infinite product is in fact an holomorphic function of $(w,f(z))$, for all 
$w\in\mathbb{C}$ and $z\in(\mathbb{C}-Z(f-f(0)))$.
\end{theorem}

\begin{corollary}\label{cor3}
Let $f$ be a non-constant entire function of a finite order, then we have the following cycle relation:
$$
\sum_{j=1}^{N}\exp{\left(F(z_j,f(z_{j+1}))\right)}
\prod_{n=1}^{\infty}\left(1-\frac{z_j}{\phi_{0n}(z_{j+1})}\right)e^{Q_{\lambda}(z_j/\phi_{0n}(z_{j+1}))}\equiv 0,
$$
for any $N$ independent variables: $z_1,\ldots,z_{N}$, where we agree that $z_{N+1}=z_1$.
\end{corollary}
\noindent
{\bf Proof.} \\
\underline{\bf Method 1:} Using the identity in Theorem \ref{thm1} we have the following:
$$
f(z_1)=f(z_2)+\exp{\left(F(z_1,f(z_2))\right)}
\prod_{n=1}^{\infty}\left(1-\frac{z_1}{\phi_{0n}(z_2)}\right)e^{Q_{\lambda}(z_1/\phi_{0n}(z_2))},
$$
$$
f(z_2)=f(z_3)+\exp{\left(F(z_2,f(z_3))\right)}
\prod_{n=1}^{\infty}\left(1-\frac{z_2}{\phi_{0n}(z_3)}\right)e^{Q_{\lambda}(z_2/\phi_{0n}(z_3))},
$$
$$ \vdots $$
$$
f(z_{N-1})=f(z_N)+\exp{\left(F(z_{N-1},f(z_N))\right)}
\prod_{n=1}^{\infty}\left(1-\frac{z_{N-1}}{\phi_{0n}(z_N)}\right)e^{Q_{\lambda}(z_{N-1}/\phi_{0n}(z_N))},
$$
$$
f(z_N)=f(z_1)+\exp{\left(F(z_N,f(z_1))\right)}
\prod_{n=1}^{\infty}\left(1-\frac{z_N}{\phi_{0n}(z_1)}\right)e^{Q_{\lambda}(z_N/\phi_{0n}(z_1))}.
$$
We plug these identities successively each in its predecessor and eventually cancel out $f(z_1)$ from
both sides of the equation. \\
\\
\underline{\bf Method 2:} The cycle relation is merely the Weierstrass factorization of each term in the
following telescopic identity:
$$
(f(z_1)-f(z_2))+(f(z_2)-f(z_3))+\ldots+(f(z_{N-1})-f(z_N))+(f(z_N)-f(z_1))\equiv 0.\,\,\,\,\,\,\,\,\,\,\qed
$$

\begin{corollary}\label{cor4}
Let $f$ be a non-constant entire function of a finite order, then we have the following chain relation:
$$
f(z_1)-f(z_{N+1})=\sum_{j=1}^{N}\exp{\left(F(z_j,f(z_{j+1}))\right)}
\prod_{n=1}^{\infty}\left(1-\frac{z_j}{\phi_{0n}(z_{j+1})}\right)e^{Q_{\lambda}(z_j/\phi_{0n}(z_{j+1}))}
$$
for any $N+1$ independent variables $z_1,z_2,\ldots,z_{N+1}$.
\end{corollary}
\noindent
{\bf Proof.} \\
It is clear how to adopt any of the two methods of the proof we gave to Corollary \ref{cor3}. $\qed $

\begin{corollary}\label{cor5}
$$
\sum_{j=1}^{N}\exp{\left(F(z_j,f(z_{j+1}))\right)}
\prod_{n=1}^{\infty}\left(1-\frac{z_j}{\phi_{0n}(z_{j+1})}\right)e^{Q_{\lambda}(z_j/\phi_{0n}(z_{j+1}))}\equiv
$$
$$
\equiv \exp{\left(F(z_1,f(z_{N+1}))\right)}
\prod_{n=1}^{\infty}\left(1-\frac{z_1}{\phi_{0n}(z_{N+1})}\right)e^{Q_{\lambda}(z_1/\phi_{0n}(z_{N+1}))}
$$
for any $N+1$ independent variables $z_1,z_2,\ldots,z_{N+1}$.
\end{corollary}

\begin{remark}\label{rem15}
The results that begin in our Corollary \ref{cor1} and end in Corollary \ref{cor5} should be carefully interpreted,
because $\exp{(g(w,z))}$ is entire in $w$, but is not known to be entire in $z$. This is because we are dealing
only with those values of the parameter $z$ for which $f(z)-f(0)\ne 0$. Thus apriori it is not clear what is the 
meaning of "elements of a normal system of maximal domains of $\exp{(g(w,z)}$ for a fixed $w$". This notion was
defined by Shimizu only for meromorphic functions, but we do not know that $\exp{(g(w,z))}$ is meromorphic in $z$.
\end{remark}

\begin{remark}\label{rem16}
The cycle relation in Corollary \ref{cor3} and the chain relation in Corollary \ref{cor4} and in Corollary \ref{cor5}
resemble the fact that the value of a path integral is independent of the path that connects the two endpoints 
for a conservative field. In our setting one may think of $f(z)$ as the "potential" of the complicated
Weierstrass products that appear within the sum. This resemblance originates in the elementary fact that
the sum of telescopic series depends only on the initial and the terminal points.
\end{remark}

\section{Conclusions from Proposition \ref{prop3} in case we have no monodromy}\label{sec4}

The conclusions from Proposition \ref{prop3} that were derived in the previous section, originated in Lemma \ref{lem1}
and Lemma \ref{lem2}. These, in turn were based on a result of Eremenko and Rubel in \cite{er}. However, their 
result used the monodromy principle that was available in their setting. What if we have no monodromy? Of course
the conclusion of Proposition \ref{prop3} is still valid but we have no composition relations between $f$ and the
Weierstrass factor $\exp{(g(w,z))}$. We will outline in this section what can we still conclude. \\
 
\begin{corollary}\label{cor0}
For any element $\phi_{ij}(z)$ of the automorphic group of $f(z)$, a non-constant entire function of a finite order,
there exists an integer $n_{ij}\in\mathbb{Z}$ so that $g(w,\phi_{ij}(z))=g(w,z)+2\pi n_{ij}\cdot i$.
\end{corollary}

\begin{remark}\label{rem13}
It would be nice to compute $g(w,z)$ for different entire functions and to check the various identities
we obtained. Later on we will carry such a computation for the exponential function.
\end{remark}

\begin{remark}\label{rem14}
The $g(w,z)$ function translates the group of automorphic functions (composition of mappings is its
binary operation), into a subgroup of $(\mathbb{Z},+)$. \\
{\bf Proof.} \\
Let us take two automorphic functions $\phi_{ij}$ and $\phi_{\alpha,\beta}$. Then: $g(w,\phi_{ij}(z))=g(w,z)+2\pi n_{ij}\cdot i$
and $g(w,\phi_{\alpha\beta}(z))=g(w,z)+2\pi n_{\alpha\beta}\cdot i$.
$$
g(w,\phi_{ij}(\phi_{\alpha,\beta}(z)))=g(w,\phi_{\alpha,\beta}(z))+2\pi n_{ij}\cdot i=
$$
$$
=(g(w,z)+2\pi n_{\alpha\beta}\cdot i)+2\pi n_{ij}\cdot i=g(w,z)+2\pi(n_{\alpha\beta}+n_{ij})\cdot i.
$$
$\qed $
\end{remark}

\begin{corollary}\label{cor05}
The automorphic group of a non-constant entire function of a finite order, is homomorphic to a subgroup of $(\mathbb{Z},+)$
(which is not always the trivial homomorphism).
\end{corollary}
\noindent
{\bf Proof.} \\
Let $f$ be a non-constant entire function. Let us denote by ${\rm Aut}(f)$ it's group of automorphic functions.
Let us define the mapping $T:\,{\rm Aut}(f)\rightarrow\mathbb{Z}$ by the formula suggested by Corollary \ref{cor0}, i.e.
$$
T(\Phi)=\frac{1}{2\pi\cdot i}\left(g(w,\Phi(z))-g(w,z)\right).
$$
Here the complex numbers $w, z\in \mathbb{C}$ are completely arbitrary within the domain of the definition of $g(w,z)$.
Then by Remark \ref{rem14} we have $\forall\,\Phi_1, \Phi_2\in {\rm Aut}(f)$ we have the identity 
$T(\Phi_1\circ\Phi_2)=T(\Phi_1)+T(\Phi_2)$. $\qed $ \\

\begin{corollary}\label{cor055}
Let $f$ be a non-constant entire function of a finite order. If $\Phi\in{\rm Aut}(f)$ is an element of a finite
order, then $\Phi\in\ker (T)$.
\end{corollary}
\noindent
{\bf Proof.} \\
The only finite subgroup of the infinite cyclic group $(\mathbb{Z},+)$ is the trivial subgroup $\{0\}$. $\qed $ \\

We deduce a family of functional-arithmetical identities from Corollary \ref{cor05}. For that we will use the
obvious short notation for repeated composition of functions.

\begin{definition}\label{def1}
Let $h_j(z)$, $j=1,\ldots,n$ be $n$ complex valued functions for which the repeated composition makes sense.
We will denote:
$$
\left(h_1\circ\ldots\circ h_n\right)=\bigcirc_{j=1}^{n}h_j.
$$
\end{definition}

\begin{corollary}\label{cor075}
Let $f$ be a non-constant entire function of a finite order, and let $\Phi_1,\ldots,\Phi_n\in{\rm Aut}(f)$ ($n\ge 2$). Then whenever
the repeated composition makes sense we have the identity:
$$
g\left(w,\left(\bigcirc_{j=1}^{n}\Phi_j\right)\left(z\right)\right)-\sum_{j=1}^{n}g\left(w,\Phi_j(z)\right)+(n-1)g(w,z)\equiv 0,
$$
$\forall\,(w,z)\in\mathbb{C}\times\left(\mathbb{C}-{\rm a}\,{\rm discrete}\,{\rm set}\right)$.
\end{corollary}
\noindent
{\bf Proof.} \\
The proof is inductive on $n\in\mathbb{Z}_{\ge 2}$. For $n=2$ we have (using the map $T$ in the proof of Corollary
\ref{cor05}: $T(\Phi_1\circ\Phi_2)=T(\Phi_1)+T(\Phi_2)$, i.e.:
$$
\frac{1}{2\pi\cdot i}\left(g(w,(\Phi_1\circ\Phi_2)(z))-g(w,z)\right)=
$$
$$
=\frac{1}{2\pi\cdot i}\left(g(w,\Phi_1(z))-g(w,z)\right)
+\frac{1}{2\pi\cdot i}\left(g(w,\Phi_2(z))-g(w,z)\right).
$$
Hence:
$$
g(w,(\Phi_1\circ\Phi_2)(z))-g(w,\Phi_1(z))-g(w,\Phi_2(z))+g(w,z)\equiv 0.
$$
This completes the case $n=2$. Similarly the case $n=3$ follows in a very similar manner from:
$T(\Phi_1\circ\Phi_2\circ\Phi_3)=T(\Phi_1)+T(\Phi_2)+T(\Phi_3)$, i.e.:
$$
\frac{1}{2\pi\cdot i}\left(g(w,(\Phi_1\circ\Phi_2\circ\Phi_3)(z))-g(w,z)\right)=
$$
$$
=\frac{1}{2\pi\cdot i}\left(g(w,\Phi_1(z))-g(w,z)\right)
+\frac{1}{2\pi\cdot i}\left(g(w,\Phi_2(z))-g(w,z)\right)
+\frac{1}{2\pi\cdot i}\left(g(w,\Phi_3(z))-g(w,z)\right).
$$
Hence:
$$
g(w,(\Phi_1\circ\Phi_2\circ\Phi_3)(z))-g(w,\Phi_1(z))-g(w,\Phi_2(z))-g(w,\Phi_3(z))+2\cdot g(w,z)\equiv 0.
$$
This completes the case $n=3$, etc... $\qed $ \\
\\
We have now a direct connection between composition arithmetic and lattice (integral) arithmetic. Here is
a straightforward example. We might think of $|T(\Phi)|$ has the distance between $g(w,\Phi(z))$ and $g(w,z)$.
Given two non trivial elements $\Phi_1,\Phi_2\in {\rm Aut}(f)$, i.e. elements for which the corresponding
distances are not $0$ ($|T(\Phi_1)T(\Phi_2)|>0$) can we find an non-trivial element with a shorter
distance. Here is a possible way to go about solving that:
\begin{corollary}\label{cor08}
Let $f$ be a non-constant entire function of a finite order and let $\Phi_1,\Phi_2\in{\rm Aut}(f)$.
Suppose that we have:
$$
m=\frac{1}{2\pi\cdot i}\left(g(w,\Phi_1(z))-g(w,z)\right),\,\,n=\frac{1}{2\pi\cdot i}\left(g(w,\Phi_2(z))-g(w,z)\right),
$$
where $m\cdot n\ne 0$. Let $d=a\cdot m+b\cdot n=(m,n)$ the lcm of the integers $m$ and $n$. Here $d$ can be the positive
or the negative lcm. Then if we define:
$$
\Phi=\Phi_1^{\circ a}\circ\Phi_2^{\circ b},
$$
then we have:
$$
d=\frac{1}{2\pi\cdot i}\left(g(w,\Phi(z))-g(w,z)\right).
$$
\end{corollary}
\noindent
{\bf Proof.} \\
We clearly define
$$
\Phi_1^{\circ a}=\left\{\begin{array}{lll} \Phi_1\circ\ldots\circ\Phi_1 & {\rm if} & a>0 \\
(\Phi_1)^{-1}\circ\ldots\circ(\Phi_1)^{-1} & {\rm if} & (-a)>0 \end{array}\right..
$$
We note that we have: $T(\Phi_1^{\circ a})=a\cdot T(\Phi_1)$. Hence $T(\Phi)=T(\Phi_1^{\circ a}\circ\Phi_2^{\circ b})=
a\cdot T(\Phi_1)+b\cdot T(\Phi_2)=a\cdot m+b\cdot n=d$. $\qed $ \\

\begin{corollary}\label{cor09}
Let $f$ be a non-constant entire function of a finite order and let $\Phi_1,\Phi_2\in{\rm Aut}(f)$. Then we have:
$$
\frac{1}{2\pi\cdot i}\left(g(w,\Phi_1(z)^{\circ T(\Phi_2)})-g(w,z)\right)=
\frac{1}{2\pi\cdot i}\left(g(w,\Phi_2(z)^{\circ T(\Phi_1)})-g(w,z)\right).
$$
\end{corollary}
\noindent
{\bf Proof.} \\
This follows by the fact that: 
$$
T(\Phi_1^{\circ T(\Phi_2)})=T(\Phi_2)T(\Phi_1)=T(\Phi_1)T(\Phi_2)=T(\Phi_2^{\circ T(\Phi_1)}).\,\,\,\,\,\qed
$$
\begin{corollary}\label{cor095}
Let $f$ be a non-constant entire function of a finite order and let $\Phi\in{\rm Aut}(f)$ satisfy the condition that the number
$$
\frac{1}{2\pi\cdot i}\left(g(w,\Phi(z))-g(w,z)\right)
$$ 
is a prime number $p$. If $\Phi=\Psi^{\circ k}$ for some $\Psi\in{\rm Aut}(f)$, then either $T(\Psi)$ equals $1$
or equals $p$.
\end{corollary}

\begin{corollary}\label{cor096}
Let $f$ be a non-constant entire function of a finite order and let $\Phi\in{\rm Aut}(f)$ satisfy the condition that the number
$$of 
\frac{1}{2\pi\cdot i}\left(g(w,\Phi(z))-g(w,z)\right)
$$ 
is a prime number $p$. If $\Phi=\Phi_1\circ\ldots\circ\Phi_n$, where $n\in\mathbb{Z}^+$, and where for $j=1,\ldots,n$,
$\Phi_j\in{\rm Aut}(f)$, then there is a single index $k$, between $1$ and $n$ such that $T(\Phi_k)=p$ while for
$j\in\{1,\ldots,n\}-\{k\}$, $T(\Phi_j)=1$.
\end{corollary}
\noindent
Corollaries \ref{cor08}, \ref{cor09}, \ref{cor095} and \ref{cor096} are all particular cases of the principle that
the arithmetic of composition of automorphic functions of a non-constant entire function have an analog
in the arithmetic of the integers, $\mathbb{Z}$. We can describe the general principle in the following:

\begin{theorem}\label{thm0}{\bf The finiteness of the decomposition of automorphic functions.}
Let $f$ be a non-constant entire function of a finite order and let $\Phi\in{\rm Aut}(f)$ satisfy the condition:
$$
\frac{1}{2\pi\cdot i}\left(g(w,\Phi(z))-g(w,z)\right)=N\in\mathbb{Z}-\{0\}.
$$
Then any decomposition of $\Phi$ into a composition of automorphic functions of $f$:
$$
\Phi=\Phi_1\circ\ldots\circ\Phi_n,\,\,\,\,\,\,\Phi_1,\ldots,\Phi_n\in{\rm Aut}(f),
$$
has the following properties: \\
1) $n$ could be any natural number with no apriori upper bound. \\
2) We have the Diophantine identity:
$$
\frac{1}{2\pi\cdot i}\left(g(w,\Phi(z))-g(w,z)\right)=
\prod_{j=1}^{n}\left(\frac{1}{2\pi\cdot i}\left(g(w,\Phi_j(z))-g(w,z)\right)\right).
$$
If we call an automorphic function $\Phi_j$ an arithmetical unit, if it satisfies:
$$
\left|\frac{1}{2\pi\cdot i}\left(g(w,\Phi_j(z))-g(w,z)\right)\right|=1,
$$
then in any such decomposition of $\Phi$, the set:
$$
\left\{\frac{1}{2\pi\cdot i}\left(g(w,\Phi_j(z))-g(w,z)\right)\ne\pm 1\right\},
$$
is a set of non-unit divisors of $N$ whose product is $N$, and all the other factors
belong to arithmetical units. In particular for any $\Phi\in{\rm Aut}(f)$, the number 
of different decompositions that differ in their non-units is bounded above by:
$$
\sum m!\cdot|\{\{k_1,\ldots,k_m\}\,|\,k_1\cdot\ldots\cdot k_m=N,\,|k_1|,\ldots,|k_m|>1\}|.
$$
The weights $m!$ must be present because composition of functions, unlike multiplication
of integers in a non-commutative binary operation.
\end{theorem}

\section{The cycle relation and the chain relation in the general case}\label{sec5}

The results in section \ref{sec3} dealt mostly with entire functions of a finite order. The key result
was Proposition \ref{prop3} and we assumed that $\lambda_n\equiv\lambda$ independent of $n$. This essentially
is the assumption that $f$ has a finite order. In this section we point out at the results if this
assumption is dropped out.

\begin{corollary}\label{cor53}
Let $f$ be a non-constant entire function, then we have the following cycle relation:
$$
\sum_{j=1}^{N}\exp{\left(g(z_j,z_{j+1})\right)}
\prod_{n=1}^{\infty}\left(1-\frac{z_j}{\phi_{0n}(z_{j+1})}\right)e^{Q_{\lambda_n}(z_j/\phi_{0n}(z_{j+1}))}\equiv 0,
$$
for any $N$ independent variables: $z_1,\ldots,z_{N}$, where we agree that $z_{N+1}=z_1$.
\end{corollary}
\noindent
{\bf Proof.} \\
\underline{\bf Method 1:} Using the first identity in Proposition \ref{prop3} where no finite order assumption is needed,
we have the following:
$$
f(z_1)=f(z_2)+\exp{\left(g(z_1,z_2)\right)}
\prod_{n=1}^{\infty}\left(1-\frac{z_1}{\phi_{0n}(z_2)}\right)e^{Q_{\lambda_n}(z_1/\phi_{0n}(z_2))},
$$
$$
f(z_2)=f(z_3)+\exp{\left(g(z_2,z_3)\right)}
\prod_{n=1}^{\infty}\left(1-\frac{z_2}{\phi_{0n}(z_3)}\right)e^{Q_{\lambda_n}(z_2/\phi_{0n}(z_3))},
$$
$$ \vdots $$
$$
f(z_{N-1})=f(z_N)+\exp{\left(g(z_{N-1},z_N)\right)}
\prod_{n=1}^{\infty}\left(1-\frac{z_{N-1}}{\phi_{0n}(z_N)}\right)e^{Q_{\lambda_n}(z_{N-1}/\phi_{0n}(z_N))},
$$
$$
f(z_N)=f(z_1)+\exp{\left(g(z_N,z_1)\right)}
\prod_{n=1}^{\infty}\left(1-\frac{z_N}{\phi_{0n}(z_1)}\right)e^{Q_{\lambda_n}(z_N/\phi_{0n}(z_1))}.
$$
We plug these identities successively each in its predecessor and eventually cancel out $f(z_1)$ from
both sides of the equation. \\
\\
\underline{\bf Method 2:} The cycle relation is merely the Weierstrass factorization of each term in the
following telescopic identity:
$$
(f(z_1)-f(z_2))+(f(z_2)-f(z_3))+\ldots+(f(z_{N-1})-f(z_N))+(f(z_N)-f(z_1))\equiv 0.\,\,\,\,\,\,\,\,\,\,\qed
$$

\begin{corollary}\label{cor54}
Let $f$ be a non-constant entire function, then we have the following chain relation:
$$
f(z_1)-f(z_{N+1})=\sum_{j=1}^{N}\exp{\left(g(z_j,z_{j+1})\right)}
\prod_{n=1}^{\infty}\left(1-\frac{z_j}{\phi_{0n}(z_{j+1})}\right)e^{Q_{\lambda_n}(z_j/\phi_{0n}(z_{j+1}))}
$$
for any $N+1$ independent variables $z_1,z_2,\ldots,z_{N+1}$.
\end{corollary}
\noindent
{\bf Proof.} \\
It is clear how to adopt any of the two methods of the proof we gave to Corollary \ref{cor53}. $\qed $

\begin{corollary}\label{cor55}
$$
\sum_{j=1}^{N}\exp{\left(g(z_j,z_{j+1})\right)}
\prod_{n=1}^{\infty}\left(1-\frac{z_j}{\phi_{0n}(z_{j+1})}\right)e^{Q_{\lambda_n}(z_j/\phi_{0n}(z_{j+1}))}\equiv
$$
$$
\equiv \exp{\left(g(z_1,z_{N+1})\right)}
\prod_{n=1}^{\infty}\left(1-\frac{z_1}{\phi_{0n}(z_{N+1})}\right)e^{Q_{\lambda_n}(z_1/\phi_{0n}(z_{N+1}))}
$$
for any $N+1$ independent variables $z_1,z_2,\ldots,z_{N+1}$.
\end{corollary}

\section{Examples (mostly the exponential function) and the role played by the assumption that we have
some summation method the infinite series:
$$
\sum_{n=1}^{\infty}Q_{\lambda_n}\left(\frac{w}{\phi_{0n}(z)}\right),
$$
for the reconstruction of $f$ from ${\rm Aut}(f)$}\label{sec6}
 
Let $f(z)=e^z$. We consider the following natural system of maximal domains of $f(z)$:
$$
\{\Omega_n=\{z\in\mathbb{C}\,|\,2\pi\cdot i\cdot n<\Im z< 2\pi\cdot i\cdot(n+1)\}\,|\,n\in\mathbb{Z}\}.
$$
This induces the infinite cyclic automorphic group:
$$
{\rm Aut}(e^z)=\{z+2\pi\cdot i\cdot n\,|\,n\in\mathbb{Z}\}=<z+2\pi\cdot i>.
$$
The discrete exceptional set of $z$ is the solution set of the equation $e^z-e^0=0$. So this is the
discrete set $\{2\pi\cdot i\cdot n\,|\,n\in\mathbb{Z}\}$. Using the representation of Proposition \ref{prop3}
we clearly can choose the sequence $\lambda_n\equiv 1$, $\forall\,n\in\mathbb{Z}$. Thus for 
$z\not\in\{2\pi\cdot i\cdot n\,|\,n\in\mathbb{Z}\}$, we have:
$$
e^w-e^z=\exp{\left(g(w,z)\right)}\prod_{n\in\mathbb{Z}}\left(1-\frac{w}{z+2\pi\cdot i\cdot n}\right)
e^{(w/(z+2\pi\cdot i\cdot n))}.
$$
We can group together symmetric pairs $n$ and $-n$, where $n\in\mathbb{Z}^+$. We compute the corresponding products:
$$
\left(1-\frac{w}{z+2\pi\cdot i\cdot n}\right)e^{(w/(z+2\pi\cdot i\cdot n))}\times
\left(1-\frac{w}{z-2\pi\cdot i\cdot n}\right)e^{(w/(z-2\pi\cdot i\cdot n))}
$$
and we can write the final result in two forms as follows:
$$
e^w-e^z=\exp{\left(g(w,z)\right)}\left(1-\frac{w}{z}\right)e^{w/z}\prod_{n=1}^{\infty}
\left(\frac{(z-w)^2+4\pi^2n^2}{z^2+4\pi^2n^2}\right)e^{2zw/(z^2+4\pi^2n^2)}=
$$
$$
=\exp{\left(g(w,z)\right)}\left(1-\frac{w}{z}\right)e^{w/z}\prod_{n=1}^{\infty}
\left(1-\frac{z^2-(z-w)^2}{z^2+4\pi^2n^2}\right)e^{2zw/(z^2+4\pi^2n^2)},
$$
$\forall\,w\in\mathbb{C}$, $\forall\,z\in\mathbb{C}-\{2\pi\cdot i\cdot n\,|\,n\in\mathbb{Z}\}$.
Here $g(w,z)$ is entire in $w$ and holomorphic in $z\not\in\{2\pi\cdot i\cdot n\,|\,n\in\mathbb{Z}\}$.
Next, we note that if we replace $z$ by $z+2\pi\cdot i\cdot k$ for some $k\in\mathbb{Z}$, then clearly:
$$
\prod_{n\in\mathbb{Z}}\left(1-\frac{w}{(z+2\pi\cdot i\cdot k)+2\pi\cdot i\cdot n}\right)
e^{(w/((z+2\pi\cdot i\cdot k)+2\pi\cdot i\cdot n))}=
$$
$$
=\prod_{n\in\mathbb{Z}}\left(1-\frac{w}{z+2\pi\cdot i\cdot(n+k)}\right)
e^{(w/(z+2\pi\cdot i\cdot(n+k)))}=
$$
$$
=\prod_{n\in\mathbb{Z}}\left(1-\frac{w}{z+2\pi\cdot i\cdot n}\right)
e^{(w/(z+2\pi\cdot i\cdot n))}.
$$
Also $e^w-e^{z+2\pi\cdot i\cdot k}=e^w-e^z$. Hence the basic Weierstrass factorization:
$$
e^w-e^z=\exp{\left(g(w,z)\right)}\prod_{n\in\mathbb{Z}}\left(1-\frac{w}{z+2\pi\cdot i\cdot n}\right)
e^{(w/(z+2\pi\cdot i\cdot n))},
$$
implies that indeed we have $\exp{(g(w,z+2\pi\cdot i\cdot k))}=\exp{(g(w,z))}$. Next, let us consider
$e^w-1$. This entire function has simple zeros at $\{2\pi\cdot i\cdot n\,|\,n\in\mathbb{Z}\}$ and only
there. So using the standard Weierstrass factorization we obtain an identity of the following form:
$$
e^w-1=\exp{(h(w))}\cdot w\cdot\prod_{n=1}^{\infty}\left(1-\frac{w^2}{4\pi^2n^2}\right).
$$
This follows by taking the symmetric order of factors in:
$$
e^w-1=\exp{(h(w))}\cdot w\cdot\prod_{n\in\mathbb{Z}-\{0\}}\left(1-\frac{w}{2\pi\cdot i\cdot n}\right)
e^{w/(2\pi\cdot i\cdot n)}.
$$
Using this identity we obtain:
$$
e^w-e^z=e^z(e^{w-z}-1)=e^z\exp{(h(w-z))}\cdot(w-z)\prod_{n=1}^{\infty}\left(1+\frac{(w-z)^2}{4\pi^2n^2}\right)=
$$
$$
=e^z\exp{(h(w-z))}\cdot(w-z)\prod_{n\in\mathbb{Z}-\{0\}}\left(1-\frac{w-z}{2\pi\cdot i\cdot n}\right)
e^{(w-z)/(2\pi\cdot i\cdot n)}.
$$
Thus we obtained two different identities:
$$
e^w-e^z=\exp{(g(w,z))}\prod_{n\in\mathbb{Z}}\left(1-\frac{w}{z+2\pi\cdot i\cdot n}\right)e^{w/(z+2\pi\cdot i\cdot n)}=
$$
$$
=e^z\exp{(h(w-z))}\cdot(w-z)\prod_{n\in\mathbb{Z}-\{0\}}\left(1-\frac{w-z}{2\pi\cdot i\cdot n}\right)
e^{(w-z)/(2\pi\cdot i\cdot n)},
$$
or by symmetric multiplication:
$$
e^w-e^z=\exp{(g(w,z))}\left(1-\frac{w}{z}\right)e^{w/z}\prod_{n=1}^{\infty}\left(1-\frac{z^2-(w-z)^2}{z^2+4\pi^2n^2}\right)
e^{(2zw)/(z^2+4\pi^2n^2)}=
$$
$$
=e^z\exp{(h(w-z))}\cdot(w-z)\prod_{n=1}^{\infty}\left(1+\frac{(w-z)^2}{4\pi^2n^2}\right).
$$
This is different from the unique factorization of polynomials. We have no uniqueness of product representation. A
well known phenomenon. Before proceeding to the computation of the Weierstrass factor $g(w,z)$, which is not
trivial even for the exponential function, let us solve first the polynomial case. We start with the following
quadratic $f(z)=z^2+z$ and we note that $f(w)-f(z)=(w-z)(w+z+1)$, so that ${\rm Aut}(f)=\{z,-z-1\}$, and the product
part is:
$$
\left(1-\frac{w}{z}\right)\left(1-\frac{w}{-z-1}\right)=\left(\frac{z-w}{z}\right)\left(\frac{w+z+1}{z+1}\right)=
\frac{f(z)-f(w)}{f(z)}.
$$
Thus we get the representation:
$$
f(w)-f(z)=(-1)\cdot f(z)\left(1-\frac{w}{z}\right)\left(1-\frac{w}{-z-1}\right)=
$$
$$
=(-1)\cdot f(z)\left(1-\frac{w}{\phi_0(z)}\right)\left(1-\frac{w}{\phi_1(z)}\right).
$$
Now, let us consider a general polynomial: $f(z)=p_d(z)=a_dz^d+a_{d-1}z^{d-1}+\ldots+a_1z+a_0$, $a_d\ne 0$. Then
${\rm Aut}(p_d)=\{\phi_0(z),\ldots,\phi_{d-1}(z)\}$, where $\phi_0(z)=z$. Clearly:
$$
p_d(w)-p_d(z)=a_d\prod_{n=0}^{d-1}(w-\phi_n(z))=(-1)^da_d\left\{\prod_{n=0}^{n-1}\phi_n(z)\right\}
\prod_{n=0}^{d-1}\left(1-\frac{w}{\phi_n(z)}\right).
$$
By $p_d(w)-p_d(z)=a_d\prod_{n=0}^{d-1}(w-\phi_n(z))$ it follows that the free term of this $w$-polynomial is
given by $p_d(0)-p_d(z)=a_d\prod_{n=0}^{d-1}(0-\phi_n(z))=(-1)^da_d\prod_{n=0}^{d-1}\phi_n(z)$. So we proved 
that the Weierstrass factorization representation of the automorphic group of a general monic polynomial is:
$$
p_d(w)-p_d(z)=(p_d(0)-p_d(z))\prod_{n=0}^{d-1}\left(1-\frac{w}{\phi_n(z)}\right).
$$
A full generalization of the quadratic case. Is that formula valid for any entire function? Unfortunately it
is not the case. One might have falsely suspected at first that we can approximate an entire $f(z)$ by the
polynomials $p_d(z)$ which are the partial sums of the power series expansion of $f$. Each $p_d$ as the above 
simple Weierstrass factorization of $p_d(w)-p_d(z)$, and then when $d\rightarrow\infty$ we clearly have
$p_d(w)-p_d(z)\rightarrow f(w)-f(z)$. We might have hoped that the automorphic functions $\phi_n^{d}$ converge when
$d\rightarrow\infty$ to the automorphic functions $\phi_{0n}$ of $f$, and if we are lucky also
$$
\lim_{d\rightarrow\infty} \prod_{n=0}^{d-1}\left(1-\frac{w}{\phi_n^{d}(z)}\right)=
\prod_{n=0}^{\infty}\left(1-\frac{w}{\phi_{0n}(z)}\right),
$$
thus proving that:
$$
f(w)-f(z)=(f(0)-f(z))\prod_{n=0}^{\infty}\left(1-\frac{w}{\phi_{0n}(z)}\right).
$$
However, this clearly is wrong for the last infinite product is usually divergent unless we multiply
each term by the corresponding normalizing Weierstrass factor $\exp{(Q_{\lambda_n}(w/\phi_{0n}(z)))}$.
This simple formula has a chance of being correct only if $f$ is of order $0$ and $\forall\,n$, $\lambda_n=0$.
For the sake of completeness let us give a concrete example which proves that this simplistic formula
is wrong. If this formula were true for $f(z)=e^z$, we would have something like the following:
$$
e^w-e^z=(1-e^z)\prod_{n\in\mathbb{Z}}\left(1-\frac{w}{z+2\pi\cdot i\cdot n}\right)e^{w/(z+2\pi\cdot i\cdot n)}=
$$
$$
=(1-e^z)\left(1-\frac{w}{z}\right)e^{w/z}\prod_{n=1}^{\infty}\left(1-\frac{z^2-(z-w)^2}{z^2+4\pi^2n^2}\right)
e^{2zw/(z^2+4\pi^2n^2)}.
$$
If this was true then:
$$
\frac{e^w-e^z}{w-z}=(e^z-1)\frac{e^{w/z}}{z}\prod_{n=1}^{\infty}\left(1-\frac{z^2-(z-w)^2}{z^2+4\pi^2n^2}\right)
e^{2zw/(z^2+4\pi^2n^2)}.
$$
Taking the limits of both sides, when $w\rightarrow z$ we get:
$$
e^z=(e^z-1)\frac{e}{z}\prod_{n=1}^{\infty}\left(1-\frac{z^2}{z^2+4\pi^2n^2}\right)e^{2z^2/(z^2+4\pi^2n^2)}.
$$
We note that the convergent infinite product has no zero, as to be expected. Thus it looks promising, till
we specialize to $z=i\pi$:
$$
-1=\frac{-2}{i\pi}e\prod_{n=1}^{\infty}\left(1+\frac{\pi^2}{4\pi^2n^2-\pi^2}\right)e^{-2\pi^2/(4\pi^2n^2-\pi^2)}.
$$
That is nonsense, of course, because the left hand side is a real number while the right hand side is a
pure imaginary number! Can we fix this wrong? Let us denote the partial sums of the power series expansions
of $f(z)=e^z$ by:
$$
p_d(z)=\sum_{n=0}^d\frac{z^n}{n!}.
$$
Let us denote the automorphic functions of $p_d(z)$ by $\phi_n^d(z)$, $n=0,\ldots,d-1$. Then we proved that:
$$
p_d(w)-p_d(z)=(1-p_d(z))\prod_{n=0}^{d-1}\left(1-\frac{w}{\phi_n^d(z)}\right).
$$
The idea now, is to mimic at the polynomial level the form of the Weierstrass factorization of the limiting
function $e^w-e^z$. This means, that we multiply the factors by the Weierstrass normalizing factors. The result
is:
$$
p_d(w)-p_d(z)=(1-p_d(z))\exp{\left(-w\sum_{n=0}^{d-1}\left(\frac{1}{\phi_n^d(z)}\right)\right)}
\prod_{n=0}^{d-1}\left(1-\frac{w}{\phi_n^d(z)}\right)e^{w/\phi_n^d(z)}.
$$
At this point we take the limit $d\rightarrow\infty$ and assume that all the automorphic functions
of the partial sums converge to those of $e^z$ and that the finite normalized products of the $p_d$'s
converge to the Weierstrass canonical product of $e^w-e^z$. Here is what we get:
\begin{equation}\label{eq5}
e^w-e^z=(1-e^z)\left(1-\frac{w}{z}\right)e^{w/z}\exp{\left(-2zw\sum_{n=1}^{\infty}\left(\frac{1}{z^2+4\pi^2n^2}\right)\right)}\times
\end{equation}
$$
\times\prod_{n=1}^{\infty}\left(1-\frac{z^2-(z-w)^2}{z^2+4\pi^2n^2}\right)e^{2zw/(z^2+4\pi^2n^2)}.
$$
In other words this formula suggests the following identity using the notation of our Theorem \ref{thm1}:
\begin{equation}\label{eq6}
\exp{\left(F(w,e^z)\right)}=(1-e^z)\exp{\left(-2zw\sum_{n=1}^{\infty}\left(\frac{1}{z^2+4\pi^2n^2}\right)\right)}.
\end{equation}
If true then it is interesting because it is not clear why the infinite sum of equation (\ref{eq6}) is a holomorphic function
of $e^z$. We can evaluate this infinite sum. The following formula is well-known:
$$
2z\sum_{n=1}^{\infty}\left(\frac{1}{z^2-n^2}\right)=\pi\cot\pi z-\frac{1}{z}.
$$
We make use of it. We let $z=iu$ below.
$$
-2zw\sum_{n=1}^{\infty}\left(\frac{1}{z^2+4\pi^2n^2}\right)=-\left(\frac{z}{2\pi}\right)\left(\frac{w}{2\pi}\right)
\sum_{n=1}^{\infty}\left(\frac{1}{(z/2\pi)^2+n^2}\right)=
$$
$$
=i\left(\frac{w}{2\pi}\right)\cdot 2\left(\frac{u}{2\pi}\right)\sum_{n=1}^{\infty}\left(\frac{1}{(u/2\pi)^2-n^2}\right)=
i\left(\frac{w}{2}\right)\left\{\cot\left(\frac{u}{2}\right)-\left(\frac{2}{u}\right)\right\}=
$$
$$
=\left(\frac{w}{z}\right)-\left(\frac{w}{2}\right)\left(\frac{e^z+1}{e^z-1}\right).
$$
The element $w/z$ seems to be an obstacle for in order to make it an holomorphic function of $e^z$ we 
might write it as $w/\log e^z$, which, at least is not singular because $e^z\ne 1$. Plugging our result
into equation (\ref{eq5}) gives us finally the following interesting identity:
\begin{equation}\label{eq7}
e^w-e^z=e^{w/z}(1-e^z)\left(1-\frac{w}{z}\right)e^{w/z}\exp{\left(-\left(\frac{w}{2}\right)\left(\frac{e^z+1}{e^z-1}\right)\right)}\times
\end{equation}
$$
\times\prod_{n=1}^{\infty}\left(1-\frac{z^2-(z-w)^2}{z^2+4\pi^2n^2}\right)e^{2zw/(z^2+4\pi^2n^2)}.
$$
We recall that the last identity was derived using the idea outlined before equation (\ref{eq5}), namely
approximating the entire function $f(z)$ by a sequence of polynomials, the partial sums of its power
series expansion, using the identity we proved for polynomials:
$$
p_d(w)-p_d(z)=(p_d(0)-p_d(z))\prod_{n=0}^{d-1}\left(1-\frac{w}{\phi_n(z)}\right).
$$
Then multiplying the last identity by the Weierstrass normalization factors that correspond to the
Weierstrass expansion of $f$ and letting $d\rightarrow\infty$ assuming we have convergence of the
automorphic functions of the polynomials $p_d$ to the automorphic functions of $f$, and also convergence
of the finite products of the $p_d(w)-p_d(z)$ to the (generically) infinite product of $f(w)-f(z)$.
Remarkably all of that actually works! We now give an independent proof of the identity (\ref{eq7})
which does not rely on any of the above "convergences assumptions". Let us write our skeleton identity
using the variables $i\pi w$ and $i\pi z$ instead of $w$ and $z$:
$$
e^{i\pi w}-e^{i\pi z}=\exp{\left(g(i\pi w,i\pi z)\right)}\left(1-\frac{w}{z}\right)e^{w/z}\times
$$
$$
\times\prod_{n=1}^{\infty}\left(1+\frac{z^2-(z-w)^2}{4n^2-z^2}\right)\exp{\left(\frac{-2zw}{4n^2-z^2}\right)}
$$
Now we use the cotangent fractional series expansion to compute:
$$
\prod_{n=1}^{\infty}\exp{\left(\frac{-2zw}{4n^2-z^2}\right)}=\exp{\left(\left(\frac{\pi w}{2}\right)
\cot\left(\frac{\pi z}{2}\right)-\left(\frac{w}{z}\right)\right)}.
$$
Next we use the well known expansion:
$$
\frac{\pi z}{\sin(\pi z)}=\prod_{n=1}^{\infty}\left(\frac{n^2}{n^2-z^2}\right),
$$
to compute the infinite product:
$$
\prod_{n=1}^{\infty}\left(1+\frac{z^2-(z-w)^2}{4n^2-z^2}\right)=\prod_{n=1}^{\infty}\left(\frac{4n^2-(z-w)^2}{4n^2}\right)
\prod_{n=1}^{\infty}\left(\frac{4n^2}{4n^2-z^2}\right)=
$$
$$
=\left(\frac{z}{z-w}\right)\frac{\sin(\pi(z-w)/2)}{\sin(\pi z/2)}.
$$
Putting together the last three identities we proved gives:
$$
e^{i\pi w}-e^{i\pi z}=\exp{\left(g(i\pi w,i\pi z)\right)}e^{w/z}\frac{\sin(\pi(z-w)/2)}{\sin(\pi z/2)}\exp{\left(\left(\frac{\pi w}{2}\right)
\cot\left(\frac{\pi z}{2}\right)-\left(\frac{w}{z}\right)\right)}.
$$
We solve for $\exp{\left(g(i\pi w,i\pi z)\right)}$ and replace $i\pi w$, $i\pi z$ by $w$ and $z$ respectively. This gives:
$$
\exp{\left(g(w,z)\right)}=(e^w-e^z)\frac{\sin(z/2i)}{\sin((z-w)/2i)}\exp{\left(-\frac{w}{2i}\cot\left(\frac{z}{2i}\right)\right)}.
$$
This concludes the proof of identity (\ref{eq7}). 

Using Proposition \ref{prop1} we deduce that if $(e^w-1)(e^z-1)\ne 0$, then:
$$
e^{w/z}(1-e^z)\left(1-\frac{w}{z}\right)e^{w/z}\exp{\left(-\left(\frac{w}{2}\right)\left(\frac{e^z+1}{e^z-1}\right)\right)}\times
$$
$$
\times\prod_{n=1}^{\infty}\left(1-\frac{z^2-(z-w)^2}{z^2+4\pi^2n^2}\right)e^{2zw/(z^2+4\pi^2n^2)}=
$$
$$
=-e^{z/w}(1-e^w)\left(1-\frac{z}{w}\right)e^{z/w}\exp{\left(-\left(\frac{z}{2}\right)\left(\frac{e^w+1}{e^w-1}\right)\right)}\times
$$
$$
\times\prod_{n=1}^{\infty}\left(1-\frac{w^2-(w-z)^2}{w^2+4\pi^2n^2}\right)e^{2wz/(w^2+4\pi^2n^2)}.
$$
It is interesting to note that both sides are entire in $w$ (left) and in $z$ (right). That agrees with the Gronwall-Hahn
Theorem. The left side is clearly $z$-holomorphic in $z\in\mathbb{C}-2\pi i\mathbb{Z}$, and the right side is
$w$-holomorphic in $w\in\mathbb{C}-2\pi i\mathbb{Z}$. Thus both sides are entire in $(w,z)$. The essential singularities
of
$$
\exp{\left(-\left(\frac{w}{2}\right)\left(\frac{e^z+1}{e^z-1}\right)\right)},
$$
and of $e^{w/z}$ are somehow canceled out by the infinite product.

We end this section by pointing at two findings that seem to emerge from our computations. The first is the extent
to which an entire function $f(z)$ is determined by a partial knowledge of its fibers. The notion of the fiber is very close
to the notion of the automorphic group, namely $\forall\,w\in\mathbb{C}$ the fiber $f^{-1}(w)=\{z_j\,|\,f(z_j)=w\}$ is the 
discrete subset of $\mathbb{C}$ (we assume that $f$ is non-constant) of all the $f$-pre-images of $w$. We note that
if $z_0\in f^{-1}(w)$, then $f^{-1}(w)$ is simply the ${\rm Aut}(f)$-orbit of $z_0$, i.e. we have the
identity $f^{-1}(w)=\{\phi(z_0)\,|\,\phi\in{\rm Aut}(f)\}$. For a general function (not necessarily holomorphic or
even continuous) the knowledge of the pairs $(w,f^{-1}(w))$ determines $f$ (uniquely). The mere knowledge of all the
fibers $f^{-1}(w)$, without knowing the $w$ itself clearly does not determine $f$. This is very close to knowing
the automorphic group of $f$ (for that partitions $\mathbb{C}$ into the $f$-fibers without the knowledge of the $w$).
So far for general functions $f$. Even if we know in advance that $f$ is continuous, the automorphic group, i.e. the
fibers $f^{-1}(w)$ do not determine $f$. If $G$ is a continuous injection then $f$ and $G\circ f$ have identical
family of fibers. But our case is very different from the continuous case. Our functions $f$ are entire and hence 
are rigid. The case that shows how this holomorphic rigidity makes the difference is the case of polynomials. we already
noted that if $p(z)=a_dz^d+\ldots+a_0$, $a_d\ne 0$, and if ${\rm Aut}(p)=\{\phi_0(z),\ldots,\phi_{d-1}(z)\}$, then
$$
p(z)=p(0)+(-1)^{d+1}a_d\prod_{j=0}^{d-1}\phi_j(z).
$$
Thus the product $\phi_0\cdot\ldots\cdot\phi_{d-1}$, i.e. the product of the $p$-fiber determines the function $p(z)$
up to a multiplicative constant $a_d$ different from $0$ and an additive constant $p(0)$. Thus we do not have to
know the fiber, just the product of its elements, a very partial information indeed, suffice to essentially reconstruct
the function.

The second finding is closely related to the first one, but here we want to handle entire not necessarily polynomials.
since in this case the group ${\rm Aut}(f)$ is usually infinite, it does not make sense to multiply its elements.
Thus in this more complicated situation we ask the following: Given ${\rm Aut}(f)$ where $f$ is a non-constant
entire function, can we reconstruct the function $f$ (up to minor parameters)? The way we outlined how to handle
the case $f(z)=e^z$ might give us the way to solve this problem.

\begin{theorem}\label{thm2}
If $f(z)-f(0)\ne 0$, then there is a function $g(w,z)$, entire in $w$ and there are non-negative integers
$\lambda_n=\lambda_n(w,\phi_{0n}(z))$ such that:
$$
f(z)=f(w)-\exp{\left(g(w,z)\right)}\prod_{n=1}^{\infty}\left(1-\frac{w}{\phi_{0n}(z)}\right)e^{Q_{\lambda_n}(w/\phi_{0n}(z))}.
$$
Here as usual:
$$
Q_{\lambda}\left(\frac{w}{\phi(z)}\right)=\left(\frac{w}{\phi(z)}\right)+\frac{1}{2}\left(\frac{w}{\phi(z)}\right)^2+\ldots
+\frac{1}{\lambda}\left(\frac{w}{\phi(z)}\right)^{\lambda}.
$$
Here the non-negative integers $\lambda_n$ are chosen so that the infinite product converges on uniformly on compact
subsets of $\mathbb{C}$. For example we might take the canonical product of the automorphic group $\{\phi_{0n}\}$.
If we can sum up by some summation method the infinite series:
$$
\sum_{n=1}^{\infty}Q_{\lambda_n}\left(\frac{w}{\phi_{0n}(z)}\right),
$$
to give an holomorphic sum, and that problem is at the moment an open problem, then:
$$
\exp{\left(g(w,z)\right)}=(f(0)-f(z))\exp{\left(-\sum_{n=1}^{\infty}Q_{\lambda_n}\left(\frac{w}{\phi_{0n}(z)}\right)\right)}.
$$
In this case we can reconstruct $f(z)$ from ${\rm Aut}(f)$ by the formula:
$$
f(z)=\frac{f(0)\cdot L-f(w)}{L-1},
$$
where
$$
L=\exp{\left(-\sum_{n=1}^{\infty}Q_{\lambda_n}\left(\frac{w}{\phi_{0n}(z)}\right)\right)}\cdot
\prod_{n=1}^{\infty}\left(1-\frac{w}{\phi_{0n}(z)}\right)e^{Q_{\lambda_n}(w/\phi_{0n}(z))},
$$
and where $w\in\mathbb{C}$ can be any number for which $L\ne 1$.
\end{theorem}

\section{Reconstruction formulas for $f(z)$ and for $f'(z)$ in terms of approximating automorphic functions. Relations
between the groups ${\rm Aut}(f)$ and ${\rm Aut}_z(g(w,z))$}\label{sec7}

The difference between the very simple reconstruction of a polynomial $p(z)=a_dz^d+\ldots+a_0$, $a_d\ne 0$ from its automorphic group
${\rm Aut}(p)=\{\phi_j(z)\,|\,j=0,\ldots,d-1\}$, $p(z)=p(0)+(-1)^{d+1}a_d\prod_{j=0}^{d-1}\phi_j(z)$ on the one hand,
and the reconstruction of a general non-constant entire function $f(z)$, in Theorem \ref{thm2}, on the other hand gives
the feeling of a possibility of a simpler reconstruction (in the general entire case). Indeed we can point at such
a formula, seemingly simpler than the one in Theorem \ref{thm2}. However, the hidden complication is in the approximating
procedure within that formula.

The setting is that we have an entire non-constant function $f(z)$ represented in terms of its Maclaurin's series:
$f(z)=\sum_{n=0}^{\infty}a_nz^n$, $\limsup_{n\rightarrow\infty}|a_n|^{1/n}=0$. We denote the sequence of the partial
sums by $f_n(z)=\sum_{k=0}^na_kz^k$. Practically we consider those partial sums for which (in the notation above) $a_n\ne 0$.
We denote the automorphic groups: ${\rm Aut}(f_n)=\{\phi_0^{(n)}(z),\ldots,\phi_{n-1}^{(n)}(z)\}$. These are all the solutions
of the automorphic equation: $f_n(\phi_j^{(n)}(z))=f_n(z)$, $j=0,\ldots,n-1$. Thus $f_n(w)-f_n(z)=a_n\prod_{j=0}^{n-1}(w-
\phi_j^{(n)}(z))$. The automorphic functions of $f_n(z)$ satisfy the Vieta identities:
$$
(-1)^ka_n\sum_{0\le i_1<i_2<\ldots<i_k\le n-1}\prod_{j=1}^k\phi_{i_j}^{(n)}(z)=\left\{\begin{array}{ccc} a_{n-k} & , & k<n \\
a_0-f_n(z) & , & k=n \end{array}\right..
$$
As usual we denote ${\rm Aut}(f)=\{\phi_0(z),\phi_1(z),\phi_2(z),\ldots\}$, and we note that:
$$
f(w)-f_n(w)=\sum_{k=n+1}^{\infty}a_kw^k,\,\,\,\,\,f(z)-f_n(z)=\sum_{k=n+1}^{\infty}a_kz^k.
$$
Hence:
$$
f(w)-f(z)=(f(w)-f_n(w))+(f_n(w)-f_n(z))-(f(z)-f_n(z))=
$$
$$
=\sum_{k=n+1}^{\infty}a_k(w^k-z^k)+a_n\prod_{j=0}^{n-1}(w-\phi_j^{(n)}(z)).
$$
If we fix $R>0$, then $\lim_{n\rightarrow\infty}\sum_{k=n+1}^{\infty}a_k(w^k-z^k)=0$ uniformly in $w, z$.
So that $f(w)-f(z)=\lim_{n\rightarrow\infty}a_n\prod_{j=0}^{n-1}(w-\phi_j^{(n)}(z))$ uniformly on compact
subsets of $\mathbb{C}$. Thus it is straightforward to prove that we can divide by $(w-z)$ and take $w\rightarrow z$
and obtain: $f'(z)=\lim_{n\rightarrow\infty}a_n\prod_{j=1}^{n-1}(z-\phi_j^{(n)}(z))$. Finally let us
fix $N$ and replace the first $N$ factors $\prod_{j=1}^{N-1}(z-\phi_j^{(n)}(z))$ when $n\rightarrow\infty$
by $\prod_{j=1}^{N-1}(z-\phi_j(z))$. We obtain:
$$
\frac{f(w)-f(z)}{a_N\prod_{j=1}^{N-1}(z-\phi_j(z))}=\lim_{n\rightarrow\infty}\frac{a_n}{a_N}
\prod_{j=N}^{n}(z-\phi_j^{(n)}(z)).
$$
This proves the following:
\begin{proposition}\label{prop4}
We have $f(w)-f(z)=\lim_{n\rightarrow\infty}a_n\prod_{j=0}^{n-1}(w-\phi_j^{(n)}(z))$ uniformly on compact
subsets of $\mathbb{C}$. Also $f'(z)=\lim_{n\rightarrow\infty}a_n\prod_{j=1}^{n-1}(z-\phi_j^{(n)}(z))$ uniformly
on compact subsets of $\mathbb{C}$ and likewise for any fixed $N$:
$$
\frac{f(w)-f(z)}{a_N\prod_{j=1}^{N-1}(z-\phi_j(z))}=\lim_{n\rightarrow\infty}\frac{a_n}{a_N}
\prod_{j=N}^{n}(z-\phi_j^{(n)}(z)).
$$
\end{proposition}

\begin{remark}\label{rem17}
The reconstruction formulas given in Proposition \ref{prop4}:
$$
f(z)=f(w)-\lim_{n\rightarrow\infty}a_n\prod_{j=0}^{n-1}(w-\phi_j^{(n)}(z)),
$$
and
$$
f'(z)=\lim_{n\rightarrow\infty}a_n\prod_{j=1}^{n-1}(z-\phi_j^{(n)}(z)),
$$
seem to be simpler than that in Theorem \ref{thm2}, but the cost lies in the sequential
limit $\lim_{n\rightarrow\infty}$, which describes an auxiliary approximation of the
entire functions by polynomials (the partial sums $f_n$).
\end{remark}
\noindent
We now bring few properties of the automorphic functions of a non-constant entire function $f$ of
a finite order. These are related to the results in section \ref{sec4}. First it will be convenient
to rephrase Proposition \ref{prop3}:

\begin{corollary}\label{cor9}
If $f$ is a non-constant entire function of a finite order, then the automorphic group of $f(z)$ is a 
subgroup of the automorphic group of $\exp{(g(w,z))}$ as a function of $z$ (for any fixed $w\in\mathbb{C}$).
In symbols: ${\rm Aut}(f(z))\subseteq{\rm Aut}_z(\exp{(g(w,z)})$, $\forall\,w\in\mathbb{C}$.
\end{corollary}
\noindent
We noticed in Corollary \ref{cor055} that any finite order element of ${\rm Aut}(f)$ ($f$ a non-constant
entire of a finite order) belongs to $\ker (T)$. Our next result gives a non-trivial characterization
of the elements in $\ker (T)$.

\begin{corollary}\label{10}
Let $f(z)$ be a non-constant entire function of a finite order. Then: $\ker (T)\equiv{\rm Aut}_z(g(w,z))$,
for any fixed $w\in\mathbb{C}$.
\end{corollary}
\noindent
{\bf Proof.} \\
The automorphic function $\Phi\in{\rm Aut}(f(z))$ belongs to $\ker (T)$, where the homomorphism $T$
was defined in Corollary \ref{cor05} $\Longleftrightarrow$ $T(\Phi)=\frac{1}{2\pi\cdot i}\left(g(w,\Phi(z))-g(w,z)\right)=0$
$\Longleftrightarrow$ $g(w,\Phi(z))=g(w,z)$, $\forall\,w\in\mathbb{C}$ $\Longleftrightarrow$
$\Phi\in{\rm Aut}_z(g(w,z))$, $\forall\,w\in\mathbb{C}$. $\qed $

\begin{corollary}\label{cor11}
Let $f$ be a non-constant entire function of a finite order. Then: \\
(a) Any $\Phi\in{\rm Aut}(f(z))-{\rm Aut}_z(g(w,z))$, is of an infinite order. \\
(b) Let $\Phi_{i_0j_0}\in{\rm Aut}(f(z))$ be such that $|T(\Phi_{i0j_0})|=\min\{|T(\Phi)|\,|\,
\Phi\in{\rm Aut}(f(z))-\ker(T)\}$. Then $T({\rm Aut}(f(z)))=<\Phi_{i_0j_0}>$. We have $\forall\,
k\in\mathbb{Z}$, $T(\Phi_{i_0j_0}^{\circ k})=k\cdot T(\Phi_{i_0j_0})$. If ${\rm Aut}(f(z))-\ker(T)
=\emptyset$, we agree to define $<\Phi_{i_0j_0}>=\{0\}$. \\
(c) If $T(\Phi_{ij})=T(\Phi_{\alpha\beta})$ then $\Phi_{ij}\circ\Phi_{\alpha\beta}^{-1}\in{\rm Aut}_z(g(w,z))$.
This could be written as $g(w,\Phi_{ij}\circ\Phi_{\alpha\beta}^{-1}(z))=g(w,z)$ or, equivalently as
$g(w,\Phi_{ij}(z))=g(w,\Phi_{\alpha\beta}(z))$.
\end{corollary}

\begin{remark}\label{rem18}
The automorphic group of a non-constant entire function of any order can contain elements of a finite
order and elements of infinite order. For example if ${\rm Aut}(f(z))$ contains elements of infinite
order, then so does the group ${\rm Aut}(f(z^2))$ but this last group contains also the order $2$ element
$\Phi(z)=-z$.
\end{remark}

\begin{remark}\label{rem19}
We point out that the construction of a Dirichlet fundamental domain for Fuchsian groups could be
used to construct a fundamental domain (maximal domain) for a non-constant entire function. Let
$f(z)$ be a non-constant entire function, $z_0\in\mathbb{C}$ a regular point of $f$ (i.e. $f'(z_0)\ne 0$)
and $\rho(\cdot,\cdot)$ a metric on $\mathbb{C}$. Mostly we have in mind the $f$-path metric, $\rho_f$ that is
induced by $f$. We recall what that is: let $z, w\in\mathbb{C}$ and let $\gamma:\,[0,1]\rightarrow\mathbb{C}$
be a continuous path from $z$ to $w$. Thus $\gamma(0)=z$ and $\gamma(1)=w$. Then the length of $\gamma$ is
given by the standard length of the $f$-image path $f\circ\gamma$. We will denote this length by $l_f(\gamma)$.
$$
l_f(\gamma)=\int_0^1\left|f'\left(\gamma(t)\right)\right|\left|\gamma'(t)\right|dt.
$$
The $f$-path metric is given by: $\rho_f(z,w)=\inf_{\gamma} l_f(\gamma)$, where the infimum is taken over all
the piecewise differentiable paths $\gamma$ from $z$ to $w$. The Dirichlet fundamental domain of $f(z)$,
centered at $z_0$, with respect to the $f$-path metric is:
$$
\{z\in\mathbb{C}\,|\,\rho_f(z,z_0)\le\rho_f(\Phi_{ij}(z),z_0)\,\,\forall\,\Phi_{ij}\in{\rm Aut}(f)\}.
$$
An alternative way to define that, which avoids using the notion of the automorphic group of $f$ is as follows:
$$
\{z\in\mathbb{C}\,|\,\rho_f(z,z_0)\le\rho_f(w,z_0)\,\,\forall\,w\in f^{-1}(f(z))\}.
$$
The interior of the set above is a domain (an open connected subset of $\mathbb{C}$), and the function $f$ is
one-to-one in this domain and the domain is maximal with respect to this property (of $f$ being injective).
\end{remark}

\section{The function $g(w,z)-g(0,z)$ is determined by the negative moments of the elements in ${\rm Aut}(f(z))$}\label{sec8}

\begin{theorem}\label{thm3}
Let $f$ be a non-constant entire function. Let us denote ${\rm Aut}(f)=\{\phi_{0n}(z)\,|\,n=0,1,2,\ldots\}$ 
($\phi_{00}\equiv {\rm id.}$) and let $g(w,z)$ be the function in the exponential of the Weierstrass (canonical)
factorization of $f(w)-f(z)$. $g(w,z)$ is entire in $w$ and holomorphic in $z\not\in f^{-1}(f(0))$. Then for
$k=1,2,3,\ldots $ we have the identities:
$$
\frac{1}{k!}\frac{\partial^k g}{\partial w^k}(0,z)=-\frac{1}{k}\sum_{\left\{\begin{array}{l} n \\ \lambda_n\ge k\end{array}\right.}
\left(\frac{1}{\phi_{0n}(z)}\right)^k.
$$
The left hand side is the $k+1$'st Maclaurin coefficient in the expansion of $g(w,z)-g(0,z)$. The right hand
side is $-1/k$ multiplying the $-k$-moment of all the relevant automorphic functions of $f$. That explains the
title of this section. Like in Theorem \ref{thm2} we assume that we have some summation method the infinite series:
$$
\sum_{n=1}^{\infty}Q_{\lambda_n}\left(\frac{w}{\phi_{0n}(z)}\right).
$$
\end{theorem}

\begin{remark}\label{rem20}
Thus the proof below is based on a vague summability assumption! It is mostly supported by formal computational
steps, and not justified. However, this already interesting sketch points to the fact that a true proof if exists will
not be an easy one, and probably it will have to utilize summability theoretical arguments.
\end{remark}
\noindent
{\bf Proof.} \\
By Theorem \ref{thm2} we have the following identity:
$$
\exp{\left(g(w,z)\right)}=(f(0)-f(z))\exp{\left(-\sum_{n=1}^{\infty}Q_{\lambda_n}\left(\frac{w}{\phi_{0n}(z)}\right)\right)}.
$$
It is assuming that we have a summability method for the infinite series:
$$
\sum_{n=1}^{\infty}Q_{\lambda_n}\left(\frac{w}{\phi_{0n}(z)}\right),
$$
which results in a holomorphic function. Here, as usual:
$$
Q_{\lambda_n}\left(\frac{w}{\phi_{0n}(z)}\right)=\left(\frac{w}{\phi_{0n}(z)}\right)+\frac{1}{2}\left(\frac{w}{\phi_{0n}(z)}\right)^2+
\ldots+\frac{1}{\lambda_n}\left(\frac{w}{\phi_{0n}(z)}\right)^{\lambda_n}.
$$
When we substitute $w=0$ into the identity above, we obtain $\exp{(g(0,z))}=f(0)-f(z)$. So we can rewrite our identity
as follows:
$$
\exp{\left(g(w,z)-g(0,z)\right)}=\exp{\left(-\sum_{n=1}^{\infty}Q_{\lambda_n}\left(\frac{w}{\phi_{0n}(z)}\right)\right)}.
$$
We deduce that there is an integer $N\in\mathbb{Z}$, such that:
$$
g(w,z)-g(0,z)=2\pi\cdot i\cdot N-\sum_{n=1}^{\infty}Q_{\lambda_n}\left(\frac{w}{\phi_{0n}(z)}\right).
$$
If we plug in $w=0$, we obtain $0=2\pi\cdot i\cdot N$, so that $N=0$ and we have:
$$
g(w,z)-g(0,z)=-\sum_{n=1}^{\infty}Q_{\lambda_n}\left(\frac{w}{\phi_{0n}(z)}\right).
$$
On the other hand $g(w,z)$ is an entire function in $w$ and so it has a power series expansion that converges
in the whole $w$-plane:
$$
g(w,z)-g(0,z)=\sum_{n=1}^{\infty}\frac{1}{n!}\frac{\partial^n g}{\partial w^n}(0,z)\cdot w^n.
$$
Hence we conclude that we have the following identity:
$$
\sum_{n=1}^{\infty}\frac{1}{n!}\frac{\partial^n g}{\partial w^n}(0,z)\cdot w^n=
-\sum_{n=1}^{\infty}Q_{\lambda_n}\left(\frac{w}{\phi_{0n}(z)}\right).
$$
Let us write the series on the right hand side, as power series in $w$. We recall that:
$$
Q_{\lambda_n}\left(\frac{w}{\phi_{0n}(z)}\right)=\left(\frac{w}{\phi_{0n}(z)}\right)+\frac{1}{2}\left(\frac{w}{\phi_{0n}(z)}\right)^2+
\ldots+\frac{1}{\lambda_n}\left(\frac{w}{\phi_{0n}(z)}\right)^{\lambda_n},
$$
when $\lambda_n\in\mathbb{Z}^+$. Otherwise $Q_0(u)\equiv 0$. We would like to compute the coefficient
of $w^k$ on the right hand side of our identity. It is the sum of all the elements of the following form:
$$
\frac{1}{k}\left(\frac{w}{\phi_{0n}(z)}\right)^k,
$$
provided, of course, that the condition $\lambda_n\ge k$ is fulfilled. Thus we obtain the following identity:
$$
-\sum_{n=1}^{\infty}Q_{\lambda_n}\left(\frac{w}{\phi_{0n}(z)}\right)=
-\sum_{k=1}^{\infty}\frac{1}{k}\sum_{\left\{\begin{array}{l} n \\ \lambda_n\ge k\end{array}\right.}
\left(\frac{1}{\phi_{0n}(z)}\right)^k\cdot w^k.
$$
By the uniqueness of the coefficients in a power series we conclude that for $k=1,2,3,\ldots$ we have:
$$
\frac{1}{k!}\frac{\partial^k g}{\partial w^k}(0,z)=-\frac{1}{k}\sum_{\left\{\begin{array}{l} n \\ \lambda_n\ge k\end{array}\right.}
\left(\frac{1}{\phi_{0n}(z)}\right)^k.
$$
Indeed we note that the sum on the right side of the last identity is the $(-k)$-moment of $\phi_{0n}(z)$, where, of
course, $\lambda_n\ge k$. $\qed $ \\

\section{An infinite product representation of $f'(w)$}\label{sec9}

In this section we will explore how the Weierstrass factorization of a non-constant entire function
induces an infinite product representation on its derivative. This will reveal a connection between
the zeros of the derivatives and the fixed-points of the elements of the automorphic group of the
function. We begin with a peculiar division property of entire functions and of their derivatives
which follows by a composition relation between these functions.

\begin{proposition}\label{prop5}
Let $f(z)$ and $g(z)$ be two non-constant entire functions. If ${\rm Aut}(f(z))\subseteq{\rm Aut}(g(z))$, then
$f'(z)$ divides $g'(z)$ over the algebra of entire functions, i.e. $\exists\,H(z)$, an entire function such
that $g'(z)=H(z)\cdot f'(z)$. In fact $\exists\,G(w,z)$, entire in $w$ and holomorphic in
$z\in\mathbb{C}-g^{-1}(g(0))-f^{-1}(f(0))$ such that $g(w)-g(z)=(f(w)-f(z))\cdot G(w,z)$.
\end{proposition}
\noindent
{\bf Proof.} \\
By Lemma \ref{lem1} it follows that there exists a function $h(w)$, holomorphic on $f(\mathbb{C})$ such
that $g(z)=h(f(z))$. Hence $g'(z)=h'(f(z))\cdot f'(z)$ which proves the first assertion with the entire
function $H(z)=h'(f(z))$. Next we denote ${\rm Aut}(f(z))=\{\phi_{0n}(z)\}$, and ${\rm Aut}(g(z))=\{\psi_{0n}(z)\}$.
By the Weierstrass factorization theorem we have:
$$
g(w)-g(z)=e^{L(w,z)}\prod_n\left(1-\frac{w}{\psi_{0n}(z)}\right)e^{Q_{\delta_n}(w/\psi_{0n}(z))},\,\,\,g(0)-g(z)\ne 0,
$$
$$
f(w)-f(z)=e^{l(w,z)}\prod_n\left(1-\frac{w}{\phi_{0n}(z)}\right)e^{Q_{\lambda_n}(w/\phi_{0n}(z))},\,\,\,f(0)-f(z)\ne 0,
$$
where $\{\phi_{0m}(z)\}\subseteq\{\psi_{0m}(z)\}$. Hence:
$$
g(w)-g(z)=(f(w)-f(z))e^{L(w,z)-l(w,z)}\prod_{\psi_{0n}\not\in\{\phi_{0m}\}}\left(1-\frac{w}{\psi_{0n}(z)}\right)
e^{Q_{\delta_n}(w/\psi_{0n}(z))}\times
$$
$$
\times\prod_n\left(1-\frac{w}{\phi_{0n}(z)}\right)e^{(Q_{\delta_n}-Q_{\lambda_n})(w/\phi_{0n}(z))}.
$$
\noindent
$\qed $ \\
\\
We recall the reconstruction formula given in Proposition \ref{prop4} for the derivative:
$$
f'(z)=\lim_{n\rightarrow\infty}a_n\prod_{j=1}^{n-1}(z-\phi_j^{(n)}(z)),
$$
This formula suggests a possible relation between the zeros of $f'(z)$ and the fixed points
of the automorphic functions of $f(z)$. Indeed we will prove that this is the case. 

Let $f(w)$ be a non-constant entire function. Let our $z$-parameter space be $\mathbb{C}-f^{-1}(f(0))$.
We consider the Weierstrass factorization of $f(w)-f(z)$ as an entire function of $w$. By our choice
of the parameter $z$, we have $f(0)-f(z)\ne 0$. Thus $0\not\in Z_w(f(w)-f(z))$ and hence:
$$
f(w)-f(z)=e^{g(w,z)}\prod_n\left(1-\frac{w}{\phi_{0n}(z)}\right)e^{Q_{\lambda_n}(w/\phi_{0n}(z))}.
$$
Here $g(w,z)$ is entire in $w$ and holomorphic for any $z\in\mathbb{C}-f^{-1}(f(0))$. Also
${\rm Aut}(f(w))=\{\phi_{0n}(w)\}_n$, and we agree that $\phi_{00}(w)\equiv w$. The numbers
$\lambda_n\in\mathbb{Z}^+\cup\{0\}$ and:
$$
Q_{\lambda}(u)=\left\{\begin{array}{lll} u+u^2/2+\ldots+u^{\lambda}/\lambda &, & \lambda\in\mathbb{Z}^+ \\
0 & , & \lambda=0 \end{array}\right.,\,\,\,\,\lambda_0=0.
$$
So:
$$
f(w)-f(z)=\left(1-\frac{w}{z}\right)e^{g(w,z)}\prod_{n\ne 0}\left(1-\frac{w}{\phi_{0n}(z)}\right)
e^{Q_{\lambda_n}(w/\phi_{0n}(z))}=
$$
$$
=\left(\frac{z-w}{z}\right)e^{g(w,z)}\prod_{n\ne 0}\left(1-\frac{w}{\phi_{0n}(z)}\right)
e^{Q_{\lambda_n}(w/\phi_{0n}(z))}.
$$
Hence, assuming that $w\ne z$, we obtain:
$$
\frac{f(w)-f(z)}{w-z}=-\frac{1}{z}e^{g(w,z)}\prod_{n\ne 0}\left(1-\frac{w}{\phi_{0n}(z)}\right)
e^{Q_{\lambda_n}(w/\phi_{0n}(z))}.
$$
Assuming that $w\in\mathbb{C}-f^{-1}(f(0))$, and that $z$ and $w$ are close enough so that in the
limit process $z\rightarrow w$, the non-negative integers $\lambda_n$ do not change, we get:
$$
f'(w)=\lim_{\begin{array}{l} z\rightarrow w \\ z\not\in f^{-1}(f(0))\end{array}}\frac{f(w)-f(z)}{w-z}=
-\frac{1}{w}e^{g(w,w)}\prod_{n\ne 0}\left(1-\frac{w}{\phi_{0n}(w)}\right)e^{Q_{\lambda_n}(w/\phi_{0n}(w))}.
$$
This shows that the zero set of the derivative function, $Z(f'(w))$ originates in three possible locations: \\
a) The fiber $f^{-1}(f(0))$might contain zeros of $f'(w)$. \\
b) Any fixed-point $w$ of a (non-identity) automorphic function $\phi_{0n}$, must be a zero of $f'(w)$. Thus:
$$
\{w\in\mathbb{C}-f^{-1}(f(0))\,|\,\exists\,n\ne 0,\,\phi_{0n}(w)=w\}\subseteq Z(f').
$$
c) The zeros (if any) of the functions $e^{Q_{\lambda_n}(w/\phi_{0n}(w)}$ for $n\ne 0$ and off
the fiber $f^{-1}(f(0))$. Thus:
$$
\bigcup_{n\ne 0}Z\left(e^{Q_{\lambda_n}(w/\phi_{0n}(w)}\right)\cap\left(\mathbb{C}-f^{-1}(f(0))\right)
\subseteq Z(f').
$$
Let us look at any automorphic equation in the domain of the definition of the corresponding automorphic
function: $f(\phi_{0n}(w))=f(w)$. By differentiation (assuming that $\phi_{0n}(w)$ has a derivative there):
$\phi_{0n}'(w)\cdot f'(\phi_{0n}(w))=f'(w)$. This implies that if $f'(w)=0$, then either $\phi_{0n}'(w)=0$
or $f'(\phi_{0n}(w))=0$. If $w$ is of type b, i.e. a fixed-point of the above automorphic function,
$\phi_{0n}(w)=w$, then clearly $f'(\phi_{0n}(w))=f'(w)=0$ and a consideration of the order of this zero
of $f'$ implies that $\phi_{0n}'(w)\ne 0$. Thus $w$ is a regular point of the automorphic function.

\begin{remark}\label{rem21}
If $f'(\phi_{0n}(w))=f'(w)=0$, then $\phi_{0n}'(w)\ne 0$ and also for its inverse $(\phi_{0n}^{-1})'(\phi_{0n}(w))\ne 0$.
\end{remark}
\noindent
If $f'(\phi_{0n}(w))\ne 0$, then necessarily $\phi_{0n}'(w)=0$. In this case (assuming it is not type b) we
either have $f(w)=f(\phi_{0n}(w))=f(0)$ (type a) or $e^{Q_{\lambda_n}(w/\phi_{0n}(w))}=0$. \\
\underline{What are the type c points?} These are zeros of $e^{Q_{\lambda_n}(w/\phi_{0n}(w))}$ outside
the fiber $f^{-1}(f(0))$. This implies that $\lambda_n>0$ and that locally the function:
$$
Q_{\lambda_n}\left(\frac{w}{\phi_{0n}(w)}\right)=\left(\frac{w}{\phi_{0n}(w)}\right)+\frac{1}{2}
\left(\frac{w}{\phi_{0n}(w)}\right)^2+\ldots+\frac{1}{\lambda_n}\left(\frac{w}{\phi_{0n}(w)}\right)^{\lambda_n},
$$
is the logarithm of a function that vanishes at $w$. Hence:
$$
Q_{\lambda_n}\left(\frac{w}{\phi_{0n}(w)}\right)=\left(\frac{w}{\phi_{0n}(w)}\right)+\frac{1}{2}
\left(\frac{w}{\phi_{0n}(w)}\right)^2+\ldots+\frac{1}{\lambda_n}\left(\frac{w}{\phi_{0n}(w)}\right)^{\lambda_n}=
$$
$$
=\log((w-w_0)\cdot h(w)).
$$
But this implies that:
$$
\lim_{w\rightarrow w_0}\left|Q_{\lambda_n}\left(\frac{w}{\phi_{0n}(w)}\right)=\left(\frac{w}{\phi_{0n}(w)}\right)+\frac{1}{2}
\left(\frac{w}{\phi_{0n}(w)}\right)^2+\ldots+\frac{1}{\lambda_n}\left(\frac{w}{\phi_{0n}(w)}\right)^{\lambda_n}\right|=+\infty.
$$
So $\lim_{w\rightarrow w_0}\phi_{0n}(w)=0$ and hence $\phi_{0n}(w_0)=0$, which implies that the point 
$w_0\in f^{-1}(f(0))$ in the forbidden fiber $f^{-1}(f(0))$. So type c points do not exist. We completed the
proof of the following:

\begin{lemma}\label{lem3}
If $f(w)$ is a non-constant entire function, then:
$$
Z(f')=\bigcup_{\phi_{0n}\in{\rm Aut}(f)-\{id\}}\{w\in\mathbb{C}-f^{-1}(f(0))\,|\,\phi_{0n}(w)=w\}\cup
\left(Z(f')\cap f^{-1}(f(0))\right).
$$
\end{lemma}
\noindent
We can now prove the simple relation that exists between the zeros of the derivative and the fixed-point
of the automorphic function.

\begin{theorem}\label{thm4}
Let $f(w)$ be a non-constant entire function. Then:
$$
Z(f')=\bigcup_{\phi_{0n}\in{\rm Aut}(f)-\{id\}}\{w\in\mathbb{C}\,|\,\phi_{0n}(w)=w\}:={\rm Fix}({\rm Aut}(f)).
$$
\end{theorem}
\noindent
{\bf Proof.} \\
Let $t\in\mathbb{C}$ be any complex number. It is clear that the number $0$ that appears in the formula
of Lemma \ref{lem3} in $f^{-1}(f(0))$ has no special significance. Indeed we could have expanded $f(w)-f(z)$
in a Weierstrass product centered at $t$ instead of $0$ and obtain in Lemma \ref{lem3} the $t$'th version:
$$
Z(f')=\bigcup_{\phi_{0n}\in{\rm Aut}(f)-\{id\}}\{w\in\mathbb{C}-f^{-1}(f(t))\,|\,\phi_{0n}(w)=w\}\cup
\left(Z(f')\cap f^{-1}(f(t))\right).
$$
The fibers $f^{-1}(f(t))$ are discrete subsets of $\mathbb{C}$ for any such a $t\in\mathbb{C}$. Even more is
true, namely: $t_1\ne t_2\,\,\Leftrightarrow\,\,f^{-1}(f(t_1))\cap f^{-1}(f(t_2))=\emptyset$. Thus the claim
of our theorem follows. $\qed $ \\

\begin{theorem}\label{thm5}
If $f$ is a non-constant entire function with only real zeros, has genus $0$ or $1$, and is real on
the real axis, then the points of ${\rm Fix}({\rm Aut}(f))$  are real and are separated
by the zeros of $f$, and the zeros of $f$ are separated by the points of ${\rm Fix}({\rm Aut}(f))$.
\end{theorem}
\noindent
{\bf Proof.} \\
This follow by Laguerre's Theorem on Separation Zeros, \cite{rc} (p. 89) and by Theorem \ref{thm4}. $\qed $ \\

\section{Common zeros of the reciprocals of almost all the automorphic functions}\label{sec10}

\begin{definition}\label{def2}
Let $f(z)$ be a non-constant entire function, let $\mathcal{P}$ be a property that an element in the
automorphic group $\phi\in{\rm Aut}(f)$ can have or does not have. We say that the property $\mathcal{P}$
is common to almost all the automorphic functions of $f(z)$ if except for a finite number of them, all
the automorphic functions $\phi\in{\rm Aut}(f)$ have the property $\mathcal{P}$.
\end{definition}
\noindent
We give in the current section a non-trivial such a property. The property will be: having a common zero
for $1/\phi$, where $\phi\in{\rm Aut}(f)$.

\begin{theorem}\label{thm6}
Let $g(z)$ be an entire function. Let $p(z)$ be a polynomial, $d:=\deg p>0$ and $Z(p)=\{\alpha_1,\ldots,\alpha_d\}
\subseteq\mathbb{C}$. Let $f(z)=p(z)e^{g(z)}$. Then $\forall\,j=1,\ldots,d$, $\alpha_j$ is a common zero of almost all the
reciprocals of the automorphic functions of $f(z)$.
\end{theorem}
\noindent
{\bf Proof.} \\
Let us denote ${\rm Aut}(f)=\{\phi_n(z)\}_n$ and we consider a Weierstrass representation:
$$
f(w)-f(z)=e^{g(w,z)}\prod_{n}\left(1-\frac{w}{\phi_{n}(z)}\right)e^{Q_{\lambda_n}(w/\phi_n(z))}.
$$
Thus we have:
$$
p(w)e^{g(w)}-p(z)e^{g(z)}=e^{g(w,z)}\prod_{n}\left(1-\frac{w}{\phi_{n}(z)}\right)e^{Q_{\lambda_n}(w/\phi_n(z))}.
$$
For the sake of simplicity, let us assume that $f(0)\ne 0$. Consider any integer $j$, $1\le j\le d$, and let
us take the parameter value $z=\alpha_j$. This is a legitimate value of the parameter $z$ for which
the above Weierstrass representation holds true. The reason is that with this parameter we have
$f(0)-f(z)=f(0)-f(\alpha_j)=f(0)-0=f(0)\ne 0$. This follows by: $f(\alpha_j)=p(\alpha_j)e^{g(\alpha_j)}=
0\cdot e^{g(\alpha_j)}=0$. Hence at least one $n$ exists for which $\phi_n(\alpha_j)=\alpha_k$ for
some $1\le k\le d$. We note that $Z(f(w))=Z(p(w))=\{\alpha_1,\ldots,\alpha_d\}$, because the only
solutions of $f(w)=0$ , i.e. $p(w)e^{g(w)}=0$ are (exactly) the solutions of $p(w)=0$ and vice verse.
So in the Weierstrass product of $f(w)-f(\alpha_j)=f(w)$ there are exactly $d$ factors. All the other
factors are "phantom" factors, i.e. 
$$
1-\frac{w}{\phi_n(\alpha_j)}\equiv 1.
$$
Thus except for $d$ factors we have:
$$
\left|\frac{w}{\phi_n(\alpha_j)}\right|=0.
$$
This means that $|\phi_n(\alpha_j)|=\infty$, for all $n$, except for $d$ of them. $\qed $ \\

\section{Sums of the derivatives of the automorphic functions}\label{sec11}

\begin{definition}\label{def3}
Let $f(z)$ be an holomorphic function in some domain $\mathcal{D}\subseteq\mathbb{C}$. A differential
monomial of $f$ is a function of the form:
$$
m_{n_1,\ldots,n_k;m_1,\ldots,m_k}(z)=a\cdot\left(f^{(n_1)}(z)\right)^{m_1}\ldots\left(f^{(n_k)}(z)\right)^{m_k}.
$$
Here $k,n_1,\ldots,n_k,m_1\ldots,m_k\in\mathbb{Z}^+$ and $a\in\mathbb{C}^{\times}$. The weight of the
monomial $m_{n_1,\ldots,n_k;m_1,\ldots,m_k}(z)$ is $w(m_{n_1,\ldots,n_k;m_1,\ldots,m_k}(z))=n_1\cdot m_1+\ldots+n_k\cdot m_k$.
\end{definition}
\noindent
In this section we will discuss the following result:
\begin{theorem}\label{thm7}
Let $f(z)$ be a non-constant entire function, ${\rm Aut}(f)=\{\phi_{0n}(z)\}_n$, $k\in\mathbb{Z}^+$ and
$R>0$. Then there is an identity (independent of $f$) of the form:
$$
2\pi\cdot i\cdot\sum_{|\phi_{0n}(z)|<R}\phi_{0n}^{(k)}(z)=\sum_{j=1}^k m_j(z)\cdot\oint_{|w|=R}\frac{dw}{(f(w)-f(z))^j},
$$
where $w(m_j(z))=k$ and in particular: $m_1(z)=f^{(k)}(z)$ and $m_k(z)=(f'(z))^k$.
\end{theorem}
\noindent
We start by writing explicit formulas for the results already obtained in section \ref{sec3}, Proposition 
\ref{prop2} and Remark \ref{rem12}.

\begin{proposition}\label{prop6}
Let $f(w)$ be a non-constant entire function, ${\rm Aut}(f(z))=\{\phi_{0n}(z)\}_n$ and let us assume that
$f(z)-f(0)\ne 0$. Then
$$
\frac{f'(w)}{f(w)-f(z)}=\frac{\partial g}{\partial w}(w,z)+\sum_n\left(\frac{w}{\phi_{0n}(z)}\right)^{\lambda_n}
\left(\frac{1}{w-\phi_{0n}(z)}\right).
$$
Here $g(w,z)$ and the $\lambda_n$ are the data of the Weierstrass presentation of $f(w)-f(z)$ (Proposition \ref{prop3}). 
\end{proposition}
\noindent
{\bf Proof.} \\
If $f(z)-f(0)\ne 0$, then there is a function $g(w,z)$, entire in $w$ and there are non-negative
integers $\lambda_n(w,z)$ which we will sometimes denote by $\lambda_n$, such that:
$$
f(w)-f(z)=\exp{\left(g(w,z)\right)}\prod_{n=1}^{\infty}E\left(\frac{w}{\phi_{0n}(z)},\lambda_n\right)=
$$
$$
=\exp{\left(g(w,z)\right)}\prod_{n=1}^{\infty}\left(1-\frac{w}{\phi_{0n}(z)}\right)e^{Q_{\lambda_n}(w/\phi_{0n}(z))},
$$
where if $\lambda_n>0$, then:
$$
Q_{\lambda_n}\left(\frac{w}{\phi_{0n}(z)}\right)=\left(\frac{w}{\phi_{0n}(z)}\right)+\frac{1}{2}
\left(\frac{w}{\phi_{0n}(z)}\right)^2+\ldots +\frac{1}{\lambda_n}\left(\frac{w}{\phi_{0n}(z)}\right)^{\lambda_n},
$$
and $Q_0(w/\phi_{0n}(z))\equiv 0$. Taking the logarithm of both sides of the identity we obtain:
\begin{equation}\label{eq8}
\log\left(f(w)-f(z)\right)=
\end{equation}
$$
=g(w,z)+\sum_n\left(\log\left(1-\frac{w}{\phi_{0n}(z)}\right)+
Q_{\lambda_n}\left(\frac{w}{\phi_{0n}(z)}\right)\right).
$$
We $\partial w$ both sides of the last identity, equation (\ref{eq8}) and obtain our result:
$$
\frac{f'(w)}{f(w)-f(z)}=\frac{\partial g}{\partial w}(w,z)+\sum_n\left(\frac{w}{\phi_{0n}(z)}\right)^{\lambda_n}
\left(\frac{1}{w-\phi_{0n}(z)}\right).\,\,\,\,\,\qed
$$
\begin{proposition}\label{prop7}
Let $f(w)$ be a non-constant entire function, ${\rm Aut}(f(z))=\{\phi_{0n}(z)\}_n$ and let us assume that
$f(z)-f(0)\ne 0$. Then
$$
\frac{f'(z)}{f(w)-f(z)}=-\frac{\partial g}{\partial z}(w,z)+\sum_n\left(\frac{w\phi_{0n}^{'}(z)}{\phi_{0n}(z)}\right)
\left(\frac{w}{\phi_{0n}(z)}\right)^{\lambda_n}\left(\frac{1}{w-\phi_{0n}(z)}\right).
$$
Here $g(w,z)$ and the $\lambda_n$ are the data of the Weierstrass presentation of $f(w)-f(z)$ (Proposition \ref{prop3}). 
\end{proposition}
\noindent
{\bf Proof.} \\
We $\partial z$ both sides of equation (\ref{eq8}) and obtain our result:
$$
\frac{f'(z)}{f(w)-f(z)}=-\frac{\partial g}{\partial z}(w,z)+\sum_n\left(\frac{w\phi_{0n}^{'}(z)}{\phi_{0n}(z)}\right)
\left(\frac{w}{\phi_{0n}(z)}\right)^{\lambda_n}\left(\frac{1}{w-\phi_{0n}(z)}\right).\,\,\,\,\,\qed
$$
\noindent
Let us write the identity of Proposition \ref{prop7} in the following form:
$$
\frac{f'(z)}{f(w)-f(z)}+\frac{\partial g}{\partial z}(w,z)=\sum_n\frac{\phi_{0n}^{'}(z)}{(\phi_{0n}(z))^{\lambda_n+1}}\cdot
\left(\frac{w^{\lambda_n+1}}{w-\phi_{0n}(z)}\right).
$$
Let $R>0$, $z\in\mathbb{C}$ be fixed so that the circle $|w|=R$ does not contain $\{\phi_{0n}(z)\}_n$. In that
event we have $f(z)-f(w)\ne 0$, $\forall\,|w|=R$ (because $f(z)-f(w)=0\Leftrightarrow w=\phi_{0n}(z)$ for some
$n\in\mathbb{Z}$). Thus we can path integrate our identity on $|w|=R$ and obtain:
$$ f'(z)\oint_{|w|=R}\frac{dw}{f(w)-f(z)}+\oint_{|w|=R}\frac{\partial g}{\partial z}(w,z)dw=
\sum_n\frac{\phi_{0n}^{'}(z)}{(\phi_{0n}(z))^{\lambda_n+1}}\cdot\oint_{|w|=R}\frac{w^{\lambda_n+1}dw}{w-\phi_{0n}(z)}.
$$
Since the function $\partial g(w,z)/\partial z$ is entire in $w\in\mathbb{C}$, it follows by the Theorem
of Cauchy that:
$$
\oint_{|w|=R}\frac{\partial g}{\partial z}(w,z)dw=0.
$$
By the generalized argument principle we have:
$$
\sum_n\frac{\phi_{0n}^{'}(z)}{(\phi_{0n}(z))^{\lambda_n+1}}\cdot\oint_{|w|=R}\frac{w^{\lambda_n+1}dw}{w-\phi_{0n}(z)}=
2\pi\cdot i\sum_{|\phi_{0n}(z)|<R}\phi_{0n}^{'}(z).
$$
Putting what we have so far together:
\begin{equation}\label{eq9}
f'(z)\oint_{|w|=R}\frac{dw}{f(w)-f(z)}=2\pi\cdot i\sum_{|\phi_{0n}(z)|<R}\phi_{0n}^{'}(z).
\end{equation}
\noindent
We just proved Theorem \ref{thm7} for the case $k=1$.

\begin{remark}\label{rem22}
If we add the assumption on $z$, that $f(z)-f(0)\ne 0$ then we have the Weierstrass presentation (Proposition \ref{prop3}):
$$
f(w)-f(z)=\exp{\left(g(w,z)\right)}\prod_n\left(\frac{\phi_{0n}(z)-w}{\phi_{0n}(z)}\right)e^{Q_{\lambda_n}(w/\phi_{0n}(z))}.
$$
In other words $Z_w(f(w)-f(z))=\{\phi_{0n}(z)\}_n$. So $1/(f(w)-f(z))$ is meromorphic in $w$, it has no zeros, and its
total set of poles are $\{\phi_{0n}(z)\}_n$. What is the residue:
$$
{\rm Res}\left(\frac{1}{f(w)-f(z)},\phi_{0n}(z)\right)?
$$
Let us assume for simplicity that all the zeros of $f(w)-f(z)$ are simple. Then this residue is:
$$
\lim_{w\rightarrow\phi_{0n}(z)}(w-\phi_{0n}(z))\cdot\frac{1}{f(w)-f(z)}=\frac{1}{f'(\phi_{0n}(z))}.
$$
We conclude that:
$$
\oint_{|w|=R}\frac{dw}{f(w)-f(z)}=2\pi\cdot i\sum_{|\phi_{0n}(z)|<R}\frac{1}{f'(\phi_{0n}(z))}.
$$
By equation (\ref{eq9}):
$$
f'(z)\left(2\pi\cdot i\sum_{|\phi_{0n}(z)|<R}\frac{1}{f'(\phi_{0n}(z))}\right)=2\pi\cdot i\sum_{|\phi_{0n}(z)|<R}\phi_{0n}^{'}(z).
$$
By the automorphic equation $f(\phi_{0n}(z))=f(z)$ and the chain rule, we get: $\phi_{0n}^{'}(z)f'(\phi_{0n}(z))=f'(z)$. Hence:
$$
\frac{1}{f'(\phi_{0n}(z))}=\frac{\phi_{0n}^{'}(z)}{f'(z)}.
$$
This agrees with our equation (\ref{eq9}). However, we proved equation (\ref{eq9}) without the extra assumption on the
simplicity of the zeros of $f(w)-f(z)$.
\end{remark}
\noindent
Let us show how to compute the next case, $k=2$, of Theorem \ref{thm7} and in fact any case follows
just as simple using inductive argument. We do the obvious and apply the operator $\partial z$:
$$
\left(\frac{f'(z)}{f(w)-f(z)}\right)'_z=\frac{f''(z)}{f(w)-f(z)}+\frac{(f'(z))^2}{(f(w)-f(z))^2}.
$$
Hence using equation (\ref{eq9}) we obtain:
\begin{equation}\label{eq10}
f''(z)\oint_{|w|=R}\frac{dw}{f(w)-f(z)}+(f'(z))^2\oint_{|w|=R}\frac{dw}{(f(w)-f(z))^2}=
\end{equation}
$$
=2\pi\cdot i\sum_{|\phi_{0n}(z)|<R}\phi_{0n}^{''}(z).
$$
\noindent
This proves the second case, $k=2$ of Theorem \ref{thm7}. Let us do one more explicit computation
and prove (explicitly) the third case, $k=3$ too. Once more we do the obvious and apply the
operator $\partial z$ to the case $k=2$:
$$
\left(\frac{f''(z)}{f(w)-f(z)}\right)'_z+\left(\frac{(f'(z))^2}{(f(w)-f(z))^2}\right)'_z=
$$
$$
=\frac{f^{(3)}(z)}{f(w)-f(z)}+\frac{3f'(z)f''(z)}{(f(w)-f(z))^2}+\frac{(f'(z))^3}{(f(w)-f(z))^3}.
$$
Using equation (\ref{eq10}) we finally get:
$$
f^{(3)}(z)\oint_{|w|=R}\frac{dw}{f(w)-f(z)}+3f'(z)f''(z)\oint_{|w|=R}\frac{dw}{(f(w)-f(z))^2}+
$$
\begin{equation}\label{eq11}
\end{equation}
$$
+(f'(z))^3\oint_{|w|=R}\frac{dw}{(f(w)-f(z))^3}=2\pi\cdot i\sum_{|\phi_{0n}(z)|<R}\phi_{0n}'''(z).
$$
This proves the third case, $k=3$ of Theorem \ref{thm7}. It is clear how to proceed inductively
by repeatedly applying the operator $\partial z$ and forming a simple wight calculation on the
resulting differential monomials. $\qed $

\section{An application of Jensen's Theorem to the automorphic group of an entire function}\label{sec12}

Here is one of the most important theorems in analysis. \\
\\
{\bf Theorem. (Jensen's Theorem)} {\it Let $f(z)$ be analytic for $|z|<R$. Suppose that $f(0)$ is not
zero, and let $r_1,r_2,\ldots,r_n,\ldots $ be the moduli of the of the zeros of $f(z)$ in the disk
$|z|<R$, arranged in a non-decreasing sequence. Then if $r_n\le r\le r_{n+1}$,
$$
\log\frac{r^n |f(0)|}{r_1r_2\ldots r_n}=\frac{1}{2\pi}\int_0^{2\pi}\log|f(re^{i\theta})|d\theta,
$$
where every zero is counted the number of its multiplicity.
} \\
\\
Let $f(w)$ be a non-constant entire function, let $z\in\mathbb{C}-f^{-1}(f(0))$ and let us apply the
Theorem of Jensen to the entire function $f(w)-f(z)$ of the variable $w$. Then indeed $f(0)-f(z)$
is not zero, and the parameter $R$ in Jensen's Theorem can be an arbitrary positive number. The zero
set of $f(w)-f(z)$ is the ${\rm Aut}(f(w))$-orbit of $z$. Thus: $Z(f(w)-f(z))=\{\phi_{0n}(z)\}_n$
and we may assume that the moduli of the zeros are arranged in a non-decreasing order. Thus
$|\phi_{00}(z)|\le |\phi_{01}(z)|\le |\phi_{02}(z)|\le\ldots $. In other words, in terms of
the notation in Jensen's Theorem we have $r_{j}=|\phi_{0,j-1}(z)|$. We conclude that if $|\phi_{0n}(z)|
\le r<|\phi_{0,n+1}(z)|$, then we have the following identity:
$$
\log\frac{r^{n+1}|f(0)-f(z)|}{|\phi_{00}(z)\phi_{01}(z)\ldots\phi_{0n}(z)|}=\frac{1}{2\pi}\int_0^{2\pi}
\log\left|f(re^{i\theta})-f(z)\right|d\theta.
$$
Equivalently:
$$
\left|\prod_{j=0}^n\phi_{0j}(z)\right|=r^{n+1}|f(0)-f(z)|\exp\left\{-\frac{1}{2\pi}\int_0^{2\pi}
\log\left|f(re^{i\theta})-f(z)\right|d\theta\right\}.
$$
If we take (as is possible) $r=|\phi_{0n}(z)|$ then:
$$
\left|\prod_{j=0}^{n-1}\phi_{0j}(z)\right|=|\phi_{0n}(z)|^n|f(0)-f(z)|\exp\left\{-\frac{1}{2\pi}\int_0^{2\pi}
\log\left|f(|\phi_{0n}(z)|e^{i\theta})-f(z)\right|d\theta\right\}.
$$
This gives a recursion between $|\phi_{0n}(z)|$ on the one hand, and the product of the previous terms
$|\prod_{j=0}^{n-1}\phi_{0j}(z)|$ on the other hand. Thus this determine $|\phi_{0n}(z)|$ uniquely? Let us
consider the following function of $r$:
$$
\psi_n(r)=r^n\left|h(0)\right|\exp\left\{-\frac{1}{2\pi}\int_0^{2\pi}\log\left|h(re^{i\theta})\right|d\theta\right\}.
$$
Here $h(z)$ is a non-constant entire function such that $h(0)\ne 0$. If this function $\psi_n(r)$ turns out to
be a strictly monotone function of $r$ (in our case, necessarily increasing), then our recursion uniquely
determines $|\phi_{0n}(z)|$ (up to multiplicity in $r$). We might go about as follows:
$$
r^n\exp\left\{-\frac{1}{2\pi}\int_0^{2\pi}\log\left|h(re^{i\theta})\right|d\theta\right\}=
e^{\log r^n}\exp\left\{-\frac{1}{2\pi}\int_0^{2\pi}\log\left|h(re^{i\theta})\right|d\theta\right\}=
$$
$$
=\exp\left\{\log r^n-\frac{1}{2\pi}\int_0^{2\pi}\log\left|h(re^{i\theta})\right|d\theta\right\}=
\exp\left\{\frac{1}{2\pi}\int_0^{2\pi}\log\left(\frac{r^n}{|h(re^{i\theta})|}\right)d\theta\right\}=
$$
$$
=\exp\left\{\frac{1}{2\pi}\int_0^{2\pi}\log\left|\frac{(re^{i\theta})^n}{h(re^{i\theta})}\right|d\theta\right\}
$$
The point is that $z^n/h(z)$ is not holomorphic for $|z|<r+\epsilon$ because of the zeros of $h(z)$
with moduli smaller than $r+\epsilon$, but, as in the proof of the Theorem of Jensen we divide out
those zeros by dividing $h(z)$ by the corresponding finite Blaschke product $B_n(z)$, without changing the modulus
of $h(z)$ on $|z|=r+\epsilon$. Thus we have for all the $z\in\mathbb{C}$, $|z|=r+\epsilon<|\phi_{0,n+1}(z)|$:
$$
\left|\frac{h(z)}{B_n(z)}\right|=\left|h(z)\right|.
$$
So on that circle of integration we have:
$$
\left|\frac{z^n}{h(z)}\right|=\left|\frac{z^nB_n(z)}{h(z)}\right|,
$$
and this function is analytic. The monotonicity now follows. A better approach:
We have (write the recursion a bit different),
$$
\left|\prod_{j=0}^{n-1}\phi_{0j}(z)\right|^{1/n}=\left|\phi_{0n}(z)\right|\left|f(0)-f(z)\right|^{1/n}\times
\exp\left\{-\frac{1}{2\pi n}\int_0^{2\pi}\log\left|f(|\phi_{0n}(z)|e^{i\theta})-f(z)\right|d\theta\right\}.
$$
So the right hand side equals the geometric mean of the sequence $|\phi_{00}(z)|,\ldots,|\phi_{0,n-1}(z)|$
and this is a part of a non-decreasing sequence so those means are also non-decreasing.

\section{A computation of the integral
$$
\frac{1}{2\pi}\int_0^{2\pi}\log\left|f(re^{i\theta})-f(z)\right|d\theta.
$$
}\label{sec13}

Assuming $z\in\mathbb{C}-f^{-1}(f(0))$, i.e. $f(0)-f(z)\ne 0$, we have:

\begin{equation}
\label{eq12}
f(w)-f(z)=e^{g(w,z)}\prod_{n=0}^{\infty}\left(1-\frac{w}{\phi_{0n}(z)}\right)e^{Q_{\lambda_n}(w/\phi_{0n}(z))}.
\end{equation}
So 
$$
\left|f(w)-f(z)\right|=e^{\Re\,g(w,z)}\prod_{n=0}^{\infty}\left|1-\frac{w}{\phi_{0n}(z)}\right|e^{\Re\,Q_{\lambda_n}(w/\phi_{0n}(z))}.
$$
Hence
$$
\log\left|f(w)-f(z)\right|=\Re\,g(w,z)+\sum_{n=0}^{\infty}\left\{\log\left|1-\frac{w}{\phi_{0n}(z)}\right|+
\Re\,Q_{\lambda_n}\left(\frac{w}{\phi_{0n}(z)}\right)\right\}.
$$
The function $\Re\,g(w,z)$ is harmonic in all of the $w$-plane. Hence:
$$
\frac{1}{2\pi}\int_0^{2\pi}\Re\,g(re^{i\theta},z)d\theta=\Re\,g(0,z).
$$
If we substitute $w=0$ in equation (\ref{eq12}) we get $f(0)-f(z)=e^{g(0,z)}$, so $|f(0)-f(z)|=e^{\Re\,g(0,z)}$ and
$\log|f(0)-f(z)|=\Re\,g(0,z)$. Hence, we proved the following:
$$
\frac{1}{2\pi}\int_0^{2\pi}\Re\,g(re^{i\theta},z)d\theta=\log|f(0)-f(z)|.
$$
Similarly (and in fact much simpler):
$$
\frac{1}{2\pi}\int_0^{2\pi}\Re\,Q_{\lambda_n}\left(\frac{re^{i\theta}}{\phi_{0n}(z)}\right)d\theta=0.
$$
We are left with the computation of:
$$
\frac{1}{2\pi}\int_0^{2\pi}\log\left|1-\frac{re^{i\theta}}{\phi_{0n}(z)}\right|d\theta.
$$
If $r<|\phi_{0n}(z)|$, then $1-(re^{i\theta})/\phi_{0n}(z)$ never vanishes on $0\le\theta<2\pi$,
and so the function $\log|1-(re^{i\theta})/\phi_{0n}(z)|$ is harmonic for $|re^{i\theta}|<|\phi_{0n}(z)|$
and again the mean value property implies that:
$$
\frac{1}{2\pi}\int_0^{2\pi}\log\left|1-\frac{re^{i\theta}}{\phi_{0n}(z)}\right|d\theta=\log|1-0|=0.
$$
Thus we are left with computing the integral
$$
\frac{1}{2\pi}\int_0^{2\pi}\log\left|1-\frac{re^{i\theta}}{\phi_{0n}(z)}\right|d\theta,
$$
for the case $|\phi_{0n}(z)|\le r$. In this case the integrand $\log|1-w/\phi_{0n}(z)|$ is singular exactly
in $w=\phi_{0n}(z)$ which lies within $|w|\le r$. Here is formula 15 on page 531 of the book \cite{gr}:

\begin{equation}\label{eq13}
\int_0^{n\pi}\log\left(1-2a\cos\theta+a^2\right)d\theta=\left\{\begin{array}{lll} 0 & , & a^2\le 1 \\
n\pi\log a^2 & , & a^2\ge 1 \end{array}\right..
\end{equation}
We need to evaluate:
$$
\frac{1}{2\pi}\int_0^{2\pi}\log\left|1-Re^{i\theta}\right|d\theta=^{(1\le R)}\frac{1}{2\pi}\int_0^{2\pi}
\log\left|\frac{1}{R}-e^{i\theta}\right|d\theta+\log R=
$$
$$
=\frac{1}{2\pi}\int_0^{2\pi}\log\left|a-e^{i\theta}\right|d\theta+\log R=^{(0<a=1/R\le 1)}
$$
$$
=\frac{1}{2\pi}\cdot\frac{1}{2}\int_0^{2\pi}\log\left(1+a^2-2a\cos\theta\right)d\theta+\log R=\log R,
$$
where in the last step we used the formula (\ref{eq13}). The conclusion is:
$$
\frac{1}{2\pi}\int_0^{2\pi}\log\left|1-\frac{re^{i\theta}}{\phi_{0n}(z)}\right|d\theta=
\left\{\begin{array}{lll} 0 & , & r<|\phi_{0n}(z)| \\ \log(r/|\phi_{0n}(z)|) & , & r\ge |\phi_{0n}(z)| \end{array}\right..
$$
Finally we get:
$$
\frac{1}{2\pi}\int_0^{2\pi}\log\left|f(re^{i\theta})-f(z)\right|d\theta=
$$
$$
=\frac{1}{2\pi}\int_0^{2\pi}\left\{
\Re\,g(re^{i\theta},z)+\sum_{n=0}^{\infty}\left(\log\left|1-\frac{re^{i\theta}}{\phi_{0n}(z)}\right|
+\Re\,Q_{\lambda_n}\left(\frac{re^{i\theta}}{\phi_{0n}(z)}\right)\right)\right\}d\theta=
$$
$$
=\log\left|f(0)-f(z)\right|+\sum_{|\phi_{0n}(z)|\le r}\log\frac{r}{|\phi_{0n}(z)|}+0=
$$
$$
=\log\left(\left|f(0)-f(z)\right|\times\prod_{|\phi_{0n}(z)|\le r}\frac{r}{|\phi_{0n}(z)|}\right).
$$
We note that in fact our computation proved the Theorem of Jensen.

\section{The product of the automorphic functions}\label{sec14}

Let $f(w)$ be a non-constant entire function, and let $z\in\mathbb{C}-f^{-1}(f(0))$. Let ${\rm Aut}(f)=\{\phi_{0n}(z')\}_n$
and let us suppose that we arranged the automorphic functions in a non-decreasing order of their orbit at $z$. Thus:
$|\phi_{00}(z)|\le|\phi_{01}(z)|\le\ldots $. Then we have the following identity:
$$
\left|\prod_{j=0}^{n-1}\phi_{0j}(z)\right|=\left|\phi_{0n}(z)\right|^n\left|f(0)-f(z)\right|\exp\left\{-\frac{1}{2\pi}
\int_0^{2\pi}\log\left|f(\left|\phi_{0n}(z)\right|e^{i\theta})-f(z)\right|d\theta\right\}.
$$
One is tempted to compare this formula with the Vieta formula that corresponds to the special case of a
non-constant polynomial of degree $d\ge 1$, $f(w)=p_d(w)$. We recall that if $p_d(w)=a_dw^d+\ldots+a_1w+a_0$,
where $a_d\in\mathbb{C}-\{0\}$, then: $p_d(w)-p_d(z)=a_d(w-\phi_{00}(z))\ldots(w-\phi_{0,d-1}(z))$. Plugging
into the last formula the value $w=0$, we obtain:

\begin{equation}\label{eq14}
\phi_{00}(z)\ldots\phi_{0,d-1}(z)=(-1)^da_d^{-1}\left(p_d(0)-p_d(z)\right).
\end{equation}
Taking absolute values we get:
$$
\left|\phi_{00}(z)\ldots\phi_{0,d-1}(z)\right|=\left|a_d\right|^{-1}\left|p_d(0)-p_d(z)\right|.
$$
Let $N$ be an integer such that $0\le N\le d-1$. Then:
$$
\left|\prod_{j=0}^{N-1}\phi_{0j}(z)\right|=\left|\prod_{k=N}^{d-1}\phi_{0k}(z)\right|^{-1}\left|a_d\right|^{-1}
\left|p_d(0)-p_d(z)\right|.
$$
We note that $\forall\,j$, $\phi_{0j}(z)\ne 0$ so that the expression on the right hand side is defined.
Comparing that to the more general expression:
$$
\left|\prod_{j=0}^{N-1}\phi_{0j}(z)\right|=\left|\phi_{0N}(z)\right|^N\left|p_d(0)-p_d(z)\right|\exp\left\{-\frac{1}{2\pi}
\int_0^{2\pi}\log\left|p_d(\left|\phi_{0N}(z)\right|e^{i\theta})-p_d(z)\right|d\theta\right\},
$$
we deduce the following identity:
$$
\exp\left\{-\frac{1}{2\pi}\int_0^{2\pi}\log\left|p_d(\left|\phi_{0N}(z)\right|e^{i\theta})-p_d(z)\right|d\theta\right\}=
\left|\prod_{k=N}^{d-1}\phi_{0k}(z)\right|\left|\phi_{0N}(z)\right|^N\left|a_d\right|.
$$
If $z$ has the smallest modulus in its ${\rm Aut}(p_d)$-orbit and $N=0$, we obtain the following interesting
identity:
$$
\exp\left\{\frac{1}{2\pi}\int_0^{2\pi}\log\left|p_d(ze^{i\theta})-p_d(z)\right|d\theta\right\}=
\left|\prod_{k=0}^{d-1}\phi_{0k}(z)\right|\left|a_d\right|.
$$
Using the identity in equation (\ref{eq14}) this proves that:
$$
\exp\left\{\frac{1}{2\pi}\int_0^{2\pi}\log\left|p_d(ze^{i\theta})-p_d(z)\right|d\theta\right\}=
\left|p_d(0)-p_d(z)\right|.
$$
This last identity could be seen to be true because of the mean value property of the harmonic function:
$\log|p_d(w)-p_d(z)|$ in the disk $|w|<|z|$. This reproduces the Vieta identity and shows that we generalized
it even in the special case of non-constant polynomials.

\section{Consequences to ${\rm Aut}(f)$ that follow from the classical theory of entire
functions}\label{sec15}

Let $f$ be a non-constant entire function of a single complex variable. The basic observation that the set
of all the zeros of the entire function of $w$, $f(w)-f(z)$, where $z\in\mathbb{C}$ is a parameter, is the
${\rm Aut}(f)$-orbit of $z$, plus the fact that for a fixed $z$, the functions $f(w)$ and $f(w)-f(z)$ differ
by a constant suggest that the elements of ${\rm Aut}(f(z))=\{\phi_{0n}(z)\}_n$ share properties with the
zeros of $f(w)+c$ for any (generic) constant. In this section we will repeatedly refer to the classical book
\cite{l}.

\begin{remark}\label{rem23}
For a fixed value of the parameter $z\in\mathbb{C}$ the $w$-entire functions $f(w)$ and $f(w)-f(z)$ have
the same order $\rho $ and the same type $\sigma $.
\end{remark}
\noindent
An immediate consequence that follows by the result on (page 16 of \cite{l}) is: \\

\begin{theorem}\label{thm8}
Let $f(w)$ be a non-constant entire function, and let ${\rm Aut}(f)=\{\phi_{0n}(z)\}_n$. Then the convergence
exponent of the sequence $\{\phi_{0n}(z)\}_n$ for any $z\in\mathbb{C}$ does not exceed the order of $f(w)$.
\end{theorem}

\begin{remark}\label{rem24}
This theorem is interesting only in the case of entire functions of a finite order $\rho<\infty$.
\end{remark}
\noindent
Using Theorem 7 on page 16 of \cite{l} and the representation formula in the first equation in the proof of
the Theorem of Wiman (page 72, \cite{l}) we deduce some interesting consequences on the automorphic group of a non-constant
entire function of order less than one.

\begin{theorem}\label{thm9}
Let $f(w)$ be a non-constant entire function of order less than one, and let ${\rm Aut}(f)=\{\phi_{0n}(w)\}_n$.
Then $\forall\,z\in\mathbb{C}$, the convergence exponent of the sequence $\{\phi_{0n}(z)\}_n$ equals the order
of $f$. In addition, if the order of $f$ is not zero, then it is of maximal, minimal or of a normal type according
to whether the upper density of any ${\rm Aut}(f)$-orbit of any $z\in\mathbb{C}$:
$$
\triangle\left(\{\phi_{0n}(z)\}_n\right)=\overline{\lim}_{r\rightarrow\infty}\frac{n_{\{\phi_{0n}(z)\}_n}(r)}{r^{\rho}},
$$
equals infinity, zero or equals a number different from zero or infinity. Here we use the notation:
$$
n_{\{\phi_{0n}(z)\}_n}(r)=\left|\{n\,|\,|\phi_{0n}(z)|<r\}\right|,
$$
i.e. the counting function of the elements in the ${\rm Aut}(f)$-orbit of $z$ of modulus less than $r$. In particular
the upper density of any ${\rm Aut}(f)$-orbit of any $z\in\mathbb{C}$ is independent of $z$.
\end{theorem}
\noindent
{\bf Proof.} \\
$\forall\,z\in\mathbb{C}$, the order $\rho $ and the type $\sigma $ of the entire function of $w$, $f(w)-f(z)$ equal
those of $f(w)$. Hence (like in the first equation in the proof above of the Theorem of Wiman) we have:
$$
f(w)-f(z)=C\cdot w^m\prod_{k=1}^{\infty}\left(1-\frac{w}{\phi_{0n}(z)}\right).
$$
Here $C=C(z)$ is a function of $z$ only, and $C(z)\ne 0$ $\forall\,z\in\mathbb{C}$. So $f(w)-f(z)$ is a canonical
product and Theorem 7 of \cite{l} applies. $\qed $ \\
\\
We will quote few results from \cite{cy}. We will use those to extract information on the function $C(z)$ that
appears in the proof of Theorem \ref{thm9}. On page 207 of \cite{cy}: \\
\\
{\bf Theorem A.2. (Clunie)} {\it
Let $f(z)$ and $g(z)$ be entire with $g(0)=0$. Let $\rho $ satisfy $0<\rho<1$ and $c(\rho)=(1-\rho)^2/48$.
Then for $R\ge 0$, $M(R,f\circ g)\ge M(c\rho M(\rho\cdot R,g),f)$.
} \\
\\
On page 208 of \cite{cy}: \\
\\
{\bf Corollary A.1. (P\"{o}lya)} {\it
Let $f(z)$, $g(z)$ and $h(z)$ be entire functions with $h(z)=f(g(z))$. If $g(0)=0$, then there exists an absolute
constant $c$, $0<c<1$, such that for all $r>0$ the following inequality holds:
$$
M(r,h)\ge M\left(c\cdot M\left(\frac{r}{2},g\right),f\right).
$$
When $g(0)\ne 0$, the corresponding inequality should read:
$$
M(r,f\circ g)\ge M\left(c\cdot M\left(\frac{r}{2},g\right)-|g(0)|,f\right).
$$
The constant $c$ can be chosen to be $1/8$.
} \\
\\
On page 209 of \cite{cy}: \\
\\
{\bf Theorem A.3.} {\it
If $f(z)$ and $g(z)$ are two entire functions such that $f\circ g$ is of a finite order (lower order), then
either: (i) $g(z)$ is a polynomial and $f(z)$ is of a finite order (lower order), or (ii) $g(z)$ is not
a polynomial but a function of a finite order (lower order) and $f(z)$ is of zero order (lower order).
} \\
\\
We conclude the following:

\begin{remark}\label{rem25}
If $f(w)$ is an entire function of order less than one, and greater than zero, then $f(w)$ has infinitely many
zeros.
\end{remark}
\noindent
{\bf A Proof of the claim in remark \ref{rem25}.} \\
If $f(w)$ has no zeros then $f(w)=e^{g(w)}$ for some entire function $g$. This corresponds in Theorem A.3. 
in \cite{cy} to case (i) with $e^w\circ g(w)$. Here $g(w)$ is a non-constant polynomial of degree $d\ge 1$.
Hence the order of $f(w)$ equals to $d\ge 1$ and this contradicts the assumption that $f(w)$ has order less
than one. If $f(a_1)=0$ then $f(w)\cdot (w-a_1)^{-1}$ is a non-constant entire function of the same order
as the order of $f(w)$. So if $f(w)$ had finitely many zeros, then we could have found a polynomial $p(w)$
such that $f(w)/p(w)$ had non zeros and was of the same order as that of $f(w)$. This contradicts the first
part of our proof and thus the proof of the claim in Remark \ref{rem25} is now completed. $\qed $ \\
\\
Finally, the function $C(z)$ in the proof of Theorem \ref{thm9} is entire because the infinite product is, and
it never vanishes and since $f(w)-f(z)$ is symmetric in $w$ and $z$, $C(z)$ must be a constant, i.e. 
independent of $z$ (because its order is less than one).

\begin{theorem}\label{thm10}
Let $f(w)$ be a non-constant entire function of a finite order $\rho $ and let ${\rm Aut}(f)=\{\phi_{0n}(w)\}_n$. 
Then $\forall\,z\in\mathbb{C}$, the upper density of the ${\rm Aut}(f)$-orbit of $z$ satisfies: \\
1. If $f(w)$ is of type not greater than $\sigma $ with respect to the order $\rho $ then,
$$
\overline{\triangle}\left(\{\phi_{0n}(z)\}_n\right)=\overline{\lim}_{r\rightarrow\infty}
\frac{n_{\{\phi_{0n}(z)\}_n}(r)}{r^{\rho}}\le\sigma\cdot e\cdot\rho.
$$
The lower density:
$$
\underline{\triangle}\left(\{\phi_{0n}(z)\}_n\right)=\underline{\lim}_{r\rightarrow\infty}
\frac{n_{\{\phi_{0n}(z)\}_n}(r)}{r^{\rho}}
$$
satisfies $\underline{\triangle}\left(\{\phi_{0n}(z)\}_n\right)\le\sigma\rho$. \\
2. We have the following two identities for these densities:
$$
\overline{\triangle}\left(\{\phi_{0n}(z)\}_n\right)=\overline{\lim}_{n\rightarrow\infty}\frac{n}{|\{\phi_{0n}(z)\}_n|^{\rho}},\,\,
\underline{\triangle}\left(\{\phi_{0n}(z)\}_n\right)=\underline{\lim}_{n\rightarrow\infty}\frac{n}{|\{\phi_{0n}(z)\}_n|^{\rho}}.
$$
\end{theorem}
\noindent
{\bf Proof.} \\
Part 1 follows by Theorem 3 on page 19 of \cite{llst}. Part 2 follows by Problem 2 on page 17 of \cite{llst}. $\qed $ \\
\\
We can refine Part 1 of Theorem \ref{thm10}:

\begin{theorem}\label{thm11}
Let $\rho>0$, $A>0$ and let $f(w)$ be a non-constant entire function for which there exists a constant $M>0$
such that $|f(w)|\le M\cdot e^{A|z|^{\rho}}$, $\forall\,w\in\mathbb{C}$. Let ${\rm Aut}(f)=\{\phi_{0n}(w)\}_n$.
We assume that $z\in\mathbb{C}-f^{-1}(f(0))$ and that the automorphic functions are arranged so that,
$|\phi_{00}(z)|\le|\phi_{01}(z)|\le|\phi_{02}(z)|\le\ldots $. Then:
$$
(\rho\cdot e\cdot A)^{-1/\rho}\cdot M\left(\frac{1+e^{A|z|^{\rho}}}{f(0)-f(z)}\right)^{-1/n}n^{1/\rho}\le\left(\prod_{k=1}^n
|\phi_{0k}(z)|\right)^{1/n}\le|\phi_{0n}(z)|,\,\,\forall\,n\in\mathbb{Z}^+.
$$
If
$$
M\left(\frac{1+e^{A|z|^{\rho}}}{f(0)-f(z)}\right)<1,
$$
then the left hand side in the double inequality is $(\rho\cdot e\cdot A)^{-1/\rho}\cdot n^{1/\rho}$.
\end{theorem}
\noindent
{\bf Proof.} \\
We define an auxiliary non-constant entire function of $w$:
$$
F(w)=\frac{f(w)-f(z)}{f(0)-f(z)}.
$$
Then $F(0)=1$, and
$$
|F(w)|\le\frac{|f(w)|+|f(z)|}{|f(0)-f(z)|}\le\frac{Me^{A|w|^{\rho}}+Me^{A|z|^{\rho}}}{|f(0)-f(z)|}=
\left(M\left(\frac{1+e^{A(|z|^{\rho}-|w|^{\rho})}}{|f(0)-f(z)|}\right)\right)e^{A|w|^{\rho}}\le
$$
$$
\le\left(M\left(\frac{1+e^{A|z|^{\rho}}}{|f(0)-f(z)|}\right)\right)e^{A|w|^{\rho}}\le M_1\cdot e^{A|w|^{\rho}},
$$
where:
$$
M_1=\max\left\{1,M\left(\frac{1+e^{A|z|^{\rho}}}{|f(0)-f(z)|}\right)\right\}.
$$
The function $F(w)$ satisfies the conditions of Proposition 1, in \cite{s}, with $M$ instead of
$M_1$. The result now follows. $\qed $ \\
\\
We recall the following result: \\
\\
{\bf Theorem 1., \cite{llst}} {\it
The convergence exponent of the zero set of an entire function $f$ of non-integer order is equal to the order of growth of $f$.
} \\
\\
This theorem implies:

\begin{theorem}\label{thm12}
Let $f(w)$ be a non-constant entire function of non-integer order $\rho $. Let the ${\rm Aut}(f)$-orbit of a point
$z\in\mathbb{C}$ be $\{\phi_{0n}(z)\}_n\}$, then the convergence exponent of this orbit equals $\rho $.
\end{theorem}
\noindent
Next we have: \\
\\
{\bf Theorem 2. \cite{llst}} {\it
If the order $\rho $ of an entire function $f(z)$ is not an integer, then its type $\sigma_f$ and the upper density
of its zeros $\overline{\triangle}_f$ simultaneously are equal either to zero, or to infinity, or to positive numbers.
} \\
\\
This immediately implies:

\begin{theorem}\label{thm13}
If the order $\rho $ of an entire function $f(w)$ is not an integer, then its type $\sigma_f$ and the upper
density of any ${\rm Aut}(f)$-orbit, $\{\phi_{0n}(z)\}_n$ ($z\in\mathbb{C}$) are equal either to zero, or to
infinity, or to positive numbers.
\end{theorem}

\section{The relations between scattering theory and automorphic functions}\label{sec20}

The book \cite{lp} deals with a discrete subgroup $\Gamma $ of the group of the fractional linear transformations,
$$
w\rightarrow\frac{aw+b}{cw+d},\,\,\,\,\,ad-bc=1,\,\,\,\,\,a,b,c,d\in\mathbb{R}.
$$
The Riemannian metric $(dx^2+dy^2)\cdot y^{-2}$ is invariant under this group of motions. The invariant Dirichlet integral
for functions is,
$$
\iint (U_x^2+U_y^2)dxdy,
$$
and the Laplace-Beltrami operator associated with this is:
$$
L_0=y^2\cdot\triangle=y^2\cdot(\partial_x^2+\partial_y^2).
$$
A function $f$ defined on the Poincar\'{e} plane $\prod $, that is the upper half plane: $y>0$, $-\infty<x<\infty $,
$w=x+iy$ is called automorphic with respect to $\Gamma $ (a discrete subgroup as above) if $f(\gamma w)=f(w)$
$\forall\,\gamma\in\Gamma $. The Laplace-Beltrami operator $L_0$ maps automorphic functions into automorphic functions.
A fundamental domain $F$ of $\Gamma $ is a sub-domain of the Poincar\'{e} plane such that every point of $\prod $ can
be carried into a point of the closure $\overline{F}$ of $F$ by a transformation in $\Gamma $ and no point of $F$ is
carried into another point of $F$ by such a transformation. $\overline{F}$ can be regarded as a manifold where those
boundary points which can be mapped by a $\gamma\in\Gamma $ to each other, are identified. The restriction of $f$ (automorphic
with respect to $\Gamma $) to the fundamental domain $F$ has to satisfy the above mentioned boundary relations imposed
by $f(\gamma w)=f(w)$. These relations serve as boundary conditions for the operator $L_0$. In fact, they define $L_0$
as a self-adjoint operator acting on $L_2(F)$, the space of functions on $F$ square integrable with respect to the invariant
measure. Our setting is different but similar. We have a non-constant entire function $f(w)$, and its automorphic group
${\rm Aut}(f)=\{\phi_{0n}(w)\}_n$ which is a discrete group. Its elements are defined by all the maximal leaves of
$f^{-1}(f(w))$. The function $f$ is automorphic with respect to the discrete group ${\rm Aut}(f)$. The normal maximal
domains of $f(w)$ are the parallels of the fundamental domains $F$ of $\Gamma $. What can be the parallel of the
Laplace-Beltrami operator $L_0$? It should be an operator that maps automorphic functions with respect to the discrete
group ${\rm Aut}(f)$, into automorphic functions. We already know that the set of all the automorphic functions are the
compositions $h\circ f$ where $h$ is a non-constant function. Thus it is reasonable to take as a parallel to $L_0$ the
right-shift operator induced by $f$, $R_f$, \cite{r1}. It has as its domain of definition the algebra $E$ of all the 
non-constant entire functions and it maps them onto the sub-class $E\circ f=\{h\circ f\,|\,h\in E\}$, i.e. onto all
the automorphic functions with respect to ${\rm Aut}(f)$. Thus:
$$
R_f\,:\,E\rightarrow E\circ f,\,\,\,\,\,R_f(h)=h\circ f.
$$
The operator $R_f$ is not only a linear operator. It is an algebra morphism. For let $h,g\in E$ and let $c\in\mathbb{C}$.
Then:
$$
R_f(c\cdot h)=(c\cdot h)\circ f=(c\cdot h)(f(w))=c\cdot h(f(w))=c\cdot(h\circ f)=c\cdot R_f(h),
$$
and
$$
R_f(g+h)=(g+h)\circ f=(g+h)(f(w))=g(f(w))+h(f(w))=(g\circ f)+(h\circ f)=R_f(g)+R_f(h),
$$
and similarly
$$
R_f(g\cdot h)=(g\cdot h)\circ f=(g\cdot h)(f(w))=g(f(w))\cdot h(f(w))=(g\circ f)\cdot(h\circ f)=R_f(g)\cdot R_f(h).
$$
Another important property of the $f$ right-shift operator is: $R_f$ is an injective mapping on $E$, \cite{r1}. \\
Also: the image of $R_f$, $R_f(E)$ is a closed subset of the topological space $(E,\tau_{cc})$. The topology $\tau_{cc}$
is the topology of compact convergence on the space $E$. For each compact $K\subset\mathbb{C}$ and for each $\epsilon>0$
and each $h\in E$ we define the corresponding open ball centered at $h$ by the standard formula:
$$
B_K(h,\epsilon)=\{ g:\mathbb{C}\rightarrow\mathbb{C}\in E\,|\,|g(z)-h(z)|<\epsilon,\,\forall\,z\in K\}.
$$
The family $\{B_K(h,\epsilon)\,|\,K\subset\mathbb{C}\,{\rm a}\,{\rm compact}, h\in E, \epsilon>0\}$ forms a sub-basis
for the topology $\tau_{cc}$. The sequence $g_n\in E$ converges to the limit $g\in E$ if and only if for restrictions
to compacta we have $g_n|_K\rightarrow g|_K$ uniformly on $K$, $\forall$ compact $K\subset\mathbb{C}$.

\begin{remark}\label{rem26}
The space $(E,\tau_{cc})$ is a path connected space, \cite{r1}. We know that the mapping $R_f\,:\,E\rightarrow E$ is a
continuous and an injective mapping, \cite{r1}. Its image $R_f(E)$ is a closed subset of $(E,\tau_{cc})$. Hence
$R_f(E)$ is also open $\Leftrightarrow$ $R_f(E)=E$. This is equivalent to $f(w)=aw+b$, for some $a\in\mathbb{C}^{\times}$
and some $b\in\mathbb{C}$.
\end{remark}
\noindent
Another interesting property of $R_f$ is $\forall\,f\in E-\{aw+b\,|\,a\in\mathbb{C}^{\times}, b\in\mathbb{C}\}$ we have
the identity $\partial R_f(E)=R_f(E)$. Thus the image $R_f(E)$ is its own boundary. \\
What can be the parallel of the hyperbolic metric? It should be a metric $d_f(\cdot,\cdot)\,:\,\mathbb{C}^2\rightarrow\mathbb{R}_{\ge 0}$
which is invariant under the action of the automorphic group ${\rm Aut}(f)$. We give it with other facts that
were mentioned above in the following:

\begin{theorem}\label{thm14}
Let $f\in E$. Then the parallel of the discrete group $\Gamma $ of hyperbolic motions in the plane is ${\rm Aut}(f)$.
The parallel of the Laplace-Beltrami operator, $L_0$, is the $f$ right-shift operator. The parallel of the hyperbolic
metric $(dx^2+dy^2)\cdot y^{-2}$ is the $f$ path-metric, $d_f(\cdot,\cdot)\,:\,\mathbb{C}^2\rightarrow\mathbb{R}_{\ge 0}$
given by the following formula: $\forall\,a,b\in\mathbb{C}$,
$$
d_f(a,b)=\inf\{l_f(\gamma)\,|\,\gamma:[0,1]\rightarrow\mathbb{C}\,{\rm is}\,{\rm a}\,{\rm smooth}\,{\rm path}\,{\rm from}\,
\gamma(0)=a\,{\rm to}\,\gamma(1)=b\,\,{\rm and}\,
$$
$$
l_f(\gamma)\,{\rm is}\,{\rm the}\,{\rm length}\,{\rm of}\,{\rm the}\,{\rm path}\,f\circ\gamma:[0,1]\rightarrow\mathbb{C}\,
{\rm from}\,(f\circ\gamma)(0)=f(a)\,{\rm to}\,(f\circ\gamma)(1)=f(b)\}.
$$
We further have the following: \\
1) $R_f\,:\,E\rightarrow E\circ f$ is continuous, injective and surjective mapping in the topological space
$(E,\tau_{cc})$. \\
2) $R_f$ is a linear operator on $E$ which also preserves multiplication of functions. \\
3)$R_f(E)$ is closed in $(E,\tau_{cc})$, and $R_f(E)=E$ if and only if $f(w)=aw+b$, for some $a\in\mathbb{C}^{\times}$
and $b\in\mathbb{C}$. \\
4) $\forall\,f\in E-\{aw+b\,|\,a\in\mathbb{C}^{\times},\,b\in\mathbb{C}\}$, $\partial R_f(E)=R_f(E)$.
\end{theorem}
\noindent
{\bf Proof.} \\
The only assertions we need to prove are those related to the metric $d_f$. Namely we need to prove two things: \\
(a) $d_f$ is a metric on $\mathbb{C}$, (b) $\forall\,\phi\in{\rm Aut}(f)$ we have invariance of $d_f$, i.e. $\forall\,a,b\in\mathbb{C}$,
$d_f(a,b)=d_f(\phi(a),\phi(b))$. \\
{\bf A proof of (a):} Let $a,b,c\in\mathbb{C}$, then clearly $d_f(a,b)=d_f(b,a)$ because each path $\gamma(t)$ from
$a$ to $b$ induces the reverse path $\gamma(1-t)$ from $b$ to $a$ and the $f$ images of both are the same hence have
equal length. \\
Also $d_f(a,a)=0$ by using the constant path. Moreover, if $d_f(a,b)=0$ then for each $\epsilon>0$ there is a path $\gamma$
from $a$ to $b$ such that the length of its $f$ image $f\circ\gamma$ is smaller than $\epsilon>0$. Since the entire function $f$
is non-constant we can find a positive but small enough number $r>0$ such that it has the following two properties: \\
(i) Any path $\gamma(t)$ from $a$ to $b$ must intersect the circle $\{w\in\mathbb{C}\,|\,|w-a|=r\}$. \\
(ii) $f^{'}(w)$ has no zero on that circle. \\
Since the circle is compact we have $\min\{|f^{'}(w)|\,|\,|w-a|=r\}=\delta>0$. This implies that the length $l_f(\gamma)$
of the image path $f\circ\gamma$ is bounded from below by some fixed number $m(\delta)>0$ and hence $d_f(a,b)\ge m(\delta>0$
which contradicts the assumption $d_f(a,b)=0$ unless $a=b$. \\
Finally, $d_f$ satisfies the triangle inequality: $d_f(a,c)\le d_f(a,b)+d_f(b,c)$, because the set of paths from $a$ to $c$
through $b$ is a subset of the set of paths from $a$ to $c$. This concludes the proof of (a). \\
{\bf A proof of (b):} We now prove the invariance of the metric $d_f(\cdot,\cdot)$ with respect to the
discrete group ${\rm Aut}(f)$. This follows directly from the definitions of $d_f(\cdot,\cdot)$ and of ${\rm Aut}(f)$, namely:
$$
d_f(a,b)=\inf\{l_f(\gamma)\,|\,\gamma:[0,1]\rightarrow\mathbb{C}\,{\rm is}\,{\rm a}\,{\rm smooth}\,{\rm path}\,{\rm from}\,
\gamma(0)=a\,{\rm to}\,\gamma(1)=b\,\,{\rm and}\,
$$
$$
l_f(\gamma)\,{\rm is}\,{\rm the}\,{\rm length}\,{\rm of}\,{\rm the}\,{\rm path}\,f\circ\gamma:[0,1]\rightarrow\mathbb{C}\,
{\rm from}\,(f\circ\gamma)(0)=f(a)\,{\rm to}\,(f\circ\gamma)(1)=f(b)\}=
$$
$$
=\inf\{l_f(\gamma)\,|\,\gamma:[0,1]\rightarrow\mathbb{C}\,{\rm is}\,{\rm a}\,{\rm smooth}\,{\rm path}\,{\rm from}\,
\gamma(0)=\phi(a)\,{\rm to}\,\gamma(1)=\phi(b)\,\,{\rm and}\,
$$
$$
l_f(\gamma)\,{\rm is}\,{\rm the}\,{\rm length}\,{\rm of}\,{\rm the}\,{\rm path}\,f\circ\gamma:[0,1]\rightarrow\mathbb{C}\,
{\rm from}\,f(\phi(a))=f(a)\,{\rm to}\,
$$
$$
f(\phi(b))=f(b)\}=d_f(\phi(a),\phi(b)).\,\,\,\,\,\,\,\,\,\,\,\,\,\,\,\,\,\,\qed
$$
\\
The next obvious step is to look for the eigenvalues of the operator $R_f$. If we think of the $f$ right composition operator
$R_f$ as a possible parallel of the Laplace-Beltrami operator, then it is natural to require about its eigenvalues,
eigenfunctions and maybe try to come up with a kind of a Selberg-trace formula. \\
We recall that in our setting the underlying linear space, $E$ contains all of the non-constant entire functions,
and the right composition operators $R_f$ of interest are those for which $f(w)$ is not an entire automorphism, i.e.
$f(w)\ne a\cdot w+b$, $\forall\,a\in\mathbb{C}^{\times}$ and $\forall\,b\in\mathbb{C}$. It turns out that the
result of this search is disappointing because the order of growth of entire functions imposes (within the class
of functions of our interest) too much rigidity. Looking in the corresponding defining equation of the eigenvalues
for $R_f$ leads us to ask for which values of $\lambda\in\mathbb{C}$, the operator $R_f-\lambda\cdot I$ is not
an invertible operator.Thus if we look for entire functions $h$ that satisfy $(R_f-\lambda\cdot I)(h)=0$, that is
$h(f(w))\equiv \lambda h(w)$, then because $f$ is not affine and due to Theorem A.2. of Clunie, on page 207 of \cite{cy}, and
to Corollary A.1. of P\"{o}lya, on page 208 of \cite{cy}, we deduce that the last equation can have only
constant solutions $9w)$. If the value of the constant is not zero, then $\lambda=1$, and if $h(w)\equiv 0$, then
$\lambda $ can be any complex number. However, non-zero constant function do not belong to $E$ by its definition.
This takes care of those cases in which the order of the growth of $h\circ f$ is larger than that of $h$. \\
If we adopt a different definition of eigenvalues, and we are interested in those $\lambda $ for which $R_f-\lambda\cdot I$
is not injective on $E$, then we are led to consider the situation where $g, h\in E$ and: $(R_f-\lambda\cdot I)(h)=
(R_f-\lambda\cdot I)(g)$, that is $h(f(w))-g(f(w))\equiv\lambda\cdot (h(w)-g(w))$, and again we deduce that
$h(w)-g(w)$ must be a constant. If the constant difference between $h(w)$ and $g(w)$ is not zero, then $\lambda=1$,
and, of course if $h(w)\equiv g(w)$ then $\lambda $ can be any complex number.

\begin{remark}\label{rem27} 1) An equation of the form $h(f(w))\equiv\lambda h(w)$ resembles very much the equation
that determines the automorphic functions of $h(w)$. When $\lambda=1$ it is exactly that equation. By Theorem 10 of
Shimizu, on page 237 of \cite{s2}, the only possible automorphic functions which are also entire functions, are
affine functions of a very special kind: $e^{\theta\pi}\cdot w+b$, where $\theta\in\mathbb{Q}$, and where $b\in\mathbb{C}$.
This is consistent with our findings prior to this remark, \ref{rem27}.
\end{remark}

\begin{remark}\label{rem28}
We know that the operator $R_f$ is injective. Hence $R_f-\lambda\cdot I$ is injective for $\lambda=0$.
\end{remark}

\section{Local groups}\label{sec27}

A good source for the theory of topological groups is the relatively new book \cite{at}.
It will be convenient to have the notion of local groups handy in our setting of the automorphic group. This is
because the automorphic functions of an entire function are generically multivalued. Thus they naturally are
defined and uniform on the corresponding Riemann surfaces. It can be useful occasionally to restrict them
to a leaf, i.e. to a sub-domain of the complex plane whose complimentary set contains no continuum. In those cases
the different restricted automorphic functions might be defined on different sub-domains of the plane, that
differ by small sets. Thus we might need the notion of a local group. We refer to \cite{tao} to definition 2.1.1
on page 26. In that book the need in local groups arise because the connection between Lie groups and Lie
algebras are local. the only portion of the Lie group which is of importance in that respect is the portion that is
close to the group identity $1$. We will adopt here the notions and the notations from \cite{tao}.

\begin{definition}\label{def4}
A local topological group $G=(G,\Omega,\Lambda,1,\cdot,()^{-1})$, or a local group for short, is a topological
space $G$, equipped with an identity element $1\in G$, a partially defined but continuous multiplication
operation $\cdot\,:\,\Omega\rightarrow G$ for some domain $\Omega\subseteq G\times G$, and a partially defined
but continuous inversion operation $()^{-1}:\,\Lambda\rightarrow G$, where $\Lambda\subseteq G$, obeying the
following axioms: \\
(1) $\Omega$ is an open neighborhood of $G\times\{1\}\cup\{1\}\times G$, and $\Lambda $ is an open neighborhood
of $1$. \\
(2) If $g,h,k\in G$ are such that $(g\cdot h)\cdot k$ and $g\cdot (h\cdot k)$ are both well-defined in $G$, then
they are equal. \\
(3) For all $g\in G$, $g\cdot 1=1\cdot g=g$. \\
(4) If $g\in G$ and $g^{-1}$ are well-defined in $G$, then $g\cdot g^{-1}=g^{-1}\cdot g=1$.
\end{definition}
\noindent
A local group is said to be symmetric if $\Lambda=G$, i.e., if every element of $G$ has an inverse $g^{-1}$ that is also 
in $G$. Clearly, every topological group is a local group. This the reason that sometimes the former are called global
topological groups. A model class of examples of a local group comes from restricting a global group to an open neighborhood
of the identity. Here is the definition from \cite{tao}:

\begin{definition}\label{def5}
If $G$ is a local group, and $U$ is an open neighborhood of the identity in $G$, then we define the restriction 
$G|_U$ of $G$ to $U$ to be the topological space $U$ with domains:
$$
\Omega|_U:=\{(g,h)\in\Omega\,|\,g,h,g\cdot h\in U\}\,\,{\rm and}\,\,\Lambda|_U:=\{g\in\Lambda\,|\,g,g^{-1}\in U\}
$$
and with the group operations $\cdot,\,()^{-1}$ being the restriction of the group operations of $G$ to $\Omega|_U$,
$\Lambda|_U$ respectively. If $U$ is symmetric (in the sense that $g^{-1}$ is well-defined and lies in $U$, 
$\forall\,g\in U$, then this restriction $G|_U$ will also be symmetric. Sometimes the notation is abused and one
refers to the local group $G|_U$ simply as $U$.
\end{definition}

\begin{remark}\label{rem29}
The natural question as to whether every local group arises as the restriction of a global group, is not simple. The answer
can be vaguely summarized as "essentially yes in certain circumstances, but not in general". That is from
\cite{tao}.
\end{remark}
\noindent
Pushing forward a topological group via a homeomorphism near the identity: Let $G$ be a global or local group
and let $Phi\,:\,U\rightarrow V$ be a homeomorphism from a neighborhood $U$ of the identity in $G$ to a neighborhood
$V$ of the origin $0$ in $\mathbb{R}^d$, such that $Phi(1)=0$. Then we can define a local group $Phi_{*}G|_U$
which is the set $V$ (viewed as a submanifold of $\mathbb{R}^d$) with the local group identity $0$, the local
group multiplication law $*$ defined by the formula: $x*y=\Phi(\Phi^{-1}(x)\cdot\Phi^{-1}(y))$ which is defined
whenever $\Phi^{-1}(x), \Phi^{-1}(y), \Phi^{-1}(x)\cdot\Phi^{-1}(y)$ are well-defined and lie in $U$, and the
local group inversion law $()^{*-1}$ defined by the formula: $x^{*-1}=\Phi(\Phi^{-1}(x)^{-1})$, defined whenever
$\Phi^{-1}(x)$, $\Phi^{-1}(x)^{-1}$ are well-defined and lie in $U$. One easily verifies that $\Phi_{*}G|_U$ is
a local group. Sometimes this group is denoted by $(V,*)$. It is different from the additive local group
$(V,+)$ arising by the restriction of $(\mathbb{R}^d,+)$ to $V$. \\
Next we generalize the notion of homomorphism.

\begin{definition}\label{def6}
A continuous homomorphism $\Phi\,:\,G\rightarrow H$ between two local groups $G, H$ is a continuous map
from $G$ to $H$ with the following properties: \\
(i) $\Phi(1_G)=1_H$, $1_G$ is the neutral element of $G$. \\
(ii) If $g\in G$ is such that $g^{-1}$ is well-defined in $G$, then $\Phi(g)^{-1}$ is well-defined in $H$,
and $\Phi(g)^{-1}=\Phi(g^{-1})$. \\
(iii) If $g,h\in G$ are such that $g\cdot h$ is well-defined in $G$, then $\Phi(g)\cdot\Phi(h)$ is well-defined
in $H$ and $\Phi(g)\cdot\Phi(h)=\Phi(g\cdot h)$.
\end{definition}

\section{The sums of the $k$'th derivatives of all the elements of the automorphic group ${\rm Aut}(f)$,
for any $f\in E$ of order $0<\rho<\frac{1}{2}$, $k=1,2,3,\ldots$}\label{sec28}

\begin{theorem}\label{thm15}
Let $f\in E$ have a positive order $\rho$, which is smaller than $\frac{1}{2}$, i.e. $0<\rho<\frac{1}{2}$.
Let ${\rm Aut}(f)=\{\phi_{0n}\}_n$. Then there exists a sequence of positive numbers, tending to infinity:
$R_1<R_2<\ldots<R_n<\ldots (R_n\rightarrow\infty)$, such that $\forall\,k\in\mathbb{Z}^+$ we have the following
identities:
$$
\lim_{j\rightarrow\infty}\sum_{|\phi_{0n}(w)|<R_j}\phi_{0n}^{(k)}(w)\equiv 0,\,\,\,\,\,\,\,\forall\,w\in\mathbb{C}.
$$
Equivalently:
$$
\lim_{j\rightarrow\infty}\frac{d^k}{dw^k}\left\{\sum_{|\phi_{0n}(w)|<R_j}\phi_{0n}(w)\right\}\equiv 0,\,\,\,\,\,\,\,
\forall\,w\in\mathbb{C}.
$$
\end{theorem}
\noindent
{\bf Proof.}
We proved that $\forall\,f\in E$, regardless of its order $\rho$, we have the following infinitely many identities: \\
$\forall\,R>0$ such that $\{|w|=R\}\cap\{\phi_{0n}(z)\}_n=\emptyset$,
\begin{equation}\label{eq15}
\sum_{|\phi_{0n}(z)|<R}\phi_{0n}^{'}(z)=\left(\frac{f'(z)}{2\pi i}\right)\oint_{|w|=R}\frac{dw}{f(w)-f(z)},
\end{equation}
$$ $$ $$ $$ $$ $$
\begin{equation}\label{eq16}
\sum_{|\phi_{0n}(z)|<R}\phi_{0n}^{''}(z)=\left(\frac{f''(z)}{2\pi i}\right)\oint_{|w|=R}\frac{dw}{f(w)-f(z)}+
\end{equation}
$$
+\left(\frac{f'(z)^2}{2\pi i}\right)\oint_{|w|=R}\frac{dw}{(f(w)-f(z))^2},
$$ 
$$ $$ $$ $$ $$ $$
\begin{equation}\label{eq17}
\sum_{|\phi_{0n}(z)|<R}\phi_{0n}^{'''}(z)=\left(\frac{f'''(z)}{2\pi i}\right)\oint_{|w|=R}\frac{dw}{f(w)-f(z)}+
\end{equation}
$$
+\left(\frac{3f'(z)f''(z)}{2\pi i}\right)\oint_{|w|=R}\frac{dw}{(f(w)-f(z))^2}+
\left(\frac{f'(z)^3}{2\pi i}\right)\oint_{|w|=R}\frac{dw}{(f(w)-f(z))^3},
$$
etc... . In the identities above, the radius $R$ of the circle of integration, is such that for a fixed value
of the parameter $z$, we have $f(w)-f(z)\ne 0$ $\forall\,|w|=R$. In other words the radius $R$ is chosen so
that the circle of integration $|w|=R$ does not contain any element of the ${\rm Aut}(f)$-orbit of $z$, $\{\phi_{0n}(z)\}_n$.
Also, since we used the one dimensional Weierstrass product representation of $f(w)-f(z)$, the fixed value of the
parameter $z$ belongs to $\mathbb{C}-f^{-1}(f(0))$. We will prove that there exists a sequence of positive
numbers tending to infinity: $R_1<R_2<R_3<\ldots<R_n<\ldots (R_n\rightarrow\infty)$ such that 
$\lim_{j\rightarrow\infty}\sum_{|\phi_{0n}(z)|<R_j}\phi_{0n}^{'}(z)=0$ for any value of the parameter $z\not\in f^{-1}(f(0))$.
It will follow by continuity arguments that the identity is also true for the "forbidden values" of $z$, i.e. on the
fiber $f^{-1}(f(0))$. Also it will be clear that a similar proof will apply for the second identity
$\lim_{j\rightarrow\infty}\sum_{|\phi_{0n}(z)|<R_j}\phi_{0n}^{''}(z)=0$ $\forall\,w\in\mathbb{C}$, and another
similar proof will apply for the third identity $\lim_{j\rightarrow\infty}\sum_{|\phi_{0n}(z)|<R_j}\phi_{0n}^{'''}(z)=0$
$\forall\,w\in\mathbb{C}$, and all with the same sequence of positive numbers $0<R_1<R_2<R_3<\ldots<R_n<\ldots 
(R_n\rightarrow\infty)$. That procedure of giving a separate proof for each identity will save us arguments about
differentiations under the limit operator and inquiring whether formulas like this:
$$
\lim_{j\rightarrow\infty}\sum_{|\phi_{0n}(z)|<R_j}\phi_{0n}^{''}(z)=\left(
\lim_{j\rightarrow\infty}\sum_{|\phi_{0n}(z)|<R_j}\phi_{0n}^{'}(z)\right)^{'},\,\,\,{\rm etc...}
$$
are valid. Up to this point restricting the value of the order $\rho$ seem not come up. However, in order to prove that
we have $\lim_{j\rightarrow\infty}\sum_{|\phi_{0n}(z)|<R_j}\phi_{0n}^{'}(z)\equiv 0$ $\forall\,z\in\mathbb{C}-f^{-1}(f(0))$
we will need the assumption $0<\rho<\frac{1}{2}$. That assumption will enable us to make a use of the $\cos\pi\rho$-Theorem
of Wiman, where the power $\cos\pi\rho$ in the inequality of this theorem satisfies $0<\cos\pi\rho<1$, by the assumption
$0<\rho<\frac{1}{2}$. The value $\rho=0$ is taken out due to another argument we need (an elementary one). A classical
reference to the Wiman's $\cos\pi\rho$-Theorem is in the book \cite{l}, Theorem 30 on page 72. Using Wiman's $\cos\pi\rho$-Theorem,
we conclude that there exists a sequence $r_1<r_2<r_3<\ldots<r_n<\ldots (r_n\rightarrow\infty)$ such that for arbitrary
$\epsilon>0$ and $n>n_{\epsilon}$ we have:
$$
m_f(r_n)>\left(M_f(r_n)\right)^{\cos\pi\rho-\epsilon},
$$
Where
$$
m_f(r)=\min_{0\le\theta<2\pi}|f(re^{i\theta})|\,\,\,{\rm and}\,\,\,M_f(r)=\max_{0\le\theta<2\pi}|f(re^{i\theta})|.
$$
One can get many more (uncountably many on each radius) such sequences $r_1<r_2<r_3<\ldots<r_n<\ldots $ by
perturbations. We recall that our first identity to be used is equation (\ref{eq15}). We are going to prove that:
$$
\lim_{j\rightarrow\infty}\left(\frac{f'(z)}{2\pi i}\right)\oint_{|w|=R_j}\frac{dw}{f(w)-f(z)}=0,
$$
for an appropriate Wiman's-type of a sequence $R_1<R_2<R_3<\ldots<R_n<\ldots (R_n\rightarrow\infty)$. This will imply
that we have:
$$
\lim_{j\rightarrow\infty}\sum_{|\phi_{0n}(z)|<R_j}\phi_{0n}^{'}(z)\equiv 0\,\,\,\,\,\forall\,z\in\mathbb{C}-f^{-1}(f(0)).
$$
Using the flexibility in choosing Wiman's-type sequences (invoking perturbations), we are going to choose sequences
of radii $\{R_j\}_j$ that satisfy the requirements $\{|w|=R_j\}\cap\{\phi_{0n}\}_n=\emptyset$, as well as the conclusion
of the Wiman's $\cos\pi\rho$-Theorem, namely that:
$$
m_f(R_j)>(M_f(R_j))^{\cos\pi\rho-\epsilon}.
$$
By the triangle inequality we have:
$$
\left|\left(\frac{f'(z)}{2\pi i}\right)\oint_{|w|=R_j}\frac{dw}{f(w)-f(z)}\right|\le
\left(\frac{|f'(z)|}{2\pi}\right)\oint_{|w|=R_j}\frac{|dw|}{|(M_f(R_j))^{\cos\pi\rho}-|f(z)||},
$$
for a large enough $j$. On the circle $|w|=R_j$ we have $w=R_je^{i\theta}$, $dw=iR_je^{i\theta}d\theta$, $|dw|=R_j d\theta$,
$0\le\theta<2\pi$. We conclude that:
$$
\left|\left(\frac{f'(z)}{2\pi i}\right)\oint_{|w|=R_j}\frac{dw}{f(w)-f(z)}\right|\le
\frac{|f'(z)|R_j}{|(M_f(R_j))^{\cos\pi\rho}-|f(z)||}.
$$
We recall that for a large $j$ we have: $M_f(R_j)\approx e^{R_j^{\rho}}$, and since we have both $\rho>0$ and $0<\cos\pi\rho$, 
by our assumption on the order, $0<\rho<\frac{1}{2}$, we conclude that:
$$
\lim_{j\rightarrow\infty}\frac{|f'(z)|R_j}{|(M_f(R_j))^{\cos\pi\rho}-|f(z)||}=0.
$$
Hence we proved our first identity on the restricted domain of the values of the parameter $z$:
$$
\lim_{j\rightarrow\infty}\sum_{|\phi_{0n}(z)|<R_j}\phi_{0n}^{'}(z)\equiv 0\,\,\,\,\,\forall\,z\in\mathbb{C}-f^{-1}(f(0)).
$$
Now the theorem follows as explained above. $\qed $ \\
\\
\begin{remark}\label{rem30}
Theorem \ref{thm15} deals with the sums $\sum_{|\phi_{0n}(w)|<R_j}\phi_{0n}^{(k)}(w)$, $\forall\,w\in\mathbb{C}$
for values of $k$ which are natural numbers. There is no claim for the value $k=0$, i.e. 
$\sum_{|\phi_{0n}(w)|<R_j}\phi_{0n}(w)$, $\forall\,w\in\mathbb{C}$. One might be expecting naively that
$\lim_{j\rightarrow\infty}\sum_{|\phi_{0n}(w)|<R_j}\phi_{0n}(w)\equiv{\rm Const.}$, $\forall\,w\in\mathbb{C}$, but
this turns out to be false also within the restricted domain of the order $\rho $.
\end{remark}
\noindent
{\bf Examples:} The first two examples are of entire functions of order $\rho $, where $\rho$ is off the domain $(0,\frac{1}{2})$. \\
\\
(1) We consider the exponential function $f(w)=e^w$. Here $\rho=1$ and ${\rm Aut}(e^w)=\{w+2\pi in\}_{n\in\mathbb{Z}}$. Thus
$\phi_{0n}^{'}(w)=1$ and hence we have 
$$
\sum_{|\phi_{0n}(w)|<R}\phi_{0n}^{'}(w)=|\{n\in\mathbb{Z}\,|\,|w+2\pi in|<R\}|\ge\left[\frac{R-|w|}{2\pi}\right]
\rightarrow_{R\rightarrow\infty}\infty.
$$
Thus clearly for the conclusions of Theorem \ref{thm15} to hold, some assumptions on the order $\rho $ are needed. \\
\\
(2) Let us consider the case of a non-constant polynomial:
$$
P_N(w)=a_Nw^N+\ldots+a_1w+a_0,\,\,\,\,\,\,\,\,\,\,a_N\ne 0,\,\,\,N\ge 1.
$$
Here $\rho=0$ and ${\rm Aut}(P_N(w))=\{w,\ldots,\phi_{N-1}(w)\}$. Hence:
$$
P_N(w)-P_N(z)=a_N(w-z)\cdot\ldots\cdot(w-\phi_{N-1}(z))=
$$
$$
=a_N(w^N-(z+\ldots+\phi_{N-1}(z))w^{N-1}+\ldots+(-1)^N\cdot(z\cdot\ldots\cdot\phi_{N-1}(z))).
$$
We conclude that:
$$
z+\ldots+\phi_{N-1}(z)=\left\{\begin{array}{cll} -\frac{a_{N-1}}{a_N} & {\rm if} & N>1 \\
z & {\rm if} & N=1\end{array}\right..
$$
Thus in the simplest case $N=1$, $\sum_{|\phi_{0n}(z)|<R}\phi_{0n}^{'}(z)\rightarrow_{R\rightarrow\infty}1$. \\
\\
Next we give an example that shows that the domain $0<\rho<\frac{1}{2}$ is sharp for the conclusion of Theorem \ref{thm15}
to be valid. We already know that necessarily $0<\rho$ (Example (2), the case $N=1$). We now show that
$\rho=\frac{1}{2}$ is out of the admissible domain. \\
\\
(3) Let $f(w)=\cos\sqrt{w}=\frac{1}{2}(e^{i\sqrt{w}}+e^{-i\sqrt{w}})$. Since $\cos w$ is an even entire
function, it follows that $f(w)$ is an entire function. We have $|f(w)|\le\frac{1}{2}(|e^{i\sqrt{w}}|+
|e^{-i\sqrt{w}}|)\le\frac{1}{2}(e^{\sqrt{|w|}}+e^{\sqrt{|w|}})=e^{\sqrt{|w|}}$. If $r>0$, then:
$f(-r)=\frac{1}{2}(e^{-\sqrt{r}}+e^{\sqrt{r}})\ge\frac{1}{2}e^{\sqrt{r}}$. The two inequalities above
prove that $\rho=\frac{1}{2}$. We next compute ${\rm Aut}(\cos\sqrt{w})$. For that purpose we first
solve for $\theta$ the following equation: $e^{i\theta}+e^{-i\theta}=e^{i\psi}+e^{-i\psi}$. This gives
the quadratic $e^{2i\theta}-(e^{i\psi}+e^{-i\psi})e^{i\theta}+1=0$. We obtain the following solution:
$$
e^{i\theta}=\left\{\begin{array}{lll} e^{i\psi} & {\rm if} & + \\ e^{-i\psi} & {\rm if} & -\end{array}\right..
$$
Hence:
$$
\left\{\begin{array}{lll} i\theta & = & i\psi+2\pi ik \\ i\theta & = & -i\psi+2\pi ik, k\in\mathbb{Z}
\end{array}\right.
$$
So:
$$
\left\{\begin{array}{lll} \theta & = & \psi+2\pi k \\ \theta & = & -\psi+2\pi k, k\in\mathbb{Z}
\end{array}\right.
$$
In our case $\cos\sqrt{w}=\cos\sqrt{z}$, so $\theta=\sqrt{w}$, $\psi=\sqrt{z}$. Thus: $\sqrt{w}=\pm\sqrt{z}+2\pi k$,
$k\in\mathbb{Z}$. Squaring: $\phi_k(z)=z+4\pi k\sqrt{z}+a\pi^2k^2,\,k\in\mathbb{Z}$. We proved that:
$$
{\rm Aut}(\cos\sqrt{w})=\left\{w+4\pi k\sqrt{w}+4\pi^2k^2\,|\,k\in\mathbb{Z}\right\}.
$$
In particular we have $f^{-1}(f(0))=\{4\pi^2k^2\,|\,k\in\mathbb{Z}\}$. Also we have:
$$
\phi_k^{'}(z)=1+\frac{2\pi k}{\sqrt{z}},\,\,\,\,\,\phi_k^{''}(z)=-\frac{\pi k}{z\sqrt{z}},\ldots .
$$
So the identities: $\sum_{|\phi_k(z)|<R}\phi_k^{(j)}(z)\rightarrow_{R\rightarrow\infty} 0,\,\,\,j\ge 2$, are
equivalent to:
$$
\lim_{R\rightarrow\infty}\sum_{|z+4\pi k\sqrt{z}+4\pi^2k^2|<R}k=0.
$$
The identity: $\sum_{|\phi_k(z)|<R}\phi_k^{'}(z)\rightarrow_{R\rightarrow\infty} 0$ is equivalent to:
$$
\lim_{R\rightarrow\infty}\sum_{|z+4\pi k\sqrt{z}+4\pi^2k^2|<R}\left(1+\frac{2\pi k}{\sqrt{z}}\right)=0.
$$
But this last identity is not consistent with the identities that correspond to $j\ge 2$. For the two identities 
together imply that:
$$
lim_{R\rightarrow\infty}\sum_{|z+4\pi k\sqrt{z}+4\pi^2k^2|<R}1\equiv 0,
$$
which is clearly false. \\
\\
(4) Our last example is a straight forward application of Theorem \ref{thm15}. \\
Let $f(w)=\frac{1}{2}(\cos w^{1/4}+\cosh w^{1/4})$. Using the identities: $\cos w=\frac{1}{2}(e^{iw}+e^{-iw})$
and $\cosh w=\frac{1}{2}(e^w+e^{-w})$, we obtain the following power series representation of $f(w)$:
$$
f(w)=\sum_{k=0}^{\infty}\frac{w^k}{(4k)!}.
$$
This shows that $f(w)$ is an entire function. Next, by the triangle inequality we obtain:
$$
|f(w)|\le\frac{1}{2}\left(\frac{1}{2}\left(|e^{iw^{1/4}}|+|e^{-iw^{1/4}}|\right)
+\frac{1}{2}\left(|e^{w^{1/4}}|+|e^{-w^{1/4}}|\right)\right)\le
$$
$$
\le\frac{1}{4}\left(e^{|w|^{1/4}}+e^{|w|^{1/4}}+e^{|w|^{1/4}}+e^{|w|^{1/4}}\right)=e^{|w|^{1/4}}.
$$
Also, for $r>0$ large enough we have:
$$
|f(r)|=\frac{1}{2}\left|\cos r^{1/4}+\frac{e^{r^{1/4}}+e^{-r^{1/4}}}{2}\right|\ge\frac{1}{2}\left|-1+\frac{1}{2}e^{r^{1/4}}\right|
\ge\frac{1}{8}e^{r^{1/4}}.
$$
The last two inequalities imply that $\rho=\frac{1}{4}$ for our entire function $f(w)$ and so Theorem \ref{thm15} applies.
Thus if for a fixed $z$ the solutions of the following equation in the unknown $w$:
$$
\cos w^{1/4}+\cosh w^{1/4}=\cos z^{1/4}+\cosh z^{1/4},
$$
are given by $\{w\}=\{\phi_{0n}(z)\}_n$, then $\forall\,k\in\mathbb{Z}^+$ we have:
$$
\lim_{j\rightarrow\infty}\sum_{|\phi_{0n}(z)|<R_j}\sum_{0n}^{(k)}(z)\equiv 0,\,\,\,\,\,\forall\,z\in\mathbb{C},
$$
for some Wiman's type sequence $0<R_1<R_2<R_3<\ldots<R_n<\ldots (R_n\rightarrow\infty)$. The task of actually
computing the automorphic functions $\phi_{0n}(z)$ for this function $f(w)$, is probably not an easy task.

\section{The circular density of the orbits of the automorphic group ${\rm Aut}(f)$, for any $f\in E$ of
order $0<\rho<\frac{1}{2}$}\label{sec29}

\begin{theorem}\label{thm16}
Let $f\in E$ have a positive order $\rho$, which is smaller than $\frac{1}{2}$, i.e. $0<\rho<\frac{1}{2}$. Let
${\rm Aut}(f)=\{\phi_{0n}\}_n$. For any $z\in\mathbb{C}$ we arrange the ${\rm Aut}(f)$-orbit of $z$,
$Z(f(w)-f(z))=\{\phi_{0n}(z)\}_n$ in a non-decreasing order of the moduli: $|\phi_{00}(z)|\le|\phi_{01}(z)|\le
|\phi_{02}(z)|\le\ldots $ and for any $r$ satisfying $|\phi_{0n}(z)|\le r\le|\phi_{0,n+1}(z)|$ we denote the
corresponding index $n=n(r,z)$. If $|\phi_{0n}(z)|=|\phi_{0,n+1}(z)|$ we may denote anything within reason, for
example $n(r,z)=n\,\,{\rm or}\,\,n+1$. Then we have the following asymptotic circular (or radial, if one prefers)
density estimate:
$$
\lim_{r_j\rightarrow\infty}\left(\left(n(r_j,z)+1\right)\log r_j-r_j^{\rho}\cos\pi\rho\right)=\infty,
$$
for some Wiman's sequence $\{r_j\}_j$ and $\forall\,\mathbb{C}$.
\end{theorem}
\noindent
{\bf Proof.} \\
We have proved (see section \ref{sec12}) that if $|\phi_{0n}(z)|\le r\le|\phi_{0,n+1}(z)|$, then by the Jensen's Theorem
applied to $f(w)-f(z)$, where as usual $f(0)-f(z)\ne 0$, we have:
$$
\left|\prod_{j=0}^n\phi_{0j}(z)\right|=r^{n+1}|f(0)-f(z)|\exp\left\{-\frac{1}{2\pi}\int_0^{2\pi}
\log\left|f(re^{i\theta})-f(z)\right|d\theta\right\}.
$$
Using Wiman's $\cos\pi\rho$-Theorem, as in the proof of Theorem \ref{thm15}, we conclude that there exists a
sequence $r_1<r_2<r_3<\ldots<r_n<\ldots (r_n\rightarrow\infty)$ such that for arbitrary $\epsilon>0$ and $n>n_{\epsilon}$
we have:
$$
m_f(r_n)>\left(M_f(r_n)\right)^{\cos\pi\rho-\epsilon}.
$$
By the triangle inequality we have for any $r_j>0$ that satisfies $(f(r_je^{i\theta})-f(z))(M_f(r_j)^{\cos\pi\rho}-|f(z)|)\ne 0$,
the following estimate:
$$
\log|f(r_je^{i\theta})-f(z)|\ge\log|m_f(r_j)-|f(z)||\ge\log|M_f(r_j)^{\cos\pi\rho}-|f(z)||.
$$
Hence:
$$
\exp\left\{-\frac{1}{2\pi}\int_0^{2\pi}\log\left|f(r_je^{i\theta})-f(z)\right|d\theta\right\}\le
\exp\left\{-\log|M_f(r_j)^{\cos\pi\rho}-|f(z)||\right\}.
$$
Since we have:
$$
M_f(r_j)^{\cos\pi\rho}-|f(z)|\approx e^{r_j^{\rho}\cos\pi\rho},
$$
for a large enough $r_j$, we conclude that for some $c=c(z)$, depending only on $z$, we have for
$|\phi_{0,n(r_j,z)}(z)|\le r_j\le|\phi_{0,n(r_j,z)+1}(z)|$, the following estimate:
$$
\left|\prod_{k=0}^{n(r_j,z)}\phi_{0k}(z)\right|\le r_j^{n(r_j,z)+1}e^{-r_j^{\rho}\cos\pi\rho}\cdot c.
$$
But $0<\rho$, so $f(w)$ is transcendental and hence:
$$
\lim_{r_j\rightarrow\infty}\left|\prod_{k=0}^{n(r_j,z)}\phi_{0k}(z)\right|=+\infty.
$$
Thus:
$$
\lim_{r_j\rightarrow\infty} r_j^{n(r_j,z)+1}e^{-r_j^{\rho}\cos\pi\rho}=+\infty.
$$
Since we have the elementary identity:
$$
r_j^{n(r_j,z)+1}e^{-r_j^{\rho}\cos\pi\rho}=e^{(n(r_j,z)+1)\log r_j-r_j^{\rho}\cos\pi\rho},
$$
it follows that:
$$
\lim_{r_j\rightarrow\infty} \left(\left(n(r_j,z)+1\right)\log r_j-r_j^{\rho}\cos\pi\rho\right)=+\infty\,\,\,\,\,\qed
$$

\section{The Vieta formulas for ${\rm Aut}(f)$, $f\in E$ of order $0\le\rho<1$}\label{sec30}

For $f\in E$ of low order $\rho $, i.e. $0\le\rho<1$, the formulas for the symmetric functions of the reciprocals
of the automorphic functions of $f$ can be derived algebraically as easy as for polynomials. The reason is the fact
that for those orders the Weierstrass canonical representations are exactly the factorization of $f(w)-f(z)$, because
the Weierstrass factors reduces to the simplest form $(1-u)$. This follows by the fact that there is no need in the
Weierstrass auxiliary exponentials that cause the convergence of the infinite product. In low order, the infinite
product converges automatically  already at the level of $(1-u)$.

\begin{theorem}\label{thm17}
Let $f(w)=\sum_{n=0}^{\infty} a_nw^n\in E$ be of order $\rho $, where $0\le\rho<1$. Let the automorphic group of $f$ be
${\rm Aut}(f)=\{\phi_{0n}\}_n$. Then for any $z\in\mathbb{C}-f^{-1}(f(0))$ and for any $n\in\mathbb{Z}^+$ we have:
$$
a_n=(-1)^n\left(f(0)-f(z)\right)\sum_{0\le i_1<i_2<\ldots<i_n}\left(\prod_{j=1}^n\phi_{0i_j}(z)\right)^{-1}.
$$
\end{theorem}
\noindent
{\bf Proof.} \\
The assumption of low order, $0\le\rho<1$ implies that:
$$
f(w)-f(z)=(f(0)-f(z))\prod_n\left(1-\frac{w}{\phi_{0n}(z)}\right)=
$$
$$
=(f(0)-f(z))\sum_{n=0}^{\infty}(-1)^n\sum_{0\le i_1<i_2<\ldots<i_n}\left(\prod_{j=1}^n\phi_{0i_j}(z)\right)^{-1}.\,\,\,\,\,\qed
$$
{\bf Example:} We computed previously, in example 3 after Theorem \ref{thm15}, that for $f(w)=\cos\sqrt{w}\in E$, we have 
$\rho=\frac{1}{2}$, ${\rm Aut}(\cos\sqrt{w})=\{w+4\pi k\sqrt{w}+4\pi^2k^2\,|\,k\in\mathbb{Z}\}$, 
$\cos\sqrt{w}=\sum_{n=0}^{\infty}(-1)^n\frac{w^n}{(2n)!}$, $a_n=\frac{(-1)^n}{(2n)!}$. The simplest case of Theorem \ref{thm17}
is the case $n=1$:
$$
\frac{-1}{2!}=(-1)^1\left(\cos\sqrt{0}-\cos\sqrt{z}\right)\left(\frac{1}{z}+\frac{1}{z+4\pi\sqrt{z}+4\pi^2}+
\frac{1}{z-4\pi\sqrt{z}+4\pi^2}+\right.
$$
$$
\left.+\frac{1}{z+4\pi\cdot 2\sqrt{z}+4\pi^2 2^2}+\frac{1}{z-4\pi\cdot 2\sqrt{z}+4\pi^2 2^2}+\ldots\right)=
$$
$$
=-\left(1-\cos\sqrt{z}\right)\left(\frac{1}{z}+\sum_{k=1}^{\infty}\left(\frac{1}{z+4\pi k\sqrt{z}+4\pi^2k^2}+
\frac{1}{z-4\pi k\sqrt{z}+4\pi^2k^2}\right)\right)=
$$
$$
=-\left(1-\cos\sqrt{z}\right)\left(\frac{1}{z}+2\sum_{k=1}^{\infty}\frac{z+4\pi^2k^2}{(z-4\pi^2k^2)^2}\right).
$$
Thus we obtain the following identity:
$$
\frac{1}{2(1-\cos\sqrt{z})}=\frac{1}{z}+2\sum_{k=1}^{\infty}\frac{z+4\pi^2k^2}{(z-4\pi^2k^2)^2}.
$$
For instance, by substituting $z=\pi^2$ we obtain:
$$
\frac{1}{4}=\frac{1}{\pi^2}+2\sum_{k=1}^{\infty}\frac{\pi^2+4\pi^2k^2}{(\pi^2-4\pi^2k^2)^2},
$$
$$
\frac{\pi^2}{4}=1+2\sum_{k=1}^{\infty}\frac{1+4k^2}{(1-4k^2)^2}.
$$
\section{Embedding the automorphic group within a larger group}\label{sec31}

Let $f\in E$. If $g\in E$ then ${\rm Aut}(f)\subseteq {\rm Aut}(g\circ f)$. For $\phi\in{\rm Aut}(f)\Rightarrow f\circ\phi=f$,
on the domain of definition of $\phi$. Hence $g\circ(f\circ\phi)=g\circ f$, on the domain of definition of $\phi$.
Thus $(g\circ f)\circ\phi=g\circ f$, on the domain of definition of $\phi$. This implies that $\phi\in{\rm Aut}(g\circ f)$.
Using this observation we deduce that if $\{g_n\}_n$ is a sequence of elements in $E$, then it induces the
following ascending sequence of discrete groups:
$$
{\rm Aut}(f)\subseteq{\rm Aut}(g_1\circ f)\subseteq{\rm Aut}(g_2\circ g_1\circ f)\subseteq\ldots .
$$
\begin{definition}\label{def7}
Let $(X_j,\tau_j)$ be a sequence of topological spaces such that for any two indices $i, j$ we have
$X_i\cap\tau_j\subseteq\tau_i$. This means that for any $V\in\tau_j$ we have $X_i\cap V\in\tau_i$.
We will define the direct limit topological space $\varinjlim (X_j,\tau_j)=(X,\tau)$ in this particular
setting by:
\begin{equation}\label{eq18}
\left\{\begin{array}{lll}X & = & \bigcup_j X_j \\ \tau & = & \{U\subseteq X\,|\,\forall\,j,\,U\cap X_j\in\tau_j\}\end{array}\right..
\end{equation}
\end{definition}

\begin{remark}\label{rem31}
The definition above is the standard definition of the final topology on the set $X=\bigcup_j X_j$ with respect
to the family of the inclusion mappings: $f_j\,:\,X_j\rightarrow X$, $f_j(x)=x$. Explicitly, a subset
$U\subseteq X$ is open in the final topology if and only if $\forall\,j,$ $f_j^{-1}(U)$ is open in $(X_j,\tau_j)$.
For in our case we have $\forall\,U\subseteq X$ and for any index $j$, $f_j^{-1}(U)=U\cap X_j$. Thus a
set $U\subseteq X$ is open in the final topology on $X$ if and only if $\forall\,j,\,U\cap X_j\in\tau_j$.
This is exactly the definition of the topology $\tau$ in equation (\ref{eq18}). We note that $(\bigcup_j X_j,\tau)$
is a discrete topological space if and only if $\forall\,j,\,(X_j,\tau_j)$ is a discrete topological space.
For $\forall\,x\in X$ we have by the definition of $\tau$: $\{x\}\in\tau$ if and only if $\forall\,j,\,\{x\}\cap X_j\in\tau_j$.
Clearly we have:
$$
\{x\}\cap X_j=\left\{\begin{array}{lll}\emptyset & {\rm if} & x\not\in X_j \\ \{x\} & {\rm if} & x\in X_j \end{array}\right..
$$
\end{remark}
\noindent
The topological groups ${\rm Aut}(g_n\circ\ldots\circ g_1\circ f)$ are discrete $\forall\,n\in\mathbb{Z}^+\cup\{0\}$.
Hence if we use the topology of definition \ref{def7}, the topological group 
$$
\bigcup_{n=0}^{\infty}{\rm Aut}(g_n\circ\ldots\circ g_1\circ f)
$$
is a discrete topological group. It will be useful to know when is it true that:
$$
\bigcup_{n=0}^{\infty}{\rm Aut}(g_n\circ\ldots\circ g_1\circ f)={\rm Aut}(F),
$$
for some $F\in E$? Let us assume that $\lim_{n\rightarrow\infty}g_n\circ\ldots\circ g_1\circ f=F$ uniformly on
compact subsets of $\mathbb{C}$. Will the answer to the question above be affirmative under this assumption?
Let $\phi\in\bigcup_{n=0}^{\infty}{\rm Aut}(g_n\circ\ldots\circ g_1\circ f)$. The for some $N\in\mathbb{Z}^+\cup\{0\}$
we have $\phi\in{\rm Aut}(g_N\circ\ldots\circ g_1\circ f)$. By the ascending property mentioned above (prior to
definition \ref{def7}) this implies that $\phi\in{\rm Aut}(g_k\circ\ldots\circ g_1\circ f)$, $\forall\,k\ge N$. This
is equivalent to $(g_k\circ\ldots\circ g_1\circ f)\circ\phi=(g_k\circ\ldots\circ g_1\circ f)$ on the domain of
definition of $\phi$, $\forall\,k\ge N$. Hence:
$$
\lim_{k\rightarrow\infty}((g_k\circ\ldots\circ g_1\circ f)\circ\phi)=\lim_{k\rightarrow\infty}(g_k\circ\ldots\circ g_1\circ f)=F,
$$
uniformly on compact subsets of $\mathbb{C}$. But:
$$
\lim_{k\rightarrow\infty}((g_k\circ\ldots\circ g_1\circ f)\circ\phi)=
\left(\lim_{k\rightarrow\infty}(g_k\circ\ldots\circ g_1\circ f)\right)\circ\phi=F\circ\phi,
$$
uniformly on compact subsets of $\mathbb{C}$. Hence $F\circ\phi=F$ on the domain of definition of $\phi$. Hence
$\phi\in{\rm Aut}(F)$. We proved the following:
\begin{proposition}\label{prop8}
If $\lim_{n\rightarrow\infty}(g_n\circ\ldots\circ g_1\circ f)=F$ uniformly on compact subsets of $\mathbb{C}$, then
$$
\bigcup_{n=0}^{\infty}{\rm Aut}(g_n\circ\ldots\circ g_1\circ f)\subseteq{\rm Aut}(F).
$$
\end{proposition}
\begin{remark}\label{rem32}
We note that by a variant of Hurwitz Theorem, the limit function $F$ in Proposition \ref{prop8} is either in $E$
or $F\equiv{\rm Const.}$ in which case it makes sense to define ${\rm Aut}(F)$ as the set of all functions
$\phi$. For in that case $F\circ\phi=({\rm Const.})\circ\phi={\rm Const.}=F$ for any $\phi$, on its domain
of definition. Thus in this case (when $F\equiv{\rm Const.}$) Proposition \ref{prop8} is clearly true.
\end{remark}
\noindent
When trying to find if equality sign can hold in Proposition \ref{prop8}, we clearly must assume that $F\in E$,
for there can be no equality if $F\equiv{\rm Const.}$. So we may assume that $\lim_{n\rightarrow\infty}
(g_n\circ\ldots\circ g_1\circ f)=F\in E$ uniformly on compact subsets of $\mathbb{C}$. Let $\psi\in{\rm Aut}(F)$.
Then $F\circ\psi=F$ on the domain of definition of $\psi$. Thus $(\lim_{n\rightarrow\infty}
(g_n\circ\ldots\circ g_1\circ f))\circ\psi=\lim_{n\rightarrow\infty}(g_n\circ\ldots\circ g_1\circ f)$. By the
continuity argument:
$$
\lim_{n\rightarrow\infty}((g_n\circ\ldots\circ g_1\circ f)\circ\psi)=\lim_{n\rightarrow\infty}(g_n\circ\ldots\circ g_1\circ f).
$$
If $\forall\,k\in\mathbb{Z}^+$, $\psi\not\in{\rm Aut}(g_k\circ\ldots\circ g_1\circ f)$, then $(g_k\circ\ldots\circ g_1\circ f)\circ\psi$
can not be extended to become an entire function. At least there seems to be no reason for such an extension
to exist, because $\psi$ may not be entire. However, $\psi$ is analytic on its domain of definition because
it is a branch of $F^{-1}(F(w))$. So we have a sequence of functions $\{(g_n\circ\ldots\circ g_1\circ f)\circ\psi\}n$,
analytic on (at least) a fixed open subset of $\mathbb{C}$ (the domain of the definition of $\psi$). The complement
of the domain of definition of $\psi$ is a closed subset of $\mathbb{C}$ that contains no continuum (by a theorem of
Julia). In that sense the open set is large. This sequence of analytic functions converges to a function $F$ which can
be extended to an entire function. Thus the limit process $\lim_{n\rightarrow\infty}((g_n\circ\ldots\circ g_1\circ f)\circ\psi)$
desingularizes the full set of singularities that originated in the automorphic function $\psi\in{\rm Aut}(F)$. By
results of Shimizu those singularities are branch points and can not include poles or algebraic poles. Thus at
each singular point the function $((g_n\circ\ldots\circ g_1\circ f)\circ\psi)$ is many valued. If there is a corresponding 
Hurwitz principle that holds true, then the limit function is either multivalued at such a singular point, or
a constant. We recall that the convergence is uniform on compact subsets of the domain of the definition
of $\psi$. Hence $\psi\not\in{\rm Aut}(g_k\circ\ldots\circ g_1\circ f)$ $\forall\,k\in\mathbb{Z}^+$, but
$\exists\,n_0$ such that for $n>n_0$, the function $(g_n\circ\ldots\circ g_1\circ f)\circ\psi$ has no singularity
and is analytic at a fixed singular point (a branch point) of $\psi$. So if the number of singular points
of $\psi$ is finite, then $\exists\,n_1$ such that for $n>n_1$, the function $(g_n\circ\ldots\circ g_1\circ f)\circ\psi$
is an entire function. So we have a sequence of entire functions $\{(g_n\circ\ldots\circ g_1\circ f)\circ\psi\}_{n>n_1}$
that converges uniformly on compact subsets of $\mathbb{C}$ to the entire function $F$. This gives a sequence
of non-zero entire functions $\{((g_n\circ\ldots\circ g_1\circ f)\circ\psi)-(g_n\circ\ldots\circ g_1\circ f)\}_{n>n_1}$
that converges uniformly on compact subsets of $\mathbb{C}$ to the zero function $F\circ\psi-F\equiv 0$.This, of
course, is a valid possibility by the Hurwitz Theorem. We recall the following: \\
\\
{\bf Lemma.}(\cite{er}){\it
$\forall\,f,g\in E,\,\,{\rm Aut}(f)\subseteq{\rm Aut}(g)\Leftrightarrow\exists\,h\in E\,\,{\rm such}\,{\rm that}\,\,
g=h\circ f$.
}
\\
\\
We conclude that if $\bigcup_{n=0}^{\infty}{\rm Aut}(g_n\circ\ldots\circ g_1\circ f)\subseteq{\rm Aut}(F)$ for
some $F\in E$, then $\forall\,n\in\mathbb{Z}^+$ $\exists\,F_n\in E$ such that $F=f_n\circ(g_n\circ\ldots\circ g_1\circ f)$.
By Proposition \ref{prop8}, this will be the case when $\lim_{n\rightarrow\infty}(g_n\circ\ldots\circ g_1\circ f)=F$
uniformly on compact subsets of $\mathbb{C}$. Hence in this case we have $\lim_{n\rightarrow\infty} F_n={\rm id.}$
uniformly on compact subsets of $\mathbb{C}$.

\begin{remark}\label{rem33}
If $F_n\in E$ and $\lim_{n\rightarrow\infty}F_n=z$ uniformly on compact subsets of $\mathbb{C}$, then
${\rm Aut}(F_n)=\{z\}\cup\{\phi_{0k}^{(n)}\}_k$ where $\forall\,k,\,\lim_{n\rightarrow\infty}\phi_{0k}^{(n)}=\infty$.
In that sense ${\rm Aut}(F_n)\rightarrow_{n\rightarrow\infty}\{z\}$.
\end{remark}
\noindent
We recall the following: \\
\\
{\bf Theorem.} (\cite{tuen}) {\it
There exists a sequence of positive real numbers $\{c_n\}_{n=1}^{\infty}$ such that the sequence of functions
$F_n(z)=(c_n e^z+z)\circ\ldots\circ(c_1 e^z+z)$ converges uniformly on compact subsets of $\mathbb{C}$ to an
entire function $F(z)$. Furthermore, for each $n\in\mathbb{Z}^+$, $F(z)=H_n\circ(c_n e^z+z)\circ\ldots\circ(c_1 e^z+z)$
for some entire function $H_n$. Hence, there is no uniform bound on the number of prime factors $c_n e^z+z$
in different decompositions of $F$ through transcendental entire functions.
}
\begin{remark}\label{rem34}
A similar result holds for factorization that go in the opposite direction, i.e. $(c_1 e^z+z)\circ\ldots\circ(c_n e^z+z)$.
\end{remark}
\noindent
As for the Riemann surface of the inverse functions that are the limits of factorizations of non-bounded
number of factors, $F^{-1}=(c_1 e^z+z)^{-1}\circ(c_2 e^z+z)^{-1}\circ\ldots$. It contains the embedded copies of
the Riemann surfaces of the factors nested one on the top of the other. One can outline the geometric 
construction of the Riemann surface of that is induced by $g_n\circ\ldots\circ g_1\circ f$. if
$\Gamma=\bigcup_{n=0}^{\infty}{\rm Aut}(g_n\circ\ldots\circ g_1\circ f)={\rm Aut}(F)$ for some $F\in E$,
then the direct limit of those Riemann surfaces will be equal to the Riemann surface of $F^{-1}$. If, however,
the discrete group $Gamma$ does not equal to an ${\rm Aut}(F)$ for some $F\in E$, then this direct limit of
Riemann surfaces will not be a Riemann surface. This structure generalizes the Riemann surfaces.

\section{Relations between the construction of the direct system of the automorphic groups,
and Weierstrass products}\label{sec32}

We recall that for $f\in E$, the elements of the automorphic group ${\rm Aut}(f)$, are the functions of
$f^{-1}\circ f$. Let $g_1\in E$ then ${\rm Aut}(f)\subseteq{\rm Aut}(g_1\circ f)$. In fact the elements
of ${\rm Aut}(g_1\circ f)$ are the functions of $(g_1\circ f)^{-1}\circ(g_1\circ f)=
f^{-1}\circ(g_1^{-1}\circ g_1)\circ f=f^{-1}\circ{\rm Aut}(g_1)\circ f$. We note that by taking the
identity element ${\rm id}$ in ${\rm Aut}(g_1)$ we get $f^{-1}\circ{\rm id}\circ f={\rm Aut}(f)$ which
explains the relation ${\rm Aut}(f)\subseteq{\rm Aut}(g_1\circ f)$. Both $f(w)-f(z)$ and $(g_1\circ f)(w)-
(g_1\circ f)(z)$ of the variable $w$, with the parameter $z\in\mathbb{C}-f^{-1}(f(0))$ in the first and
$z\in\mathbb{C}-(g_1\circ f)^{-1}((g_1\circ f)(0))$ in the second, have Weierstrass products representation
that are based on the product:
$$
\prod_{\phi_{0n}\in{\rm Aut}(f)}\left(1-\frac{w}{\phi_{0n}(z)}\right),
$$
for $f(w)-f(z)$ and on the product:
$$
\prod_{\psi_{0n}\in{\rm Aut}(g_1\circ f)}\left(1-\frac{w}{\psi_{0n}(z)}\right),
$$
for $(g_1\circ f)(w)-(g_1\circ f)(z)$. Since ${\rm Aut}(f)\subseteq{\rm Aut}(g_1\circ f)$, any factor of the first product
is also a factor of the second product. In that sense the first product divides the second one. We will
denote that by standard notation:
$$
\prod_{\phi_{0n}\in{\rm Aut}(f)}\left(1-\frac{w}{\phi_{0n}(z)}\right)\left|
\prod_{\psi_{0n}\in{\rm Aut}(g_1\circ f)}\left(1-\frac{w}{\psi_{0n}(z)}\right)\right..
$$
If we actually divide the full detailed second product by the full detailed first product, we obtain
a Weierstrass product type representation for the meromorphic function of $w$,
$$
\frac{(g_1\circ f)(w)-(g_1\circ f)(z)}{f(w)-f(z)}.
$$
We will denote that by:
$$
\frac{\prod_{\phi_{0n}\in{\rm Aut}(f)}\left(1-\frac{w}{\phi_{0n}(z)}\right)}
{\prod_{\psi_{0n}\in{\rm Aut}(g_1\circ f)}\left(1-\frac{w}{\psi_{0n}(z)}\right)}\thicksim
\frac{(g_1\circ f)(w)-(g_1\circ f)(z)}{f(w)-f(z)}.
$$
We note that the meromorphic function of $w$, $\frac{(g_1\circ f)(w)-(g_1\circ f)(z)}{f(w)-f(z)}$, is in fact
an entire function of $w$, because it has only removable singularities and not poles in the finite plane. Symbolically
we have the following assignment:
$$
f(w)-f(z)\rightarrow\prod_{{\rm Aut}(f)},\,\,\,\,\,\,\,g_1(f(w))-g_1(f(z))\rightarrow\prod_{{\rm Aut}(g_1\circ f)},
$$
$$
\left\{\frac{g_1(f(w))-g_1(f(z))}{f(w)-f(z)}\right\}\rightarrow\prod_{{\rm Aut}(g_1\circ f)-{\rm Aut}(f)}.
$$
Similarly we can go on:
$$
g_2(g_1(f(w)))-g_2(g_1(f(z)))\rightarrow\prod_{{\rm Aut}(g_2\circ g_1\circ f)},
$$
$$
\left\{\frac{g_2(g_1(f(w)))-g_2(g_1(f(z)))}{g_1(f(w))-g_1(f(z))}\right\}\prod_{{\rm Aut}(g_2\circ g_1\circ f)-{\rm Aut}(g_1\circ f)},
$$
$$
\left\{\frac{g_2(g_1(f(w)))-g_2(g_1(f(z)))}{f(w)-f(z)}\right\}\prod_{{\rm Aut}(g_2\circ g_1\circ f)-{\rm Aut}(f)}.
$$
We note the consistency:
$$
\left\{\frac{g_2(g_1(f(w)))-g_2(g_1(f(z)))}{f(w)-f(z)}\right\}=\left\{\frac{g_2(g_1(f(w)))-g_2(g_1(f(z)))}{g_1(f(w))-g_1(f(z))}\right\}
\times \left\{\frac{g_1(f(w))-g_1(f(z))}{f(w)-f(z)}\right\}\rightarrow
$$
$$
\prod_{{\rm Aut}(g_2\circ g_1\circ f)-{\rm Aut}(f)}=\prod_{{\rm Aut}(g_2\circ g_1\circ f)-{\rm Aut}(g_1\circ f)}\times
\prod_{{\rm Aut}(g_1\circ f)-{\rm Aut}(f)}.
$$
If we denote union by plus: $+$, then it corresponds to multiplication. This is in agreement with the fact that
minus: $-$ corresponded to division. In this notation we have:
$$
{\rm Aut}(g_2\circ g_1\circ f)-{\rm Aut}(f)=({\rm Aut}(g_2\circ g_1\circ f)-{\rm Aut}(g_1\circ f))+
({\rm Aut}(g_1\circ f)-{\rm Aut}(f).
$$
It is clear that in general we have:
\begin{proposition}\label{prop9}
If $g_n, f\in E,\,\,\forall\,n\in\mathbb{Z}^+$, then:
$$
(g_n\circ\ldots\circ g_1\circ f)(w)-(g_n\circ\ldots\circ g_1\circ f)(z)\rightarrow {\rm Aut}(g_n\circ\ldots\circ g_1\circ f),
$$
and $\forall\,n>m\ge 1$ in $\mathbb{Z}^+$:
$$
\left\{\frac{(g_n\circ\ldots\circ g_1\circ f)(w)-(g_n\circ\ldots\circ g_1\circ f)(z)}
{(g_m\circ\ldots\circ g_1\circ f)(w)-(g_m\circ\ldots\circ g_1\circ f)(z)}\right\}\rightarrow
{\rm Aut}(g_n\circ\ldots\circ g_1\circ f)-{\rm Aut}(g_m\circ\ldots\circ g_1\circ f).
$$
\end{proposition}
\noindent
another suggestive assignment which is natural, is the exponential and the logarithmic notations:
$$
\left\{\begin{array}{ll} \log\left((g_n\circ\ldots\circ g_1\circ f)(w)-(g_n\circ\ldots\circ g_1\circ f)(z)\right)\rightarrow
{\rm Aut}(g_n\circ\ldots\circ g_1\circ f) \\
\exp\left({\rm Aut}(g_n\circ\ldots\circ g_1\circ f)\right)\rightarrow (g_n\circ\ldots\circ g_1\circ f)(w)-
(g_n\circ\ldots\circ g_1\circ f)(z)\end{array}\right.
$$
Next we note that the discrete group $\Gamma=\bigcup {\rm Aut}(g_n\circ\ldots\circ g_1\circ f)$ can also be denoted by
$\Gamma=\sum_{n=0}^{\infty}{\rm Aut}(g_n\circ\ldots\circ g_1\circ f)$ and formally it should be assigned to the
Weierstrass type product:
$$
\prod_{\theta_{0n}\in\Gamma}\left(1-\frac{w}{\theta_{0n}(z)}\right).
$$
On the other hand $\forall\,n\in\mathbb{Z}^+$ we have:
$$
{\rm Aut}(g_n\circ\ldots\circ g_1\circ f)=\sum_{k=1}^{n-1}\left({\rm Aut}(g_{k+1}\circ\ldots\circ g_1\circ f)-
{\rm Aut}(g_k\circ\ldots\circ g_1\circ f)\right)+
$$
$$
+({\rm Aut}(g_1\circ f)-{\rm Aut}(f))+{\rm Aut}(f).
$$
That corresponds to:
$$
(g_n\circ\ldots\circ g_1\circ f)(w)-(g_n\circ\ldots\circ g_1\circ f)(z)=
$$
$$
=\prod_{k=1}^{n-1}\left\{\frac{(g_{k+1}\circ\ldots\circ g_1\circ f)(w)-(g_{k+1}\circ\ldots\circ g_1\circ f)(z)}
{(g_k\circ\ldots\circ g_1\circ f)(w)-(g_k\circ\ldots\circ g_1\circ f)(z)}\right\}\times
$$
$$
\times\left\{\frac{(g_1\circ f)(w)-(g_1\circ f)(z)}{f(w)-f(z)}\right\}\times(f(w)-f(z)).
$$
If indeed $\Gamma={\rm Aut}(F)$ for some $F\in E$, then:
$$
F(w)-F(z)\rightarrow\prod_{{\rm Aut}(F)}=\prod_{\Gamma}\rightarrow
$$
$$
\rightarrow\prod_{k=1}^{n-1}\left\{\frac{(g_{k+1}\circ\ldots\circ g_1\circ f)(w)-(g_{k+1}\circ\ldots\circ g_1\circ f)(z)}
{(g_k\circ\ldots\circ g_1\circ f)(w)-(g_k\circ\ldots\circ g_1\circ f)(z)}\right\}\times
$$
$$
\times\left\{\frac{(g_1\circ f)(w)-(g_1\circ f)(z)}{f(w)-f(z)}\right\}\times(f(w)-f(z))=
$$
$$
=\lim_{n\rightarrow\infty}((g_n\circ\ldots\circ g_1\circ f)(w)-(g_n\circ\ldots\circ g_1\circ f)(z)).
$$
This implies:
\begin{theorem}\label{thm18}
Let $g_n, f\in E\,\,\,\forall\,n\in\mathbb{Z}^+$. Then there exists an $F\in E$ such that:
$$
\bigcup_{n=0}^{\infty}{\rm Aut}(g_n\circ\ldots\circ g_1\circ f)={\rm Aut}(F),
$$
if and only if $G=\lim_{n\rightarrow\infty}(g_n\circ\ldots\circ g_1\circ f)$ exists uniformly
on compact subsets of $\mathbb{C}$ and is not a constant.In the case the limit $G\not\equiv{\rm Const.}$
then $G\in E$ and we can take $F(w)=a\cdot G(w)+b$ $\forall\,a\in\mathbb{C}^{\times}$ and $\forall\,b\in\mathbb{C}$.
\end{theorem}
\noindent
{\bf Example:} If $F(z)=\lim_{n\rightarrow\infty}(c_ne^z+z)\circ\ldots\circ(c_1e^z+z)$ is a Tuen Wai NG entire
function then we have the identity:
$$
{\rm Aut}(F)=\bigcup_{n=1}^{\infty}{\rm Aut}((c_ne^z+z)\circ\ldots\circ(c_1e^z+z)),
$$
where the entire functions $c_ke^z+z$ ($c_k>0$) are primes. We also have the functional identity:
$$
F(w)-F(z)=
$$
$$
=\prod_{n=1}^{\infty}\left\{\frac{((c_{n+1}e^w+w)\circ\ldots\circ(c_1e^w+w))
-((c_{n+1}e^z+z)\circ\ldots\circ(c_1e^z+z))}{((c_ne^w+w)\circ\ldots\circ(c_1e^w+w))-
((c_ne^z+z)\circ\ldots\circ(c_1e^z+z))}\right\}\times
$$
$$
\times((c_1e^w+w)-(c_1e^z+z)).
$$
\begin{theorem}\label{thm19}
Let $g_n,f\in E\,\,\forall\,n\in\mathbb{Z}^+$. If the limit $G=\lim_{n\rightarrow\infty}(g_n\circ\ldots\circ g_1\circ f)$
exists uniformly on compact subsets of $\mathbb{C}$ and is not a constant, then there exists a path metric
$\rho:\,\mathbb{C}\times\mathbb{C}\rightarrow\mathbb{R}_{\ge 0}$ which is invariant for ${\rm Aut}(g_n\circ\ldots\circ g_1
\circ f)$, $\forall\,n=0,1,2,\ldots$. In other words $\forall\,n\ge 0$, $\forall\,\phi\in{\rm Aut}(g_n\circ\ldots\circ g_1
\circ f)$, $\phi$ is a $\rho$-isometry on a domain of definition of a leaf of $\phi$, i.e. $\forall\,z,w$ we have:
$\rho(\phi(z),\phi(w))=\rho(z,w)$.
\end{theorem}
\noindent
{\bf Proof.} \\
We have by the assumption on $\lim_{n\rightarrow\infty}(g_n\circ\ldots\circ g_1\circ f)$ the containment:
$$
\bigcup_{n=0}^{\infty}{\rm Aut}(g_n\circ\ldots\circ g_1\circ f)\subseteq{\rm Aut}(G).
$$
By the way, we do not need Theorem \ref{thm18} for this. Now take the path metric on $\mathbb{C}$ induced
by $G$, $\rho=\rho_G:\mathbb{C}\times\mathbb{C}\rightarrow\mathbb{R}_{\ge 0}$. We know that any $G$-automorphic
function $\phi\in{\rm Aut}(G)$ is a $\rho_G$-isometry in the sense of the theorem. See Theorem \ref{thm14}. $\qed $ \\
\\
\begin{theorem}\label{thm20}
Let $h_n,\,g_n,\,f\in E\,\,\forall\,n\in\mathbb{Z}^+$. If $\lim_{n\rightarrow\infty}(g_n\circ\ldots\circ g_1\circ f)=
\lim_{n\rightarrow\infty}(h_n\circ\ldots\circ h_1\circ f)$ exist uniformly on compact subsets of $\mathbb{C}$
and the limit function is not a constant, then:
$$
\bigcup_{n=0}^{\infty}{\rm Aut}(g_n\circ\ldots\circ g_1\circ f)=\bigcup_{n=0}^{\infty}{\rm Aut}(h_n\circ\ldots\circ h_1\circ f).
$$
\end{theorem}
\noindent
{\bf Proof.} \\
Let $G$ be the limit function of the two sequences $\{g_n\circ\ldots\circ g_1\circ f\}_n$ and
$\{h_n\circ\ldots\circ h_1\circ f\}_n$ of functions in $E$. Then $G\in E$ and by Theorem \ref{thm18}
the unions of the automorphic groups, both, are equal to ${\rm Aut}(G)$. $\qed $ \\
\\
So far our construction gives under the appropriate conditions the identity:
$$
{\rm Aut}(F)=\bigcup_{n=0}^{\infty}{\rm Aut}(g_n\circ\ldots\circ g_1\circ f),
$$
where the functions $F,\,g_n,\,f\in E,\,\,\forall\,n\in\mathbb{Z}_{\ge 0}$. Hence all the automorphic groups
that are involved are discrete groups and are countable. Trivially any discrete group is a locally compact
Hausdorff group. The compact subsets in a discrete group are the finite subsets, and the Haar measure up
to a multiplication by a positive constant is the counting measure.

\begin{definition}\label{def8}
If $H$ is a subgroup of the topological group $G$, then it induces two relations on $G$: \\
(a) The $H$-right relation: $\gamma_1\sim_{H-{\rm right}}\gamma_2\Leftrightarrow\exists\,\delta\in H\,{\rm such}\,\,
{\rm that}\, \gamma_1=\gamma_2\cdot\delta$. \\
(b) The $H$-left relation: $\gamma_1\sim_{H-{\rm left}}\gamma_2\Leftrightarrow\exists\,\delta\in H\,{\rm such}\,\,
{\rm that}\, \gamma_1=\delta\cdot\gamma_2$.
\end{definition}

\begin{proposition}\label{prop10}
Let $H$ be a subgroup of the topological group $G$, then both $H$-right and $H$-left are equivalence relations on $G$.
\end{proposition}
\noindent
{\bf Proof.} \\
This is straight forward from Definition \ref{def8}. $\qed $ \\

\begin{definition}\label{def9}
We will denote the equivalence classes of $G$ with respect to the equivalence relation $H$-right by $(G/H)_{\rm right}$.
Similarly $(G/H)_{\rm left}$ will denote the family of equivalence classes with respect to $H$-left.
\end{definition}

\begin{theorem}\label{thm21}
Let $h,\,f\in E$. Then: \\
(a) $\forall\,\gamma_1,\,\gamma_2\in{\rm Aut}(h\circ f)$ we have $\gamma_1\sim_{{\rm Aut}(f)-{\rm left}}\gamma_2
\Leftrightarrow f(\gamma_1)=f(\gamma_2)$. \\
(b) $\forall\,\gamma_1,\,\gamma_2\in{\rm Aut}(h\circ f)$ we have $\gamma_1\sim_{{\rm Aut}(f)-{\rm right}}\gamma_2
\Leftrightarrow f(\gamma_1^{-1})=f(\gamma_2^{-1})$.
\end{theorem}
\noindent
{\bf Proof.} \\
(a) $\gamma_1\sim_{{\rm Aut}(f)-{\rm left}}\gamma_2\Leftrightarrow\gamma_2=\psi\circ\gamma_1$ for some $\psi\in{\rm Aut}(f)$
$\Leftrightarrow f(\gamma_2)=f(\psi\circ\gamma_1)=(f\circ\psi)\circ\gamma_1=f(\gamma_1)$. \\
(b) $\gamma_1\sim_{{\rm Aut}(f)-{\rm right}}\gamma_2\Leftrightarrow\gamma_2=\gamma_1\circ\psi$ for some $\psi\in{\rm Aut}(f)$
$\Leftrightarrow\gamma_2^{-1}=\psi^{-1}\circ\gamma_1^{-1}\Leftrightarrow\gamma_1^{-1}\sim_{{\rm Aut}(f)-{\rm left}}
\gamma_2^{-1}\Leftrightarrow f(\gamma_1^{-1})=f(\gamma_2^{-2})$ where in the last step we made a use in (a). $\qed $ \\

\begin{theorem}\label{thm22}
The cardinality of the equivalence classes in $({\rm Aut}(h\circ f)/{\rm Aut}(f))_{\rm left}$ and in
$({\rm Aut}(h\circ f)/{\rm Aut}(f))_{\rm right}$ is equal to the cardinality of ${\rm Aut}(f)$, and hence
are at most $\aleph_0$.
\end{theorem}
\noindent
{\bf Proof.} \\
By Definition \ref{def8} it follows that for any $[\gamma]\in({\rm Aut}(h\circ f)/{\rm Aut}(f))_{\rm left}$
we have: $[\gamma]=\{\psi\circ\gamma\,|\,\psi\in{\rm Aut}(f)\}$. Since $\psi\circ\gamma=\psi_1\circ\gamma
\Leftrightarrow\psi=\psi_1$ it follows that the mapping:
$$
[\gamma]\rightarrow{\rm Aut}(f),\,\,\,\,\,\psi\circ\gamma\rightarrow\psi,
$$
is a bijection. Hence $|[\gamma]|=|{\rm Aut}(f)|$. A similar argument works for the ${\rm Aut}(f)$-right equivalence
relation. $\qed $ \\

\begin{remark}\label{rem35}
If $f(z)\in E$ is a transcendental entire function then the equivalence classes in both left and right ${\rm Aut}(f)$
equivalence relations have cardinality $\aleph_0$.
\end{remark}

\begin{theorem}\label{thm23}
If $h,\,f\in E$ then both topological spaces $({\rm Aut}(h\circ f)/{\rm Aut}(f))_{\rm left}$, 
$({\rm Aut}(h\circ f)/{\rm Aut}(f))_{\rm right}$ are discrete and Hausdorff.
\end{theorem}
\noindent
{\bf Proof.} \\
This follows by the following well known facts: \\
If $H$ is a subgroup of $G$ then $(G/H)$ is discrete if and only if $H$ is open in $G$. \\
$(G/H)$ is Hausdorff if and only if $H$ is closed in $G$. \\
In our case ${\rm Aut}(h\circ f)$ is a discrete group and hence ${\rm Aut}(f)$ is both open and closed
in ${\rm Aut}(h\circ f)$. $\qed $

\section{Continuity properties of the automorphic groups}\label{sec33}

In this section we study the following: Let $f\in E$ and let $f_n\in E\,\,\forall\,n\in\mathbb{Z}^+$. Suppose
that $\lim_{n\rightarrow\infty} f_n=f$ uniformly on compact subsets of $\mathbb{C}$. Is it true that the automorphic
groups ${\rm Aut}(f_n)$ become closer to the automorphic group ${\rm Aut}(f)$? If the answer is affirmative, in
what sense? \\
Clearly an attractive situation is the one in which $f_n(z)\in\mathbb{C}[z]$, i.e. the approximating sequence
is a sequence of polynomials. For example the partial sums of the power series expansion of $f$:
$$
f_n(z)=\sum_{k=0}^{n}\frac{f^{(k)}(0)}{k!},\,\,\,\,\,n\in\mathbb{Z}^+.
$$
Since the functions that constitute ${\rm Aut}(f)$ are $f^{-1}\circ f$, it makes sense to find if in some sense
the many valued functions $f_n^{-1}$ approach $f^{-1}$. We recall once more the following well-known: \\
\\
{\bf Theorem (The generalized argument principle).} {\it
Let $F$ be a meromorphic function in the simply connected domain $D$, $a_j$ the zeros of $F$, $b_k$ the poles of
$F$ in $D$ and $\gamma$ a closed curve in $D$ avoiding the $a_j,\,b_k$. Then $\forall\,G\in C^{\omega}(D)$ we have:
$$
\sum_j G(a_j)\cdot n(\gamma,a_j)-\sum_k G(b_k)\cdot n(\gamma,b_k)=\frac{1}{2\pi i}\oint_{\gamma}G(z)\cdot\frac{F'(z)}{F(z)}dz.
$$
Here we have, for any $a\not\in\gamma$:
$$
n(\gamma,a)=\frac{1}{2\pi i}\oint_{\gamma}\frac{dz}{z-a},
$$
is the index or the winding number of the closed curve $\gamma$ with respect to the point $a\not\in\gamma$.
}

\begin{theorem}\label{thm24}
Let $f\in E$ and let $f_n\in E\,\,\forall\,n\in\mathbb{Z}^+$. Suppose that $\lim_{n\rightarrow\infty} f_n=f$
uniformly on compact subsets of $\mathbb{C}$. Let us denote the automorphic group’s elements by 
${\rm Aut}(f)=\{\phi_{0k}\}_k$, ${\rm Aut}(f_n)=\{\phi_{0k}^{[n]}\}_k$. Then for any $R>0$ and for any
$\epsilon>0$, $\exists\,N=N(R,\epsilon)$ such that: \\
(1) $\forall\,n>N(R,\epsilon)$, the number of $\phi_{0k}^{[n]}(z)$ for a fixed $z$, such that $|\phi_{0k}^{[n]}(z)|<R$
equals (counting with multiplicity) the number of those $\phi_{0k}(z)$ for which $|\phi_{0k}(z)|<R$. \\
(2) The indexing of the $\phi_{0k}$ and of the $\phi_{0k}^{[n]}$, can be arranged, so that $\forall\,n>N(R,\epsilon)$
and $\forall\,z$ such that $|\phi_{0k}(z)|<R$, we have $|\phi_{0k}(z)-\phi_{0k}^{[n]}(z)|<\epsilon$.
\end{theorem}
\noindent
{\bf Proof.} \\
We recall some elementary facts from the algebra of polynomials in one variable: Let $\{\alpha_1,\alpha_2,\ldots,\alpha_m\}$
be a finite sequence of complex numbers and let $\{\{\alpha_1^{[n]},\alpha_2^{[n]},\ldots,\alpha_{m}^{[n]}\}\}_n$ be an
infinite sequence of finite sequences over $\mathbb{C}$ of the same length $m$ as the first. We denote the moments by
$m_k(\alpha_1,\ldots,\alpha_m)=\alpha_1^k+\ldots+\alpha_m^k$, and similarly $m_k(\alpha_1^{[n]},\ldots,\alpha_m^{[n]})=
(\alpha_1^{[n]})^k+\ldots+(\alpha_m^{[n]})^k,\,\,k\in\mathbb{Z}^+$. If $\forall\,k$ we have:
$$
\lim_{n\rightarrow\infty}m_k(\alpha_1^{[n]},\ldots,\alpha_m^{[n]})=m_k(\alpha_1,\ldots,\alpha_m),
$$
then the indexing of the $\alpha_i$ and of the $\alpha_i^{[n]}$, $1\le i\le m$ can be arranged, so that:
$$
\lim_{n\rightarrow\infty}\alpha_i^{[n]}=\alpha_i,\,\,\,\,\,1\le i\le m.
$$
The reason for this is the following: Let us denote the symmetric functions of the sequence
$\{\alpha_1,\ldots,\alpha_m\}$ by 
$$
S_k(\alpha_1,\ldots,\alpha_m)=\sum_{1\le i_1<\ldots<i_k\le m}\alpha_{i_1}\alpha_{i_2}\ldots\alpha_{i_k},
\,\,\,1\le k\le m.
$$
Similarly $S_k(\alpha_1^{[n]},\ldots,\alpha_m^{[n]})$, $1\le k\le m$, denote the symmetric functions of the
other sequences. Then each moment $m_k$ can be written as a polynomial over $\mathbb{Q}$, with fixed coefficients
for a given $k$, of the symmetric functions and vice versa. These are known as Newton's identities. They start
as follows: $m_1=S_1$, $m_2=S_2-S_1^2,\ldots$ and $S_1=m_1$, $S_2=m_1^2-m_2\ldots $. By the assumption
$\lim_{n\rightarrow\infty}m_k(\alpha_1^{[n]},\ldots,\alpha_m^{[n]})=m_k(\alpha_1,\ldots,\alpha_m)$, it follows
that $\lim_{n\rightarrow\infty}S_k(\alpha_1^{[n]},\ldots,\alpha_{m}^{[n]})=S_k(\alpha_1,\ldots,\alpha_m)$ for
$1\le k\le m$. We note that for the monic polynomials that have as their zero sets the negatives of these
sequences we have:
$$
P(w)=(w+\alpha_1)\ldots(w+\alpha_m)=w^m+S_1(\alpha)w^{m-1}+\ldots+S_k(\alpha)w^{m-k}+\ldots+S_m(\alpha),
$$
$$
P^{[n]}(w)=(w+\alpha_1^{[n]})\ldots(w+\alpha_m^{[n]})=
w^m+S_1(\alpha^{[n]})w^{m-1}+\ldots+S_k(\alpha^{[n]})w^{m-k}+\ldots+S_m(\alpha^{[n]}).
$$
Hence $\lim_{n\rightarrow\infty}P^{[n]}(w)=P(w)$ uniformly on compact subsets of $\mathbb{C}$. Hence after
an appropriate indexing the zeros of the $P^{[n]}$ (multiplicity included) approach as their limits when
$n\rightarrow\infty$ the zeros of $P$ and we proved the claimed fact. \\
Coming back to our entire functions, we recall that by using the Weierstrass representation as a canonical product
of the function $f(w)-f(z)$ which is entire in $w$ where $z$ is a fixed parameter, and differentiating
$\log(f(w)-f(z))$ with respect to $w$, we obtain:
$$
\frac{f'(w)}{f(w)-f(z)}=\frac{\partial g}{\partial w}(w,z)+\sum_k\left(\frac{w}{\phi_{0k}(z)}\right)^{\lambda_k}
\left(\frac{1}{w-\phi_{0k}}\right).
$$
Hence we have (using the generalized argument principle):
$$
\frac{1}{2\pi i}\oint_{|w|=R}\frac{f'(w)dw}{f(w)-f(z)}=
$$
\begin{equation}\label{eq19}
\end{equation}
$$
=\frac{1}{2\pi i}\oint_{|w|=R}\frac{\partial g}{\partial w}(w,z)dw+
\sum_k\frac{1}{2\pi i}\oint_{|w|=R}\left(\frac{w}{\phi_{0k}(z)}\right)^{\lambda_k}\left(\frac{dw}{w-\phi_{0k}}\right).
$$
We have:
$$
\frac{\partial g}{\partial w}(w,z)\in C^{\omega}(\mathbb{C},w)\Rightarrow\frac{1}{2\pi i}\oint_{|w|=R}
\frac{\partial g}{\partial w}(w,z)dw=0,
$$
$$
\sum_k\frac{1}{2\pi i}\oint_{|w|=R}\left(\frac{w}{\phi_{0k}(z)}\right)^{\lambda_k}\left(\frac{dw}{w-\phi_{0k}}\right)=
\sum_{|\phi_{0k}(z)|<R}\frac{1}{(\phi_{0k}(z))^{\lambda_k}}\cdot(\phi_{0k}(z))^{\lambda_k}=
$$
$$
=|\{k|\,|\phi_{0k}(z)|<R\}|.
$$
Next we have:
$$
\frac{1}{2\pi i}\oint_{|w|=R}w^l\cdot\frac{f'(w)dw}{f(w)-f(z)}=
$$
\begin{equation}\label{eq20}
\end{equation}
$$
=\frac{1}{2\pi i}\oint_{|w|=R}w^l\frac{\partial g}{\partial w}(w,z)dw+
\sum_k\frac{1}{2\pi i}\oint_{|w|=R}w^l\left(\frac{w}{\phi_{0k}(z)}\right)^{\lambda_k}\left(\frac{dw}{w-\phi_{0k}}\right).
$$
Once more, by the Cauchy Theorem:
$$
\frac{1}{2\pi i}\oint_{|w|=R}w^l\frac{\partial g}{\partial w}(w,z)dw=0,
$$
and by the generalized argument principle:
$$
\sum_k\frac{1}{2\pi i}\oint_{|w|=R}w^l\left(\frac{w}{\phi_{0k}(z)}\right)^{\lambda_k}\left(\frac{dw}{w-\phi_{0k}}\right)=
\sum_{|\phi_{0k}(z)|<R}\phi_{0k}^l(z).
$$
Hence we proved the following integral identity for the moments of the automorphic functions:
$$
m_l(\phi_{0k}(z)||\phi_{0k}(z)|<R)=\frac{1}{2\pi i}\oint_{|w|=R}w^l\cdot\frac{f'(w)dw}{f(w)-f(z)}.
$$
Similarly we have for any $n\in\mathbb{Z}^+$:
$$
m_l(\phi_{0k}^{[n]}(z)||\phi_{0k}^{[n]}(z)|<R)=\frac{1}{2\pi i}\oint_{|w|=R}w^l\cdot\frac{f_n'(w)dw}{f_n(w)-f_n(z)}.
$$
By the assumption: $\lim_{n\rightarrow\infty}f_n=f$ uniformly on compact subsets of $\mathbb{C}$ and by the
Cauchy estimates: $\lim_{n\rightarrow\infty}f'_n=f'$ uniformly on compact subsets of $\mathbb{C}$. This implies
that:
$$
\lim_{n\rightarrow\infty}\frac{1}{2\pi i}\oint_{|w|=R}w^l\cdot\frac{f_n'(w)dw}{f_n(w)-f_n(z)}=
\frac{1}{2\pi i}\oint_{|w|=R}w^l\cdot\frac{f'(w)dw}{f(w)-f(z)}.
$$
We proved:
$$
\lim_{n\rightarrow\infty}m_l(\phi_{0k}^{[n]}(z)||\phi_{0k}^{[n]}(z)|<R)=m_l(\phi_{0k}(z)||\phi_{0k}(z)|<R).
$$
Now the assertions of our theorem follow by the first part of our proof. $\qed $ \\
\\One way to interpret Theorem \ref{thm24} is that the automorphic functions of the approximating functions
$f_n$ to the entire function $f\in E$, converge themselves to the automorphic functions of $f$. This convergence
is very ordered and not chaotic. By that we mean that from a certain index $n_0$ and on it is unambiguous for
certain of the automorphic functions of $f$ which of the automorphic functions of $f_n$ (for $n$ large enough)
correspond to them. This happens because when we fix the value $z$ of the complex parameter in $f(w)-f(z)$
and in $f_n(w)-f_n(z)$ and consider the disk $B(0,R)$ and only those automorphic functions $\phi_{0k}$ of $f$,
$\phi_{0k}\in{\rm Aut}(f)$ whose $z$-image lies inside that disk, i.e. $|\phi_{0k}(z)|<R$ and take in Theorem
\ref{thm24} the positive $\epsilon$, small enough, then for values of the index $n>n_0$ it is clear which which
of the automorphic functions of $f_n$ is the one that corresponds to a particular $\phi_{0k}$. We changed the 
indices so that $|\phi_{0k}(z)-\phi_{0k}^{[n]}(z)|<\epsilon$. In other words the values $\phi_{0k}^{[n]}(z)$
for $n>n_0$ (in Theorem \ref{thm24} we denoted $n_0=N(R,\epsilon)$), are trapped inside a small circle of a
radius $\epsilon$ centered at $\phi_{0k}(z)$. The un-ambiguity follows because for a small enough $\epsilon>0$,
the disks $B(\phi_{0k}(z),\epsilon)$ for $|\phi_{0k}(z)|<R$ have disjoint closures. We can achieve this by
choosing $\epsilon<\frac{1}{2}\min\{|\phi_{0,k_1}(z)-\phi_{0,k_2}(z)|||\phi_{0,k_1}(z)|,|\phi_{0,k_2}(z)|<R,
\phi_{0,k_1}(z)\ne\phi_{0,k_2}(z)\}$. The minimum exists because the set $\{\phi_{0,k}(z)||\phi_{0,}(z)|<R\}$
is a finite set. The number $\epsilon$ should also be smaller than $\min\{R-|\phi_{0,k}(z)|||\phi_{0,k}(z)|<R\}$.
Every value $z$ of complex parameter determines such a configuration as the one described above. Thus those 
configurations (that geometrically look like an open disk of radius $R$ punctured by finitely many small disks
of radius $\epsilon$ that have disjoint closures and that stay away from $\partial B(0,R)$) are determined
by three quantities:
$$
(z,R,\epsilon)\in\left(\mathbb{C}-f^{-1}(f(0))\cup\bigcup_n f_n^{-1}(f_n(0))\right)\times\mathbb{R}^+\times
(0,\delta(z,R)),
$$
where we have the formula:
$$
\delta(z,R)=\min\left\{\frac{1}{2}\min\{R-|\phi_{0,k}(z)|||\phi_{0,k}(z)|<R\},\right.
$$
$$
\left.\frac{1}{2}\min\{|\phi_{0,k_1}(z)-\phi_{0,k_2}(z)|||\phi_{0,k_1}(z)|,|\phi_{0,k_2}(z)|<R,
\phi_{0,k_1}(z)\ne\phi_{0,k_2}(z)\}\right\}.
$$
In the sequel we will be interested in such configurations determined by $(z,R,\epsilon)$ for which $R\rightarrow+\infty$ 
and $\epsilon\rightarrow 0^+$.

\section{Amenability of the automorphic group}\label{sec34}

Let us assume that the sequence $\{f_n\}_n\subseteq E$ satisfies the following: \\
(a) $f_n\rightarrow f\in E$ uniformly on compact subsets of $\mathbb{C}$. \\
(b) The discrete groups ${\rm Aut}(f_n)$ are amenable. \\
\\
{\bf Example:} For $f(z)=\sum_{j=0}^{\infty}a_jz^j\in E$ we can take $f_n(z)=\sum_{j=0}^na_jz^j\in E$.
Then $f_n(z)\in\mathbb{C}[z]$, polynomials, and hence for each $n$ ${\rm Aut}(f_n)$ is a finite group
(of order $\deg f_n$). Hence ${\rm Aut}(f_n)$ are amenable $\forall\,n\in\mathbb{Z}^+$, for which $f_n\in E$. \\
\\
One might try to construct a F\o lner sequence in order to prove amenability of ${\rm Aut}(f),\,f\in E$.
Let us recall few notions and results.

\begin{definition}\label{def10}
A discrete group $G$ satisfies the F\o lner condition if for every finite subset $A\subseteq G$ and every
$\epsilon>0$ there exists a finite nonempty subset $F\subseteq G$ such that $\forall\,a\in A$ we have:
$$
\frac{|aF\vartriangle F|}{|F|}\le\epsilon.
$$
If $G$ is locally compact we use the same definition but $A$ is a compact subgroup, $F$ is a Borel set
with positive Haar measure and we use Haar measure instead of cardinality.
\end{definition}
\noindent
{\bf Example:} All finite (or compact in the locally compact case) groups satisfy the F\o lner condition, by
simply taking $F=G$ ($aF\vartriangle F=aG\vartriangle G=\emptyset$). \\

\begin{definition}\label{def11}
For a discrete and countable (resp. locally compact) group $G$, a F\o lner sequence is a sequence $\{F_n\}$
of nonempty finite (resp. compact) subsets of $G$ such that:
$$
\frac{|gF_n\vartriangle F_n|}{|F_n|}\rightarrow_{n\rightarrow\infty} 0\,\,\,\,\,\left({\rm resp.}\,\,
\frac{\mu(gF_n\vartriangle F_n)}{\mu(F_n)}\rightarrow_{n\rightarrow\infty} 0\right)\,\,\forall\,g\in G.
$$
\end{definition}
\noindent
The following lemma is well-known.

\begin{lemma}\label{lem4} {\rm\bf (\cite{pat})}
A group satisfies the F\o lner condition, if and only if it has a F\o lner sequence.
\end{lemma}
\noindent
{\bf Example:} The group $\mathbb{Z}$ has a F\o lner sequence, namely $F_n=\{-n,\ldots,n\}$. \\
\\
The usefulness of Definition \ref{def10} comes from the following well-known theorem.

\begin{theorem}\label{thm25} {\rm\bf (\cite{pat})}
A group satisfies the F\o lner condition, if and only if it is amenable.
\end{theorem}
\noindent
Coming back to our setting where $f\in E$, $f_n\in E$ are polynomials, we fix $z\in\mathbb{C}$. We take a
sequence $0<R_1<R_2<\ldots<R_n<\ldots (R_n\rightarrow\infty)$, and for each pair $(z,R_n)$ we take an
$\epsilon_n$ so that $0<\epsilon_n<\delta(z,R_n)$ and $\epsilon_n\rightarrow 0^+$. We take $f_{m(n)}$ such
that, using the notations of Theorem \ref{thm24}, $m(n)>N(R_n,\epsilon_n)$. We define a sequence of finite
subsets of ${\rm Aut}(f)$ by:
$$
F_n=\{\phi_{0k}|\,|\phi_{0k}(z)|<R_n\},\,\,\,\,n\in\mathbb{Z}^+.
$$
We fix an automorphic function $\phi_{0l}\in{\rm Aut}(f)$ and we consider:
$$
\frac{|\phi_{0l}\circ F_n\vartriangle F_n|}{|F_n|}.
$$
By the choice $m(n)>N(R_n,\epsilon_n)$ there is (as explained after Theorem \ref{thm24}) a canonical bijection
between $F_n$ and $F_n(f_{m(n)})=\{\phi_{0k}^{[m(n)]}|\,|\phi_{0k}^{[m(n)]}(z)|<R_n\}$. Moreover, if $n$ is
large enough, then $\phi_{0l}\in F_n$ and so it is canonically corresponding to $\phi_{0l}^{[m(n)]}$. Hence:
$$
\frac{|\phi_{0l}\circ F_n\vartriangle F_n|}{|F_n|}=
\frac{|\phi_{0l}^{[m(n)]}\circ F_n(f_{m(n)})\vartriangle F_n(f_{m(n)})|}{|F_n(f_{m(n)})|}.
$$
By $\phi_{0l}\in F_n$ we clearly have $\phi_{0l}^{[m(n)]}\in F_n(f_{m(n)})$ and in fact when $n\rightarrow\infty$,
we have $\phi_{0l}^{[m(n)]}(z)\rightarrow\phi_{0l}(z)$. Thus $|\phi_{0l}^{[m(n)]}(z)|$ is bounded for $n\rightarrow\infty$
and gets closer as we please to $|\phi_{0l}(z)|$.

\begin{theorem}\label{thm26}
If
$$
\lim_{n\rightarrow\infty}\frac{|F_n(f_{m(n)})|}{|{\rm Aut}(f_{m(n)})|}=1,
$$
then $\{F_n\}$ is a F\o lner sequence and hence ${\rm Aut}(f)$ is amenable.
\end{theorem}
\noindent
{\bf Proof.} \\
Clearly $\forall\,n\in\mathbb{Z}^+$ we have $|F_n(f_{m(n)})|\le|{\rm Aut}(f_{m(n)})|$ simply because
$F_n(f_{m(n)})\subseteq{\rm Aut}(f_{m(n)})$. Let us denote $|F_n(f_{m(n)})|=(1-\epsilon_n)|{\rm Aut}(f_{m(n)})|$.
Then $0\le\epsilon_n\le 1$, and by our assumption:
$$
1=\lim_{n\rightarrow\infty}\frac{|F_n(f_{m(n)})|}{|{\rm Aut}(f_{m(n)})|}=\lim_{n\rightarrow\infty}(1-\epsilon_n).
$$
Thus $\lim_{n\rightarrow\infty}\epsilon_n=0$. Clearly, we have the following straight forward estimate:
$$
0\le\frac{\phi_{0l}^{[m(n)]}\circ F_n(f_{m(n)})\vartriangle F_n(f_{m(n)})|}{|F_n(f_{m(n)})|}\le
\frac{2\epsilon_n}{1-\epsilon_n}.
$$
Hence:
$$
0\le\frac{|\phi_{0l}\circ F_n\vartriangle F_n|}{|F_n|}\le\frac{2\epsilon_n}{1-\epsilon_n}.
$$
This implies that:
$$
\lim_{n\rightarrow\infty}\frac{|\phi_{0l}\circ F_n\vartriangle F_n|}{|F_n|}=0,
$$
and our theorem follows. $\qed$ \\
\\
We can give another condition on $f(w)$ that implies that ${\rm Aut}(f)$ is amenable. This time it is a geometrical
condition. We start with the following:

\begin{definition}\label{def12}
Let $f\in E$. Suppose that the $z$-plane is tiled up by a system of fundamental domains $\{\Omega_j\}_j$ of the
entire function $w=f(z)$. We say that two fundamental domains $\Omega_1$ and $\Omega_2$ are neighboring if 
$\Omega_1\cap\Omega_2=\emptyset$, $\partial\Omega_1\cap\partial\Omega_2\ne\emptyset$.
\end{definition}

\begin{definition}\label{def13}
Let $f\in E$. Suppose that the $z$-plane is tiled up by a system of fundamental domains $\{\Omega_j\}_j$ of the
entire function $w=f(z)$. Let $\Omega_0$ be one of the fundamental domains in the system and let us denote by
$G_1(\Omega_0)$ the family of all the neighboring domains of $\Omega_0$. We will sometimes denote the members
of $G_1(\Omega_0)=\{\Omega_{1j}\}_j$ and call $G_1(\Omega_0)$ the first generation about $\Omega_0$.
\end{definition}

\begin{definition}\label{def14}
Let $f\in E$. Suppose that the $z$-plane is tiled up by a system of fundamental domains $\{\Omega_j\}_j$ of the
entire function $w=f(z)$. Let $\Omega_0$ be one of the fundamental domains in the system. Let $n\in\mathbb{Z}_{\ge 2}$.
The $n$'th generation about $\Omega_0$ is denoted by $G_n(\Omega_0)=\{\Omega_{nj}\}_j$ and is defined recursively
by the following recursive equation:
$$
G_n(\Omega_0)=\bigcup_{\Omega\in G_{n-1}(\Omega_0)}G_1(\Omega)-\{\Omega_0\}\cup\bigcup_{j=1}^{n-1}G_j(\Omega_0).
$$
The counting function of the generations about $\Omega_0$ is defined by: $g(\Omega_0,n)=|G_n(\Omega_0)|$.
\end{definition}
\noindent
{\bf Examples:} \\
1) Let $f(z)=z^N$ for some $N\in\mathbb{Z}_{\ge 2}$. Then a natural system of fundamental domains are:
$$
\Omega_j=\left\{z\in\mathbb{C}|\,\frac{2\pi j}{N}<\arg z<\frac{2\pi(j+1)}{N}\right\},\,\,\,j=0,1,\ldots,N-1.
$$
Then $\forall\,j$ $G_1(\Omega_j)=\{\Omega_0,\Omega_1,\ldots,\Omega_{N-1}\}-\{\Omega_j\}$, and $G_n(\Omega_j)=\emptyset$,
$\forall\,n>1$. So:
$$
g(\Omega_j,n)=\left\{\begin{array}{lll} N-1 & {\rm if} & n=1 \\ 0 & {\rm if} & n>1\end{array}\right..
$$
2) Let $f(z)=e^z$. A natural system of fundamental domains are:
$$
\Omega_j=\left\{z\in\mathbb{C}|\,2\pi j<\Im z<2\pi(j+1)\right\},\,\,\,j\in\mathbb{Z}.
$$
Here we have: $G_n(\Omega_0)=\{\Omega_{-n},\Omega_n\}$, and hence $g(\Omega_0,n)=2$.

\begin{theorem}\label{thm27}
Let $f\in E$ have a system of fundamental domains. Let $\Omega_0$ be a fundamental domain in the system and
let $G_1(\Omega_0)=\{\Omega_{1j}\}_j$. Let denote by (as usual) by $\phi_{0j}:\Omega_0\rightarrow\Omega_{1j}$
the corresponding automorphic function of $f$. Then $\{\phi_{0j}\}_j$ is a generating set of the automorphic
group, ${\rm Aut}(f)$.
\end{theorem}
\noindent
{\bf Proof.} \\
This is immediate from the definitions. The automorphic function $\phi_{12}:\Omega_{11}\rightarrow\Omega_{12}$
is given by the composition: $\phi_{02}\circ\phi_{01}^{-1}$ which maps as follows:
$$
\Omega_{11}\overset{\phi_{01}^{-1}}{\rightarrow}\Omega_0\overset{\phi_{02}}{\rightarrow}\Omega_{12}.
$$
If, for instance, the curve $a\curvearrowright b$ is common to $\partial\Omega_{12}$ and to $\partial\Omega_{24}$
and the automorphic function $\phi_{02}:\,\Omega_0\rightarrow\Omega_{12}$ maps the curve $a'\curvearrowright b'$
which is common to $\partial\Omega_0$ and t $\partial\Omega_{13}$ to the curve $a\curvearrowright b$, then
$\phi_{0(24)}:\Omega_0\rightarrow\Omega_{24}$ is given by the composition: $\phi_{0(24)}=\phi_{03}\circ\phi_{02}$, etc... $\qed $ \\
\\
In particular we have:

\begin{corollary}\label{cor12}
Let $f\in E$, have a system of fundamental domains. If ${\rm Aut}(f)$ is not a finitely generated group, then for
any system $\{\Omega_j\}_j$ of fundamental domains and for any $j$, the first generation $G_1(\Omega_j)$ is an infinite
family.
\end{corollary}

\begin{remark}\label{rem36}
We recall that according to Shimizu's definition in \cite{s1}, the boundaries of a fundamental system of an
entire function have no accumulation point in the finite plane. Moreover, not every entire function has a system
of fundamental domains. Gross constructed an entire function which has all the points of $\mathbb{C}$ as its
asymptotic values. In \cite{s1} Shimizu proved that the Gross function has no system of fundamental domains.
\end{remark}
\begin{theorem}\label{thm28}
Let $f\in E$ have a system of fundamental domains $\{\Omega_j\}_j$ having the property that $\forall\,j$ we have:
$$
\lim_{n\rightarrow\infty}\left\{\frac{g(\Omega_j,n)}{\sum_{m=1}^{n}g(\Omega_j,m)}\right\}=0.
$$
In particular the $g(\Omega_j,m)$ are always finite! Then ${\rm Aut}(f)$ is amenable.
\end{theorem}
\noindent
{\bf Proof.} \\
One can check that the finite sets of automorphic functions $\phi:\,\Omega_0\rightarrow\Omega$ where
$\Omega\in\bigcup_{k=1}^n G_k(\Omega_0)$, which we denote by $F_n$ form a F\o lner sequence for ${\rm Aut}(f)$. $\qed $

\begin{corollary}\label{cor13}
Let $f\in E$ have a system of fundamental domains $\{\Omega_j\}_j$ such that $\forall\,j$ there is a polynomial $P_j(x)$
of degree $d_j$ for which $g(\Omega_j,n)\in\Omega(P_j(n))$, i.e. there are two positive numbers $0<c_j<C_j$ such
that $\forall\,n\in\mathbb{Z}^+$, $c_j\cdot P_j(n)\le g(\Omega_j,n)\le C_j\cdot P_j(n)$. Then ${\rm Aut}(f)$ is amenable.
\end{corollary}
\noindent
{\bf Proof.} \\
We will use the following well-known estimate of the moments of the natural numbers:
$$
1^d+2^d+\ldots+n^d=\frac{n^{d+1}}{d+1}+\frac{n^d}{2}+\frac{dn^{d-1}}{12}+\mathcal{O}(n^{d-3}).
$$
By this estimate we obtain:
$$
\lim_{n\rightarrow\infty}\frac{n^d}{1^d+2^d+3^d+\ldots+n^d}=0.
$$
By the assumption on the counting function $g(\Omega_j,n)$ and by Theorem \ref{thm28} the result follows. $\qed $ \\

\begin{remark}\label{rem37}
Theorem \ref{thm28} does not imply anything in the case where $g(\Omega_j,n)=\Omega(q^n)$, i.e. a geometric growth.
For:
$$
\lim_{n\rightarrow\infty}\frac{q^n}{1+q+q^2+\ldots+q^n}=1.
$$
\end{remark}

\noindent
{\it Ronen Peretz \\
Department of Mathematics \\ Ben Gurion University of the Negev \\
Beer-Sheva , 84105 \\ Israel \\ E-mail: ronenp@math.bgu.ac.il} \\ 
 
\end{document}